\pdfoutput=1
\documentclass[a4paper, 11pt, abstract=true, titlepage, oneside]{scrartcl}
\makeatletter

\usepackage[T1]{fontenc}
\usepackage[utf8]{inputenc}
\usepackage[onehalfspacing]{setspace}
\usepackage[headsepline]{scrlayer-scrpage}
\usepackage{layout}
\usepackage{geometry}
\usepackage{adjustbox}
\usepackage{authblk}
\usepackage[titletoc,title]{appendix}
\usepackage[hyphens]{url}
\usepackage{color}

\usepackage{amsmath,amssymb,amsfonts,amsthm, mathtools, bm}
\usepackage{mathrsfs}
\usepackage{accents}
\usepackage{thmtools, thm-restate}
\usepackage{nicefrac}
\usepackage{array}

\usepackage[normal]{threeparttable}
\usepackage{threeparttablex}
\usepackage{tabularx}
\usepackage{longtable}
\usepackage{booktabs}
\usepackage{colortbl}
\usepackage{multirow}
\usepackage{makecell}
\usepackage{algorithm,algpseudocode}
\usepackage[bottom,hang,multiple]{footmisc}
\usepackage[acronym]{glossaries}

\usepackage{subcaption}
\usepackage{import}
\usepackage{soul}
\usepackage{graphicx}
\usepackage{xcolor} % For \textcolor command
\usepackage{placeins}

\usepackage{tikz}
\usetikzlibrary{arrows.meta}
\usetikzlibrary{math}
\usetikzlibrary{shapes}
\usetikzlibrary{patterns,shapes.arrows}
\usepackage{pgfplots}
\pgfplotsset{compat=1.18} % Adjust this to the version you're targeting
\usepgfplotslibrary{groupplots}
\usepgfplotslibrary{groupplots,dateplot}
\pgfplotsset{compat=newest}
\usepgfplotslibrary{dateplot}
\usetikzlibrary{backgrounds}
\pgfdeclarelayer{background}
\pgfdeclarelayer{foreground}
\pgfsetlayers{background,main,foreground} 
\usepackage{tikzscale}

\usepackage{enumerate}
\usepackage{multicol}
\usepackage[author=Patrick]{pdfcomment}
\usepackage{xparse}
\usepackage{twoopt}
\usepackage{etoolbox}

\usepackage{eurosym}
\usepackage{ellipsis}
\usepackage{natbib}
\usepackage{paralist}
\usepackage{ulem}
\AtBeginDocument{
	\newgeometry{
		textheight=9in,
		textwidth=6.5in,
		top=1in,
		headheight=16.49675pt,
		headsep=25pt,
		footskip=30pt
	}
}
\newsavebox{\abstractbox}
\renewenvironment{abstract}
{\begin{lrbox}{0}\begin{minipage}{\textwidth}
			\begin{center}\normalfont\sectfont\abstractname\end{center}\quotation}
		{\endquotation\end{minipage}\end{lrbox}%
	\global\setbox\abstractbox=\box0 }

\makeatletter
\expandafter\patchcmd\csname\string\maketitle\endcsname
{\vskip\z@\@plus3fill}
{\vskip\z@\@plus2fill\box\abstractbox\vskip\z@\@plus1fill}
{}{}
\makeatother

\allowdisplaybreaks

\addtokomafont{captionlabel}{\footnotesize\bfseries} % style for caption labels
\addtokomafont{caption}{\footnotesize\bfseries} % and the caption
\bibpunct[, ]{(}{)}{,}{a}{}{,}%
\newtheoremstyle{upright} % name
  {3pt}                    % Space above
  {3pt}                    % Space below
  {\upshape}               % Body font
  {}                       % Indent amount
  {\bfseries}              % Theorem head font
  {.}                      % Punctuation after theorem head
  {.5em}                   % Space after theorem head
  {}                       % Theorem head spec

\theoremstyle{upright}
\newtheorem{remark}{Remark}
\newtheorem{result}{Result}
\newtheorem{assumption}{Assumption}
\newtheorem{definition}{Definition}
\newtheorem{theorem}{Theorem}

\def\Halmos{\mbox{\quad$\square$}}% \square <--> 
\let\proof\relax 
\let\endproof\relax
\let\theoremstyle\relax

\counterwithin*{equation}{section}

\newenvironment{myfont}{\fontfamily{ccr}\selectfont}{\par}
\DeclareTextFontCommand{\textmyfont}{\myfont}

\AtBeginDocument{%
	\abovedisplayskip=7.5pt plus 0pt minus 0pt
	\abovedisplayshortskip=7.5pt plus 0pt minus 0pt
	\belowdisplayskip=7.5pt plus 0pt minus 0pt
	\belowdisplayshortskip=7.5pt plus 0pt minus 0pt
}
\newcolumntype{L}[1]{>{\raggedright\let\newline\\\arraybackslash\hspace{0pt}}p{#1}}
\newcolumntype{C}[1]{>{\centering\let\newline\\\arraybackslash\hspace{0pt}}p{#1}}
\newcolumntype{R}[1]{>{\raggedleft\let\newline\\\arraybackslash\hspace{0pt}}p{#1}}
\renewcommand{\emph}[1]{\textit{#1}}

\algtext*{EndIf} % Suppress "end if" text
\algtext*{EndWhile} % Suppress "end while" text
\algtext*{EndFor} % Suppress "end while" text

\begin{document}
\emergencystretch 3em
% first put your acronyms and abbrevs.
% a
\newcommand{\Arc}[0]{a}
\newcommand{\LastArc}{\bar{\Arc}^\Trip}
\newcommand{\FirstArc}{\ubar{\Arc}^\Trip}
\newacronym{amod}{AMoD}{autonomous mobility-on-demand}
% b
\newacronym{bpr}{BPR}{Bureau of Public Roads}
\newacronym{beb}{B\&B}{branch and bound}
% c
\newcommand{\ArcCapacity}[1]{c_{#1}}
% d
% e
\newcommand{\TripReleaseTime}[1]{e^{#1}}
\newcommand{\TripEarliestDepartureTime}[2]{e_{#1}^{#2}}
\newcommand{\EarliestEntry}[2]{\ubar{e}_{#1}^{#2}}
\newcommand{\EarliestEntryVector}[0]{\ubar{\boldsymbol{e}}}
\newcommand{\ReleaseTimesVector}[0]{{\boldsymbol{e}}}
\newcommand{\LatestEntry}[2]{\bar{e}_{#1}^{\hspace{0.5pt}#2}}
% f
\newcommand{\Flow}[2]{f_{#1}^{#2}}
\newcommand{\MinFlow}[2]{\ubar{f}_{#1}^{#2}}
\newcommand{\MaxFlow}[2]{\bar{f}_{#1}^{#2}}
\newcommand{\ThFlow}[2]{\hat{f}_{#1}^{#2}}
% g
% h
\newacronym{hc}{HC}{high-congestion}
% i
% j
% k
% l
\newacronym{lc}{LC}{low-congestion}
% ell (small calligraph l)
\newcommand{\TripLatestExitTime}[2]{\ell_{#1}^{#2}}
\newcommand{\TripDeadline}[1]{\ell^{#1}}
\newcommand{\EarliestExit}[2]{\ubar{\ell}_{#1}^{#2}}
\newcommand{\LatestExit}[2]{\bar{\ell}_{#1}^{#2}}
\newcommand{\LatestExitVector}[0]{\bar{\boldsymbol{\ell}}}
% m
\newacronym{milp}{MILP}{mixed integer linear program}
\newacronym{mpc}{MPC}{model predictive control}
\newacronym{mps}{MPSs}{mobility service providers}
% n
\newacronym{nyc}{NYC}{New York City}
% o
% p
\newcommand{\TripPath}[1]{p^{#1}}
\newcommand{\TripPathBar}[1]{\overline{\TripPath{#1}}}
\newcommand{\TripPathStartingFromArc}[2]{\TripPath{#2}(#1)}
% q
% r
\newcommand{\Trip}[0]{r}
\newcommand{\SecondTrip}[0]{\Trip'}
% s
\newcommand{\TripDepartureOnArc}[2]{s_{#1}^{#2}}
% t
\newcommand{\Arrival}[0]{t}
% u
% v
\newcommand{\Node}[0]{v}
% w
% x
% y
% z

%%%%%%%%%%%%%%%%%%%%%%%%%%%%%%%%%%%%%%%%%%%%%%%%%%%%%%%%%%%%%%%%%%%%%%%%%%%%%%%%%%%%%%%%%%%%%%%%%%%%%%%%%%%%%%%%%%%%%%%%%
% Commands
% \newcommand{}{}
%-Large-Arabic-------------------------------------------------------------------------------
% A
%\newcommand{}{A}
% B
%\newcommand{}{B}
% C
\newcommand{\Counter}[0]{\textsc{C}}
%\newcommand{}{C}
% D
%\newcommand{}{D}
% E
%\newcommand{}{E}
% F
%\newcommand{}{F}
% G
%\newcommand{}{G}
% H
\newcommand{\HCOfflineMedianReduction}[0]{66}
\newcommand{\HCOfflineMeanReduction}[0]{66}
\newcommand{\HCOnlineMeanReduction}[0]{56}
\newcommand{\HCOfflineMinReduction}[0]{50}
\newcommand{\HCOfflineMaxReduction}[0]{80}
\newcommand{\HCMaximumStaggeringApplied}[0]{\text{six}}
%\newcommand{}{H}
% I
%\newcommand{}{I}
% J
\newacronym{j3uptnwt}{J3-UPT-NWT}{three machines job shop scheduling problem with unit processing times and no-wait constraints}
\newcommand{\JSUTNW}[0]{J_{3}|p_{i,j}= 1,\text{nwt}|C_{\text{max}}}
%\newcommand{}{J}
% K
\newcommand{\NumPieces}{K}
% L
\newcommand{\LowerBound}[0]{\text{LB}}
%---
\newcommand{\LCOfflineMedianReduction}[0]{99}
\newcommand{\LCOfflineMeanReduction}[0]{94}
\newcommand{\LCOfflineMinReduction}[0]{80}
\newcommand{\LCOnlineMeanReduction}[0]{90}
\newcommand{\LCMaximumStaggeringApplied}[0]{\text{two}}
%---
%\newcommand{}{L}
% M
\newcommand{\BigM}[1]{M_{#1}}
% N
\newcommand{\NumberOfEpochs}[0]{N}
% O
%\newcommand{}{O}
% P
%\newcommand{}{P}
% Q
\newcommand{\Queue}[0]{\textsc{Q}}
%\newcommand{}{Q}
% R
%\newcommand{}{R}
% S
%\newcommand{}{S}
% T
%\newcommand{}{T}
% U
%\newcommand{}{U}
% V
%\newcommand{}{V}
% W
%\newcommand{}{W}
% X
%\newcommand{}{X}
% Y
%\newcommand{}{Y}
% Z
\newcommand{\TotalDelayUncSolution}[0]{\bar{\text{Z}}}
\newcommand{\TotalDelay}[1]{\text{Z}(#1)}
%-Small-Arabic-------------------------------------------------------------------------------
% a
%\newcommand{}{a}
% b
%\newcommand{}{b}
% c
%\newcommand{}{c}
% d
%\newcommand{}{d}
% e
%\newcommand{}{e}
% f
%\newcommand{}{f}
% g
%\newcommand{}{g}
% h
%\newcommand{}{h}
% i
%\newcommand{}{i}
% j
%\newcommand{}{j}
% k
%\newcommand{}{k}
% l
%\newcommand{}{l}
% m
%\newcommand{}{m}
% n
%\newcommand{}{n}
% o
%\newcommand{}{o}
% p
%\newcommand{}{p}
% q
%\newcommand{}{q}
% r
%\newcommand{}{r}
% s
%\newcommand{}{s}
% t
%\newcommand{}{t}
% u
%\newcommand{}{u}
% v
%\newcommand{}{v}
% w
%\newcommand{}{w}
% x
\newcommand{\DelayOnArc}[2]{x_{#1}^{#2}}
\newcommand{\MinDelay}[2]{\ubar{x}_{#1}^{\hspace{0.7pt}#2}}
\newcommand{\MinDelayVector}[0]{\ubar{\boldsymbol{x}}}
\newcommand{\MaxPotentialDelay}[2]{\check{x}_{#1}^{\hspace{1pt}#2}}
%\newcommand{}{x}
% y
%\newcommand{}{y}
% z
%\newcommand{}{z}
%-Large-Calligraphic-------------------------------------------------------------------------
% A
\newcommand{\SetArcs}[0]{\mathcal{A}}
\newcommand{\SetConflictingArcs}[0]{\hat{\SetArcs}}
% B
%\newcommand{}{\mathcal{B}}
% C
%\newcommand{}{\mathcal{C}}
% D
%\newcommand{}{\mathcal{D}}
% E
%\newcommand{}{\mathcal{E}}
% F
%\newcommand{}{\mathcal{F}}
% G
\newcommand{\Graph}[0]{\mathcal{G}}
%\newcommand{}{\mathcal{G}}
% H
%\newcommand{}{\mathcal{H}}
% I
%\newcommand{}{\mathcal{I}}
% J
%\newcommand{}{\mathcal{J}}
% K
%\newcommand{}{\mathcal{K}}
% L
%\newcommand{}{\mathcal{L}}
% M
%\newcommand{}{\mathcal{M}}
% N
%\newcommand{}{\mathcal{N}}
% O
\newcommand{\BigO}[1]{\mathcal{O}\left( #1 \right)}
% P
\newcommand{\SetOfPairs}[1]{\mathcal{P}_{#1}}
\newcommand{\SetOfPartition}[2]{\mathcal{P}_{#1}^{#2}}
%\newcommand{}{\mathcal{P}}
% Q
%\newcommand{}{\mathcal{Q}}
% R
\newcommand{\SetTrips}[0]{\mathcal{R}}
\newcommand{\SetTripsOnArc}[1]{\SetTrips_{#1}}
\newcommand{\SetTripsEpoch}[1]{\SetTrips^{#1}}
\newcommand{\TransferTrips}[0]{\SetTrips^{T}}
\newcommand{\DummyTrips}[0]{\SetTrips^{D}}
% S
\newcommand{\ConflictingSet}[2]{\mathcal{S}_{#1}^{#2}}
%\newcommand{}{\mathcal{S}}
% T
\newcommand{\SetOfArcArrivals}[1]{\mathcal{T}_{#1}}
%\newcommand{}{\mathcal{T}}
% U
%\newcommand{}{\mathcal{U}}
% V
\newcommand{\SetNodes}[0]{\mathcal{V}}
% W
%\newcommand{}{\mathcal{W}}
% X
% Y
%\newcommand{}{\mathcal{Y}}
% Z
%\newcommand{}{\mathcal{Z}}
%-Large-Scr-------------------------------------------------------------------------
% A
\newcommand{\PolyAlgo}{\mathscr{A}}
% B
%\newcommand{}{\mathscr{B}}
% C
%\newcommand{}{\mathscr{C}}
% D
%\newcommand{}{\mathscr{D}}
% E
%\newcommand{}{\mathscr{E}}
% F
%\newcommand{}{\mathscr{F}}
% G
%\newcommand{}{\mathscr{G}}
% H
%\newcommand{}{\mathscr{H}}
% I
%\newcommand{}{\mathscr{I}}
% J
%\newcommand{}{\mathscr{J}}
% K
%\newcommand{}{\mathscr{K}}
% L
%\newcommand{}{\mathscr{L}}
% M
%\newcommand{}{\mathscr{M}}
% N
%\newcommand{}{\mathscr{N}}
% O
%\newcommand{}{\mathscr{O}}
% P
%\newcommand{}{\mathscr{P}}
% Q
%\newcommand{}{\mathscr{Q}}
% R
%\newcommand{}{\mathscr{R}}
% S
%\newcommand{}{\mathscr{S}}
% T
%\newcommand{}{\mathscr{T}}
% U
%\newcommand{}{\mathscr{U}}
% V
%\newcommand{}{\mathscr{V}}
% W
%\newcommand{}{\mathscr{W}}
% X
%\newcommand{}{\mathscr{X}}
% Y
%\newcommand{}{\mathscr{Y}}
% Z
%\newcommand{}{\mathscr{Z}}
%-Large-Greek--------------------------------------------------------------------------------
% Alpha
%\newcommand{}{A}
% Beta
%\newcommand{}{B}
% Gamma
%\newcommand{}{\Gamma}\\
% Delta
\newcommand{\Overlap}[0]{\Delta}
\newcommand{\EpochLength}[0]{\Delta t}
%\newcommand{}{\Delta}
% Epsilon
%\newcommand{}{E}
% Zeta
%\newcommand{}{Z}
% Eta
%\newcommand{}{H}
% Theta
%\newcommand{}{\Theta}
% Iota
%\newcommand{}{I}
% Kappa
%\newcommand{}{K}
% Lambda
%\newcommand{}{\Lambda}
% Mu
%\newcommand{}{M}
% Nu
%\newcommand{}{N}
% Xi
%\newcommand{}{\Xi}
% Omicron
%\newcommand{}{O}
% Pi
\newcommand{\FeasibleSet}[0]{\Pi}
%\newcommand{}{\Pi}
% Rho
%\newcommand{}{P}
% Sigma
%\newcommand{}{\Sigma}
% Tau
%\newcommand{}{T}
% Upsilon
\newcommand{\ConflictsList}[0]{\Upsilon}
%\newcommand{}{\Upsilon}
% Phi
%\newcommand{}{\Phi}
% Chi
%\newcommand{}{X}
% Psi
%\newcommand{}{\Psi}
% Omega
%\newcommand{}{\Omega}
%-Small-Greek--------------------------------------------------------------------------------
% alpha
\newcommand{\VarAlpha}[3]{\alpha_{#1}^{#2,#3}}
% beta
\newcommand{\VarBeta}[3]{\beta_{#1}^{#2,#3}}
% gamma
\newcommand{\VarGamma}[3]{\gamma_{#1}^{#2,#3}}
% delta
\newcommand{\Successor}[2]{\delta^{#2}({#1})}
% epsilon
\newcommand{\SmallConstant}[0]{\varepsilon}
% zeta
%\newcommand{}{\zeta}
% eta
\newcommand{\TimeLimitMILP}[0]{\eta}
% theta
%\newcommand{}{\vartheta}
% iota
%\newcommand{}{\iota}
% kappa
%\newcommand{}{\kappa}
% lambda

% mu
\newcommand{\PWLSlope}[2]{\mu_{#1}^{#2}}
%\newcommand{}{\mu}
% nu
%\newcommand{}{\nu}
% xi
%\newcommand{}{\xi}
% omicron
%\newcommand{}{o}
% pi
\newcommand{\Solution}[0]{\pi}
% rho
\newcommand{\ForwardShift}[1]{\rho^{\text{#1}}}
% sigma
\newcommand{\StaggeringApplied}[1]{\sigma^{#1}}
\newcommand{\MaxStaggering}[1]{\bar{\sigma}^{#1}}
\newcommand{\TimePercentage}[0]{\varsigma^{\text{MAX}}}
% tau
\newcommand{\ArcNominalTravelTime}[1]{\tau_{#1}}
\newcommand{\ArcTravelTime}[2]{\tau^{#2}_{#1}(\Flow{#1}{#2})}
%\newcommand{}{\tau}
% upsilon
%\newcommand{}{\upsilon}
% phi
\newcommand{\BackwardShift}[1]{\phi^{\text{#1}}}
% chi
%\newcommand{}{\chi}
% psi
%\newcommand{}{\psi}
% omega
%\newcommand{}{\omega}
%-Large-Blackboard-Bold--------------------------------------------------------------------
% A
%\newcommand{}{\mathbb{A}}
% B
%\newcommand{}{\mathbb{B}}
% C
%\newcommand{}{\mathbb{C}}
% D
%\newcommand{}{\mathbb{D}}
% E
%\newcommand{}{\mathbb{E}}
% F
%\newcommand{}{\mathbb{F}}
% G
%\newcommand{}{\mathbb{G}}
% H
%\newcommand{}{\mathbb{H}}
% I
%\newcommand{}{\mathbb{I}}
% J
%\newcommand{}{\mathbb{J}}
% K
%\newcommand{}{\mathbb{K}}
% L
%\newcommand{}{\mathbb{L}}
% M
%\newcommand{}{\mathbb{M}}
% N
\newcommand{\NaturalNumbers}[0]{\mathbb{N}}
%\newcommand{}{\mathbb{N}}
% O
%\newcommand{}{\mathbb{O}}
% P
%\newcommand{}{\mathbb{P}}
% Q
%\newcommand{}{\mathbb{Q}}
% R
\newcommand{\PositiveRealNumbers}[0]{\mathbb{R}^{\geq 0}}
%\newcommand{}{\mathbb{R}}
% S
%\newcommand{}{\mathbb{S}}
% T
%\newcommand{}{\mathbb{T}}
% U
%\newcommand{}{\mathbb{U}}
% V
%\newcommand{}{\mathbb{V}}
% W
%\newcommand{}{\mathbb{W}}
% X
%\newcommand{}{\mathbb{X}}
% Y
%\newcommand{}{\mathbb{Y}}
% Z
%\newcommand{}{\mathbb{Z}}
%-Other-------------------------------------------------------------------
% definition of the under bar symbol
\newcommand{\ubar}[1]{\underaccent{\bar}{#1}}
% Figure parameters
\newcommand{\GapPlotsHeight}{3.5cm}
\newcommand{\TotalDelayBarplotWidth}{0.5\columnwidth}
\newcommand{\ShiftAppliedWidth}{0.5\columnwidth}
\newcommand{\DelayReductionWidth}{0.45\columnwidth}
\newcommand{\DelayReductionHeight}{4cm}

\title{\large Staggered Routing in Autonomous
Mobility-on-Demand Systems}

% and the authors
\author[1]{\normalsize Antonio Coppola}
\author[1]{\normalsize Gerhard Hiermann}
\author[2]{\normalsize Dario Paccagnan}
\author[3]{\normalsize Maximilian Schiffer}
\affil{\small 
	School of Management, Technical University of Munich, Germany
	
	\scriptsize antonio.coppola@tum.de,
	\scriptsize gerhard.hiermann@tum.de
	
	\small
\textsuperscript{2}
        Department of Computing, Imperial College London, U.K.

        \scriptsize d.paccagnan@imperial.ac.uk

 \small
	\textsuperscript{3}School of Management \& Munich Data Science Institute,
	
	Technical University of Munich, Germany
	
	\scriptsize schiffer@tum.de}

% if you like - a date
\date{}

% in case you have a headline - otherwise outcomment
\lehead{\pagemark}
%\rehead{\normalfont\scriptsize\textbf{Schiffer et. al.:} \textit{Perspectives for electric commercial vehicles}}
%\lohead{\normalfont\scriptsize\textbf{Schiffer et. al.:} \textit{Perspectives for electric commercial vehicles}}
\rohead{\pagemark}

% finally your abstract
\begin{abstract}
\begin{singlespace}
{\small\noindent In autonomous mobility-on-demand systems, effectively managing vehicle flows to mitigate induced congestion and ensure efficient operations is imperative for system performance and positive customer experience. Against this background, we study the potential of staggered routing, i.e., purposely delaying trip departures from a system perspective, in order to reduce congestion and ensure efficient operations while still meeting customer time windows. We formalize the underlying planning problem and show how to efficiently model it as a mixed integer linear program. Moreover, we present a matheuristic that allows us to efficiently solve large-scale real-world instances both in an offline full-information setting and its online rolling horizon counterpart. We conduct a numerical study for Manhattan, New York City, focusing on low- and highly-congested scenarios. Our results show that in low-congestion scenarios, %i.e., scenarios in which only the AMoD fleet induces congestion but the system without the fleet would remain uncongested, 
staggering trip departures allows mitigating, on average, $\LCOfflineMeanReduction$\% of the induced congestion in a full information setting. In a rolling horizon setting, our algorithm allows us to reduce $\LCOnlineMeanReduction$\% of the induced congestion. 
In high-congestion scenarios, %i.e., scenarios in which the system already encounters exogenous congestion and the AMoD fleet induces further congestion, 
we observe an average reduction of $\HCOfflineMeanReduction$\% as the full information bound and an average reduction of $\HCOnlineMeanReduction$\% in our online setting. Surprisingly, we show that these reductions can be reached by shifting trip departures by a maximum of six minutes in both the low and high-congestion scenarios.\\

\smallskip}
{\footnotesize\noindent \textbf{Keywords:} autonomous mobility-on-demand; trip-staggering}
\end{singlespace}
\end{abstract}

% don't forget to make the tile
\maketitle
% and your chapters
\section{Introduction}\label{sec:introduction}

Congestion in cities is soaring. In 2022, the average US driver spent 51 hours stuck in traffic, resulting in a country-wide loss of \$81 billion when accounting solely for lost working hours~\citep{Pishue2023}. Beyond financial implications, congestion is responsible for emissions of noxious gases, including carbon dioxide, which fuel climate change and cause significant health issues \citep{LevyBuonocoreEtAl2010}. Against this backdrop, new mobility paradigms such as \gls{amod} hold the promise to mitigate the negative externalities of congestion by allowing for more efficient transport solutions. Such an \gls{amod} system consists of a centrally controlled fleet of self-driving vehicles that provide ride-hailing services to passengers \citep{Pavone2015}. By observing the system's complete state in real time and by leveraging centralized vehicle dispatching and rebalancing algorithms, a fleet operator may provide transportation solutions that are more efficient than passengers traveling either in their own cars or via uncoordinated ride-hailing services. However, no consensus exists yet as to whether \gls{amod} systems will ultimately reduce congestion. Indeed, while some experts argue for reduced congestion via improved fleet control, others claim that \gls{amod} systems may lead to an increase in demand, possibly resulting in higher transport volumes during peak {hours~\citep{OhSeshadriEtAl2020}.}  

To effectively reduce congestion, particularly in light of potentially increased demand, a number of approaches have been proposed, including ride-sharing~\citep{RuchLuEtAl2020}, allocation of parking space~\citep{ZhangGuhathakurta2017}, and strategic coordination with public transport \citep{SalazarLanzettiEtAl2019}, to name a few. Amongst them, \emph{congestion pricing} and \emph{fleet-optimal routing} have emerged as the most prominent ones both within academia and practice. 
To this end, congestion pricing assumes a system-wide authority that decides on a congestion fee for certain districts or roads to influence driver behavior and, consequently, the traffic flow \citep{Pigou1920}. Depending on the tolling scheme, individuals may be incentivized to choose an alternative route \citep{Roughgarden2005, PaccagnanGairing2021} or shift their travel to a different time \citep{ArnottdePalmaEtAl1990}. 
Accordingly, congestion pricing approaches aim to balance traffic flow across two dimensions, \emph{space} and \emph{time}. 
Fleet-optimal routing, instead, assumes control over a fleet of vehicles and aims at reducing congestion by coordinating their routing decisions \citep{JahnMoehringEtAl2005, Patriksson2015}. This approach can be applied locally by a fleet operator or via central control, i.e., traffic guidance by a municipality \citep{Levin2017,HoushmandWollensteinBetechEtAl2019}. Notably, a significant bulk of the existing literature in this domain has focused on balancing traffic flow across the \emph{space} dimension \citep{JalotaPaccagnanEtAl2023}. On the contrary, fleet-optimal approaches to balance traffic across the \emph{time} dimension remain unexplored. Nonetheless, with the advent of \gls{amod} systems, orchestrating vehicle flows over time provides an additional degree of freedom, which may prove valuable in taming congestion. %
Indeed, similar approaches have already provided factual benefits in related fields, e.g., in the context of ground delay programs for air traffic control \citep[cf.][]{Jacquillat2022}, thus further motivating the pursuit of this direction.

 Against these backdrops, we introduce the concept of \emph{staggered routing}, where an \gls{amod} operator delays the departure of trips over time to reduce local congestion bottlenecks so long as all trips' drop-off deadlines are met. Specifically, we take a first step to unravel the \emph{potential} of staggered routing in \gls{amod} systems. %Within this space, we hope our work will open new perspectives for future research. % instead of setting the final word on the problem.

\subsection{Contribution}
Our work aims to inform on the impact of staggering trip departures in reducing congestion. Towards this goal, we present three key contributions. First, we introduce the problem of \emph{staggered routing}, i.e., the problem of staggering the fleet departure times so as to minimize the resulting congestion while meeting drop-off deadlines. We then formulate this problem as a \gls{milp} and show that it is NP-hard. Second, we provide an effective algorithmic framework to solve the staggered routing problem at scale. Specifically, we i) show how a large number of the resulting big-M constraints can be dropped without altering the optimal solutions, and ii) develop a matheuristic that allows us to solve large-scale instances. This matheuristic bases on a construction heuristic, a problem-specific local search for intensification, while relying on our MILP to escape local optima and verify optimality. Finally, we apply our methodology to a real-world case study for the Manhattan area in New York City. Specifically, we study the impact that staggered routing has in both an offline and an online decision-making setting.

Our results show that in \gls{lc} scenarios, i.e., scenarios in which only the AMoD fleet induces congestion but the system without the fleet would remain uncongested, staggering trip departures allow to mitigate on average $\LCOfflineMeanReduction$\% of the induced congestion in a full-information setting. While not entirely reaching this upper bound, our rolling horizon approach still allows us to reduce $\LCOnlineMeanReduction$\% of the induced congestion in the online \gls{lc} scenario. In \gls{hc} scenarios, i.e., scenarios in which the system already encounters exogenous congestion and the AMoD fleet induces further congestion, we observe an average congestion reduction of $\HCOfflineMeanReduction$\% and $\HCOnlineMeanReduction$\% in the full-information and rolling-horizon settings. Surprisingly, we show that these reductions can be reached by shifting trip departures of a maximum of $\LCMaximumStaggeringApplied$ minutes in the LC scenario and by $\HCMaximumStaggeringApplied$ minutes in the \gls{hc} scenario while meeting each trip's arrival time.

\subsection{Technical challenges and our approach.}
\label{sec:our-approach}
Solving the staggered routing problem at scale requires tackling at least three intertwined challenges. 

First, following the commonly employed approach whereby one keeps track of the congestion level on each arc at predetermined time intervals (see, e.g., time-expanded construction in \citet{KoehlerLangkauEtAl2002}) is simply not possible as the size of the resulting optimization problem grows too quickly, making the solution of even moderately sized instances beyond reach. We tackle this challenge by taking a so-called Lagrangian perspective, that is, we ``follow'' each trip along their path and only keep track of the congestion this trip encounters at the time when it \emph{enters} each arc on its path. This allows us to formalize an optimization problem that is significantly more compact. 

Second, while the Lagrangian perspective we adopt allows to reduce the problem dimensionality, it makes it more difficult to determine the congestion level each trip encounters along the travelled arcs, leading to a set of non-linear constraints. Although easily handled through big-M reformulations, the number of such constraints grows quickly such that the resulting continuous relaxation performs poorly for large problem instances. We resolve this challenge by showing how to significantly reduce their number without affecting optimality. We do so by carefully reformulating the problem over a multigraph.  

Third, the resulting \gls{milp} can still not be solved at scale by state-of-the-art solvers. We tackle this challenge by developing a local-search algorithm that we embed within the branch and bound tree generated when solving the \gls{milp}. Specifically, whenever we find a new incumbent by solving the \gls{milp}, we stop and improve this solution through our local search before feeding it back to warmstart the \gls{milp}. While our local search follows a greedy approach based on delay severity, its development requires careful considerations as efficiently evaluating the impact that modifying one departure time has on all other trips is nontrivial.

Complimentary to this, scaling our solution to complex real-world networks introduces additional computational challenges due to the large number of road segments involved. By utilizing contraction hierarchies \citep[cf.][]{GeisbergerSandersEtAl2008}, we reduce the network's complexity while preserving shortest path lengths. This strategy significantly improves the model's applicability and computational efficiency, making it suitable for studying staggered routing in real-world scenarios.

\subsection{Related work}\label{subsec:related_literature}
Our work connects with various streams of literature for the optimization of \gls{amod} systems, an area that has recently attracted significant attention. Given the sheer volume of literature in this space, providing a comprehensive overview is beyond the scope of this work, and we refer the interested reader to \citet{ZardiniLanzettiEtAl2022} and \citet{NarayananChaniotakisEtAl2020}. Instead, we focus on two streams that are closest to our work: fleet control at the \emph{strategic} and \emph{operational} level. 
 
A significant bulk of the existing literature for \emph{strategic} fleet control focuses on time-invariant mesoscopic formulations that utilize network flow models to describe vehicle routing. Various contributions have focused on such steady-state regimes, applying either capacity constraints to model congestion \citep[see, e.g.,][]{RossiZhangEtAl2018, RossiIglesiasEtAl2020, EstandiaSchifferEtAl2021} or alternatively utilizing a piece-wise approximation of a \gls{bpr} function \citep[see, e.g.,][]{SalazarLanzettiEtAl2019,BahramiRoorda2020,BangMalikopoulos2022}. However, none of these paradigms is suitable for our study as they remain time-invariant and, therefore, do not allow the study of staggered routing.

Existing \emph{operational} control approaches are instead based on fine-grained models that do incorporate a time dimension as they aim at devising an online control policy. Several works have been published in this domain, focusing on planning tasks such as vehicle-to-request dispatching, vehicle rebalancing, and request pooling. In this context, various methodological approaches have been studied, ranging from model predictive control for vehicle dispatching \citep[see, e.g.,][]{IglesiasRossiEtAl2018b, TsaoIglesiasEtAl2018, LiuSamaranayake2022} and request pooling \citep[see, e.g.,][]{AlonsoMoraWallarEtAl2017}, to deep reinforcement learning for rebalancing \citep[see, e.g.,][]{JiaoTangEtAl2021, GammelliYangEtAl2021, SkordilisHouEtAl2021, LiangWenEtAl2021} and dispatching \citep[see, e.g.,][]{XuLiEtAl2018, LiQinEtAl2019, TangQinEtAl2019, ZhouJinEtAl2019, SadeghiEshkevariTangEtAl2022, EndersHarrisonEtAl2023}, to optimization-augmented learning pipelines \citep[see, e.g.,][]{ZhangHuEtAl2017, JungelParmentierEtAl2023}. However, these approaches focus on optimizing the routing decisions over space but do not focus on the possibility of staggering the departures over time, which is the focus of this work. Surprisingly, even studies that explicitly model congestion are scarce as existing research often bases travel times on precise forecasts or deterministic assumptions. 

Finally, our work draws inspiration from the seminal work of \citet{Vickrey1969} and its extensions \citep{LiHuangEtAl2020}, albeit with significant differences. Indeed, while Vickrey's pioneering bottleneck model lays the theoretical groundwork to analyze commuters' departure time decisions to avoid congestion, their model has a purely \emph{descriptive} purpose. To the contrary, our work ventures beyond descriptive modeling of traffic behaviors by actively optimizing departure times to minimize congestion.
\subsection{Roadmap}
The remainder of this paper is structured as follows. Section~\ref{sec:problem_setting} outlines the problem setting for the staggered routing problem. Section~\ref{sec:methodology} details our matheuristic developed for solving the staggered routing problem. In Section~\ref{sec:case_study}, we provide an overview of our experimental design, while in Section~\ref{sec:results}, we present the results obtained from our numerical study. Finally, we conclude in Section~\ref{sec:conclusion} with some remarks and future research directions.

\section{Problem setting}\label{sec:problem_setting}

We begin by considering an offline problem wherein an \gls{amod} operator coordinates an autonomous vehicle fleet to serve a given set of trips, e.g., trips arising within a day or a specific time window within a city. Each trip is operated by a single vehicle and can correspond to both on-duty, i.e., customer delivery, or off-duty, i.e., rebalancing activities, of the vehicle. Importantly, vehicles operate trips along a \emph{prespecified} route, e.g., the shortest, within the road network.
Within this context, we focus on enhancing operational efficiency \emph{after} vehicle routing decisions are taken. That is, we assume that the origin, destination, and route for each trip (equivalently, vehicle) are fixed.
To do so, we aim to stagger trip departures, i.e., postpone trip departure times beyond their originally requested times, to reduce congestion caused by route overflows and thus minimize the fleet's travel time. 

We model the road network within which the fleet operates as a directed graph ${\Graph=(\SetNodes, \SetArcs)}$. The nodes~$\SetNodes$ represent trip origins, destinations, and road intersections, while the arcs $\SetArcs$ are the road segments connecting these nodes.  Given a set of trips $\SetTrips$, let a quadruple $(\TripPath{\Trip}, \TripReleaseTime{\Trip}, \TripDeadline{\Trip}, \MaxStaggering{\Trip})$ define a trip $\Trip \in \SetTrips$, with $\TripPath{\Trip}$ denoting the sequence of arcs composing the trip's fixed route through the network, $\TripReleaseTime{\Trip}$ being the earliest departure from the trip's origin, $\TripDeadline{\Trip}$ being the latest acceptable arrival at the destination, and $\MaxStaggering{\Trip}$ being the maximum staggering time applicable to trip departure. 

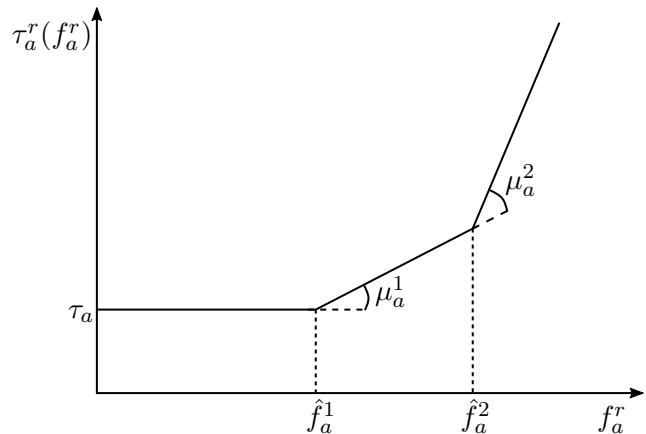
\begin{figure}[!b]    
    \centering
    %% Creator: Inkscape 1.1.1 (3bf5ae0d25, 2021-09-20), www.inkscape.org
%% Accompanies image file 'latency.pdf' (pdf, eps, ps)
\begingroup%
  \makeatletter%
  \providecommand\color[2][]{%
    \errmessage{(Inkscape) Color is used for the text in Inkscape, but the package 'color.sty' is not loaded}%
    \renewcommand\color[2][]{}%
  }%
  \providecommand\transparent[1]{%
    \errmessage{(Inkscape) Transparency is used (non-zero) for the text in Inkscape, but the package 'transparent.sty' is not loaded}%
    \renewcommand\transparent[1]{}%
  }%
  \providecommand\rotatebox[2]{#2}%
  \newcommand*\fsize{\dimexpr\f@size pt\relax}%
  \newcommand*\lineheight[1]{\fontsize{\fsize}{#1\fsize}\selectfont}%
  \ifx\svgwidth\undefined%
    \setlength{\unitlength}{235.36627101bp}%
    \ifx\svgscale\undefined%
      \relax%
    \else%
      \setlength{\unitlength}{\unitlength * \real{\svgscale}}%
    \fi%
  \else%
    \setlength{\unitlength}{\svgwidth}%
  \fi%
  \global\let\svgwidth\undefined%
  \global\let\svgscale\undefined%
  \makeatother%
  \begin{picture}(1,0.76723485)%
    \lineheight{1}%
    \setlength\tabcolsep{0pt}%
    \put(-0.00268863,0.6224658){\makebox(0,0)[lt]{\lineheight{1.25}\smash{\begin{tabular}[t]{l}$\ArcTravelTime{\Arc}{\Trip}$\end{tabular}}}}%
    \put(0.0865341,0.18376114){\makebox(0,0)[lt]{\lineheight{1.25}\smash{\begin{tabular}[t]{l}$\ArcNominalTravelTime{\Arc}$\end{tabular}}}}%
    \put(0.57309445,0.21233158){\makebox(0,0)[lt]{\lineheight{1.25}\smash{\begin{tabular}[t]{l}$\PWLSlope{\Arc}{1}$\end{tabular}}}}%
    \put(0.77900341,0.39002863){\makebox(0,0)[lt]{\lineheight{1.25}\smash{\begin{tabular}[t]{l}$\PWLSlope{\Arc}{2}$\end{tabular}}}}%
    \put(0.924868,0.01100346){\makebox(0,0)[lt]{\lineheight{1.25}\smash{\begin{tabular}[t]{l}$\Flow{\Arc}{\Trip}$\end{tabular}}}}%
    \put(0.46354327,0.01100346){\makebox(0,0)[lt]{\lineheight{1.25}\smash{\begin{tabular}[t]{l}$\ThFlow{\Arc}{1}$\end{tabular}}}}%
    \put(0.70853818,0.01100346){\makebox(0,0)[lt]{\lineheight{1.25}\smash{\begin{tabular}[t]{l}$\ThFlow{\Arc}{2}$\end{tabular}}}}%
    \put(0,0){\includegraphics[width=\unitlength,page=1]{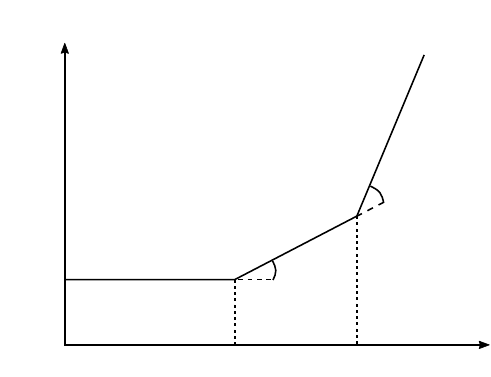}}%
  \end{picture}%
\endgroup%

    \caption{Travel time of trip $\Trip$ on arc $\Arc$ as a function of the number of \gls{amod} trips $\Flow{\Arc}{\Trip}$.
    \label{fig:pwl_func}}
\end{figure}

A trip $\Trip$ on arc $\Arc$ incurs a travel time $\ArcTravelTime{\Arc}{\Trip}$ that is the sum of a free-flow travel time $\ArcNominalTravelTime{\Arc}$ and a congestion-induced delay, which we determine by the number of \gls{amod} trips $\Flow{\Arc}{\Trip}$ on arc $\Arc$ when trip $\Trip$ starts traversing the arc.\footnote{Note that our formulation easily allows to consider the background presence of non-\gls{amod} traffic, see Section~\ref{sec:case_study} for details.} %\footnote{As anticipated in Section~\ref{sec:our-approach}, when formulating our problem, we do \emph{not} keep track of the number of trips at fixed locations over different time steps (the commonly adopted Eulerian approach), but instead follow each trip and keep track only of the number of trips encountered by such vehicle when entering a new arc (the Lagrangian approach).} 
We model such a congestion-induced delay as a piece-wise linear, non-decreasing, and convex function, resulting in a travel time over each arc such as that represented in Figure~\ref{fig:pwl_func}. Without loss of generality, the total travel time encountered by trip $\Trip$ on arc $\Arc$ can thus be written as 
\begin{align}\label{eq:travel_time_approximation}
\ArcTravelTime{\Arc}{\Trip} = \ArcNominalTravelTime{\Arc} +  
\max_{k\in \NumPieces} \{0,\PWLSlope{\Arc}{k}\cdot(\Flow{\Arc}{\Trip}-\ThFlow{\Arc}{k})\},
\end{align}
where $\ArcNominalTravelTime{\Arc}>0$, $\PWLSlope{\Arc}{k}>0$ and $\ThFlow{\Arc}{k}>0$

\begin{remark}
Two comments are in order. 
First, while other congestion models exist, e.g., simpler capacity models \citep{EstandiaSchifferEtAl2021} or more accurate microscopic models \citep{LevinKockelmanEtAl2017}, we believe that our approach to modeling congestion provides a good balance between accuracy and computational complexity that is appropriate for the mesoscopic nature of our study.
In this respect, it is important to remark that, while in mesoscopic studies, the relation between travel time and congestion is often taken to be a fourth-order polynomial \citep[see, e.g.,][]{USBoPR1964}, our approach not only approximates such a function in a computationally convenient way, but it also allows to accommodate other, potentially different, relations.
Second, when detecting the congestion level for a vehicle entering the arc, we ignore whether the other vehicles have just begun traversing the arc or are nearing completion. We believe that the resulting trade-off between slightly overestimating congestion and deriving a tractable model is reasonable for the scope of our study.
\end{remark}

In this setting, a solution $\Solution$ is a vector with a departure time for each trip $\Trip$. As we assume no vehicle idles once its trip has started, this information suffices to determine both the time as well as the number of vehicles $\Flow{\Arc}{\Trip}$ encountered by every vehicle when entering arc $a$ of its route $\TripPath{\Trip}$. In turn, this sequence determines the travel times over all arcs and, thus, the trips' arrival times.

We require a solution $\Solution$ to satisfy the following constraints, which we collect in the feasible set $\FeasibleSet$:
\begin{enumerate}
\item[i)] The departure time of each trip $\Trip$ must occur after $\TripReleaseTime{\Trip}$.
\item[ii)] The departure time of each trip $\Trip$ can only be staggered for a maximum amount of time $\MaxStaggering{\Trip}$.
\item[iii)] The departures in $\Solution$ must ensure that the arrival at the destination of each trip $\Trip$ occurs before $\TripDeadline{\Trip}$.
\end{enumerate}

In this work, we seek a solution $\Solution^{\star}\in\FeasibleSet$ that minimizes the total fleet travel time, i.e., the sum of the travel times each trip incurs over all arcs of its path
\begin{equation}\label{eq:problem}
\begin{aligned}
\min_{\Solution \in \Pi} &\; \TotalDelay{\Solution} = \sum_{\Trip \in \SetTrips}\sum_{\Arc\in \TripPath{\Trip}} \tau_{\Arc}^{\Trip}(\pi),
\end{aligned}
\end{equation}
where, with slight abuse of notation, we explicitly denote the fact that the travel time $\tau_{\Arc}^{\Trip}$ encountered on arc $\Arc$ by trip $\Trip$ depends on the entire vector of departure times.

Finally, while \eqref{eq:problem} describes the offline version of the problem we will consider, one can easily transform this offline setting into an online problem by discretizing time and taking decisions on batches of trips that arrive in the system according to a point process. In our numerical studies, we will comment on the difference between the maximum potential of staggering identified by solving the offline problem and the improvement attained by solving its online counterpart.
\section{Methodology}\label{sec:methodology}
 The staggered routing problem as introduced in Section~\ref{sec:problem_setting} is NP-hard, which we show in Appendix~\ref{sec:appendix_hardness}. This motivates us to develop an efficient matheuristic that allows to solve large-scale instances in the following. First, we introduce a \gls{milp} that allows us to obtain solutions for Problem~\eqref{eq:problem} in Section~\ref{subsec:milp_formulation}. Towards this goal, we first introduce a basic model that uses an indicator function to quantify the aggregated flow, before we elaborate on how to model the respective indicator constraints efficiently. Second, we present a matheuristic that combines our \gls{milp} with a local search scheme to solve large real-world instances in Section~\ref{subsec:matheuristic}. Finally, Section~\ref{subsec:online_algorithm} adapts the previously developed algorithm to solve the online counterpart of our problem in a rolling-horizon fashion.
 \vspace{-.5cm}
\subsection{Mixed integer linear program}\label{subsec:milp_formulation}
Our objective is to formulate \eqref{eq:problem} as a \gls{milp}. Towards this goal, note that minimizing the total travel time is equivalent to minimizing the congestion-induced delay as we have no control over the arcs' free-flow travel times $\ArcNominalTravelTime{\Arc}$. With this in mind, we perform an epigraphic reformulation of the objective function which we achieve by introducing $\DelayOnArc{\Arc}{\Trip}$ and requesting that $\DelayOnArc{\Arc}{\Trip}\ge \PWLSlope{\Arc}{k} \cdot (\Flow{\Arc}{\Trip} -\ThFlow{\Arc}{k}$) for all arcs, trips, and $k \in \NumPieces$.
To represent the entry time of trip $\Trip$ on arc $\Arc$ on its path, we introduce the continuous variable $\TripDepartureOnArc{\Arc}{\Trip}$. Further, we let $\Successor{\Arc}{\Trip}$ denote the successor arc of arc $\Arc$ on route $\TripPath{\Trip}$, and $\FirstArc$, $\LastArc$ denote the first and last arcs of trip $\Trip$,  respectively. Finally, we let $\SetTripsOnArc{\Arc}$ denote the set of all trips that transit through arc $\Arc$.
%we let $\TripPathBar{\Trip}$ be the route of trip $\Trip$ without the last arc, 
%To introduce our basic model, we use notation as introduced in Section~\ref{sec:problem_setting}. For a summarized view of the variable definitions, we refer the reader to Appendix~\ref{sec:appendix_vars_table}. 
Then, problem \eqref{eq:problem} can be written as:
%also
\begin{subequations}\label{milp:srp}
\begin{align}
	\min \; &\sum_{\Trip \in \SetTrips}\sum_{\Arc \in \TripPath{\Trip}} \DelayOnArc{\Arc}{\Trip}, \label{basic_model:obj_func_model}\\
	\text{s.t. }& \nonumber\\ 
 	&{\DelayOnArc{\Arc}{\Trip}} \geq \PWLSlope{\Arc}{k} \cdot \left( \Flow{\Arc}{\Trip} - \ThFlow{\Arc}{k} \right),  &\quad \forall k \in \NumPieces, \Trip \in \SetTrips, \Arc \in \TripPath{\Trip}\label{basic_model:delay_constraint_first_model}\\ 
    & \TripDepartureOnArc{\Arc'}{\Trip} = \TripDepartureOnArc{\Arc}{\Trip}  + \ArcNominalTravelTime{\Arc} + \DelayOnArc{\Arc}{\Trip}, &\quad \forall \Trip \in \SetTrips, \Arc \in \TripPath{\Trip}\setminus \LastArc, \Arc' =\Successor{\Arc}{\Trip}\label{basic_model:continuity_trips_model}\\
      & \Flow{\Arc}{\Trip} = \sum_{
      \SecondTrip \in \SetTripsOnArc{\Arc}\setminus\{\Trip\}
      }\mathbb{I}(\TripDepartureOnArc{\Arc}{\SecondTrip} \leq \TripDepartureOnArc{\Arc}{\Trip} < \TripDepartureOnArc{\Arc}{\SecondTrip} +  \ArcNominalTravelTime{\Arc} + \DelayOnArc{\Arc}{\SecondTrip}), &\quad \forall \Trip \in \SetTrips, \Arc \in \TripPath{\Trip}\label{basic_model:indicator_function}\\
    &  \TripEarliestDepartureTime{}{\Trip} \leq \TripDepartureOnArc{a}{\Trip} \leq  \TripEarliestDepartureTime{}{\Trip} + \MaxStaggering{\Trip},  &\quad \forall \Trip \in \SetTrips, \Arc=\FirstArc \label{basic_model:time_windows1}\\
  & \TripDepartureOnArc{\Arc}{\Trip} + \ArcNominalTravelTime{\Arc} + \DelayOnArc{\Arc}{\Trip}\leq \TripLatestExitTime{}{\Trip}, &\quad \forall \Trip \in \SetTrips, \Arc=\LastArc  \label{basic_model:time_windows3}\\
  &\DelayOnArc{\Arc}{\Trip}, \TripDepartureOnArc{\Arc}{\Trip}, \Flow{\Arc}{\Trip} \geq 0. &\quad \forall \Trip \in \SetTrips, \Arc \in \TripPath{\Trip} \label{basic_model:positive_departures_model}
\end{align}
\end{subequations}
%\todo{can one turn 3g in a way so that the one does not need to take additionally lamba in the big M?}
Our objective~\eqref{basic_model:obj_func_model} jointly with constraints~\eqref{basic_model:delay_constraint_first_model} minimizes the total delay over all trips and arcs. Constraints~\eqref{basic_model:continuity_trips_model} propagate forward the time at which trip $\Trip$ enters arc $\Arc$ along its route $\TripPath{\Trip}$, while constraints~\eqref{basic_model:indicator_function} calculate the number of vehicles encountered by each trip $\Trip$ on arc $\Arc$. This is achieved using an indicator function $\mathbb{I}(\cdot)$, which is set to one if trip $\Trip$ enters arc $\Arc$ when $\SecondTrip$ is traversing the same arc. Constraints~\eqref{basic_model:time_windows1}-\eqref{basic_model:time_windows3} ensure that the start times and drop-off times meet the requirements. 
Constraints~\eqref{basic_model:positive_departures_model} state the variable domains.

To convert the latter model into a \gls{milp}, we replace the indicator function in Constraints~\eqref{basic_model:indicator_function} with big-M constraints that replicate its logic. To do so, we introduce binary variables $\VarAlpha{\Arc}{\Trip}{\SecondTrip},\VarBeta{\Arc}{\Trip}{\SecondTrip}$ and $\VarGamma{\Arc}{\Trip}{\SecondTrip}$ for each pair of distinct trips $(\Trip,\SecondTrip)$ traversing the same arc, whose logic is as follows:
\begin{equation*}\label{eq:binaries_conditions}
    \begin{aligned}
        \VarAlpha{\Arc}{\Trip}{\SecondTrip} &=\begin{cases}
            0, &\hspace{-0.2cm}\text{if } \TripDepartureOnArc{\Arc}{\Trip} < \TripDepartureOnArc{\Arc}{\SecondTrip}\\
            1, &\hspace{-0.2cm}\text{otherwise} 
        \end{cases}
    \end{aligned}
    \quad 
    \begin{aligned}
        \VarBeta{\Arc}{\Trip}{\SecondTrip} &=\begin{cases}
            0, &\hspace{-0.2cm}\text{if } \TripDepartureOnArc{\Arc}{\Trip}\ge \TripDepartureOnArc{\Arc}{\SecondTrip} + \ArcNominalTravelTime{\Arc} + \DelayOnArc{\Arc}{\SecondTrip} \\
            1, &\hspace{-0.2cm}\text{otherwise}
        \end{cases}
    \end{aligned}
    \quad 
    \begin{aligned}
        \VarGamma{\Arc}{\Trip}{\SecondTrip} &=\begin{cases}
            0, &\hspace{-0.2cm}\text{if } \VarAlpha{\Arc}{\Trip}{\SecondTrip} + \VarBeta{\Arc}{\Trip}{\SecondTrip} < 2\\
            1, &\hspace{-0.2cm}\text{otherwise}
        \end{cases}
    \end{aligned}
\end{equation*}

We enforce this logic and therefore evaluate $\Flow{\Arc}{\Trip}$ by replacing constraints~\eqref{basic_model:indicator_function} with the following:
\vspace{-.6cm}
\addtocounter{equation}{-1} % Adjust the equation counter to continue from the last equation of the 
\begin{subequations}\label{set_constr:check_conflicts}
\setcounter{equation}{7}
\begin{align}
    &\TripDepartureOnArc{\Arc}{\Trip} - \TripDepartureOnArc{\Arc}{\SecondTrip} + \SmallConstant \leq \BigM{1} \cdot \VarAlpha{\Arc}{\Trip}{\SecondTrip}, &&\forall \Arc \in \SetArcs, \Trip, \SecondTrip \in \SetTripsOnArc{\Arc}, \Trip\neq\SecondTrip  \label{constr:original_alpha_constraint_first}\\
    &\TripDepartureOnArc{\Arc}{\SecondTrip} - \TripDepartureOnArc{\Arc}{\Trip}  \leq \BigM{2} \cdot (1 - \VarAlpha{\Arc}{\Trip}{\SecondTrip}), &&\forall \Arc \in \SetArcs, \Trip, \SecondTrip \in \SetTripsOnArc{\Arc}, \Trip\neq\SecondTrip  \label{constr:original_alpha_constraint_second}\\
    &\TripDepartureOnArc{\Arc}{\SecondTrip} + \ArcNominalTravelTime{\Arc} + \DelayOnArc{\Arc}{\SecondTrip} - \TripDepartureOnArc{\Arc}{\Trip}  \leq \BigM{3} \cdot \VarBeta{\Arc}{\Trip}{\SecondTrip}, &&\forall \Arc \in \SetArcs, \Trip, \SecondTrip \in \SetTripsOnArc{\Arc}, \Trip\neq\SecondTrip \label{constr:original_beta_constraint_first}\\
    &\TripDepartureOnArc{\Arc}{\Trip} - \TripDepartureOnArc{\Arc}{\SecondTrip} - \ArcNominalTravelTime{\Arc} - \DelayOnArc{\Arc}{\SecondTrip} + \SmallConstant  \leq \BigM{4} \cdot (1 -\VarBeta{\Arc}{\Trip}{\SecondTrip}), &&\forall \Arc \in \SetArcs, \Trip, \SecondTrip \in \SetTripsOnArc{\Arc}, \Trip\neq\SecondTrip \label{constr:original_beta_constraint_second}\\
   & \VarGamma{\Arc}{\Trip}{\SecondTrip} - \VarAlpha{\Arc}{\Trip}{\SecondTrip} - \VarBeta{\Arc}{\Trip}{\SecondTrip}  \geq - 1,  &&\forall \Arc \in \SetArcs, \Trip, \SecondTrip \in \SetTripsOnArc{\Arc}, \Trip\neq\SecondTrip \label{constr:original_gamma_constraint_first}\\
    &2\VarGamma{\Arc}{\Trip}{\SecondTrip} - \VarAlpha{\Arc}{\Trip}{\SecondTrip} - \VarBeta{\Arc}{\Trip}{\SecondTrip}  \leq 0, &&\forall \Arc \in \SetArcs, \Trip, \SecondTrip \in \SetTripsOnArc{\Arc}, \Trip\neq\SecondTrip \label{constr:original_gamma_constraint_second}\\
    & \Flow{\Arc}{\Trip} = \sum_{
      \SecondTrip \in \SetTripsOnArc{\Arc}\setminus\{\Trip\}} \VarGamma{\Arc}{\Trip}{\SecondTrip} && \forall \Arc \in \SetArcs, \Trip \in \SetTripsOnArc{\Arc} \label{constr:compute_number_of_conflicts}
\end{align} 
\end{subequations}
\noindent where $\SmallConstant$ is a small constant and  $\BigM{i}$ are sufficiently large numbers.\footnote
{Note that the use of a small constant $\SmallConstant$ is necessary. Indeed, in its absence, given a pair of distinct trips $(\Trip,\SecondTrip)$ jointly traversing arc $\Arc$, for $\TripDepartureOnArc{\Arc}{\Trip}$ being equal to $\TripDepartureOnArc{\Arc}{\SecondTrip}$, setting $\VarAlpha{\Arc}{\Trip}{\SecondTrip}$ to either zero or one would satisfy Constraints~\eqref{constr:original_alpha_constraint_first}-\eqref{constr:original_alpha_constraint_second}. Details on our choice of  $\BigM{i}$ can be found in Section~\ref{sec:constr_redefinition}.}
%\color{red}\noindent [DP: let's de-emphasize this. just a small comment on the footnote is enough and no need to give all details. The deatils can be given when you explain how the big M is chosen if you wich.]
First, observe that the resulting program is indeed a \gls{milp}. However, this formulation presents significant computational challenges due to the number of big-M constraints required to %\textcolor{red}{[DP: I would rather avoid the jargon ``conflict'' until necessary, it's unclear what this even means. For this purpose, I have removed all the references and rephrased them here. To complete the job, one needs to rephrase the additional two occurrences in red on this page, and then once and for all introduce what ``conflict'' means in the beginning of 3.1.1] \st{manage conflicts between trip pairs}} 
check if a given pair of vehicles travels simultaneously on the same arc. %{[AC: added to next section definition of conflict]}.
Specifically, each constraint in \eqref{constr:original_alpha_constraint_first}-\eqref{constr:original_gamma_constraint_second} needs to be repeated for all pairs of trips that transit through each given arc.
%Specifically, the number of necessary constraints grows as %$\sum_a \sum_{r\in \SetTripsOnArc{\Arc}} |\SetTripsOnArc{\Arc}|\cdot |\SetTripsOnArc{\Arc}|-1$
% \textcolor{yellow}{$\BigO{|\SetArcs| \cdot |\SetTripsOnArc{\Arc}|^2}$},  which approximates the number of ways to select two distinct elements within the set $\SetTripsOnArc{\Arc}$ for every arc in $\SetArcs$. 
This formulation results, unfortunately, in prohibitive solution times, even for modest problem sizes, as the continuous relaxation used in the solution of the \gls{milp} performs poorly. To enhance computational efficiency, we leverage the problem structure to identify and remove as many redundant big-M constraints as possible in the following.

A key observation is that not every distinct pair of trips traversing the same arc may influence each other's travel time, largely because each trip is constrained by specific time windows. However, the original problem formulation introduced in Section~\ref{sec:problem_setting} only defines these windows for each entire trip, not for individual arcs. Based on these trip-specific time windows, our approach constructs arc-specific time windows to identify trips that, by traversing the arc simultaneously, may incur a delay. 

We then leverage these bounds to build a representation of the problem on a multidigraph, in which parallel arcs represent identical road segments, but each arc is associated with either a set of potentially conflicting trips, i.e., trips that may influence each other's travel time and thus contain conflicting arcs, i.e., arcs on which trip induced congestion may arise, or trips that will not delay each other, containing only non-conflicting arcs. This differentiation allows us to apply big-M constraints selectively, focusing only on conflicting arcs and relevant trip pairs, thereby significantly reducing the number of big-M constraints appearing in the \gls{milp}. Additionally, we explore methods to merge or exclude non-conflicting arcs from the model, further minimizing its complexity and the computational overhead of the local search method discussed later in the section.

In the rest of this section, we will proceed as follows. In Section~\ref{sec:arc_dependent_time_windows}, we detail how to determine arc-specific time windows for each trip. In Section~\ref{sec:conflicting_sets} and~\ref{sec:graph_processing}, we detail how we derive the multigraph representation introduced above. Building upon this, we are finally able to reduce the number of big-M constraints needed by our \gls{milp} model in Section~\ref{sec:constr_redefinition}.

\subsubsection{Determining arc-dependent time windows}\label{sec:arc_dependent_time_windows}
%{To ease the discussion, we introduce the concept of \emph{trip conflict}. We say that a trip $\Trip$ conflicts with $\SecondTrip$ on arc $\Arc$ if ${\TripDepartureOnArc{\Arc}{\SecondTrip}} \leq {\TripDepartureOnArc{\Arc}{\Trip}} < {\TripDepartureOnArc{\Arc}{\SecondTrip} + \ArcNominalTravelTime{\Arc} + \DelayOnArc{\Arc}{\SecondTrip}}$, i.e., if trip $\Trip$ begins traversing arc $\Arc$ after $\SecondTrip$ enters and before $\SecondTrip$ exits the arc. This definition is asymmetric: if $\Trip$ conflicts with $\SecondTrip$ on $\Arc$, then the inverse is not true, i.e. $\SecondTrip$ does not conflict with $\Trip$ on $\Arc$.}
To distinguish trips that may conflict from trips that do not, we construct upper and lower bounds on the entry and exit times for each given arc on a trip. At a high level, the goal is to make these bounds as tight as possible while ensuring that every feasible solution satisfies them, which allows to reduce the number of big-M constraints. 

Specifically, we let $\EarliestEntry{\Arc}{\Trip}$, $\LatestEntry{\Arc}{\Trip}$ represent the earliest and latest \emph{entry} time of trip $\Trip$ in arc $\Arc$. Similarly, we use $\EarliestExit{\Arc}{\Trip}$, $\LatestExit{\Arc}{\Trip}$ to represent earliest and latest \emph{exit} times. Note, however, that knowledge of the earliest entry time $\EarliestEntry{\Arc}{\Trip}$ and latest exit time $\LatestExit{\Arc}{\Trip}$ for all arcs on a trip suffices to reconstruct also the other quantities, i.e., the latest entry and the earliest exit over all the arcs on that trip. In fact, by continuity, the latest entry on an arc $\LatestEntry{\Arc}{\Trip}$ must equal the latest exit from the previous arc on that trip, while the earliest exit from an arc $\EarliestExit{\Arc}{\Trip}$ must equal the earliest entry on the subsequent arc on that trip. This is the case for all arcs in the trip except for the very first and last arc. For the first arc $\Arc = \FirstArc$, we set the latest entry as $\LatestEntry{\Arc}{\Trip} = \TripReleaseTime{\Trip} + \MaxStaggering{\Trip}$, while for the last arc $\Arc = \LastArc$ we determine the earliest exit by letting the vehicle travel in free flow through that arc so that $\EarliestExit{\Arc}{\Trip}=\EarliestEntry{\Arc}{\Trip}+\ArcNominalTravelTime{\Arc}$. 

Thus, in the following, we focus on determining only the earliest entry $\EarliestEntry{\Arc}{\Trip}$ and the latest exit $\LatestExit{\Arc}{\Trip}$ times. The idea is to initialize these quantities by propagating the earliest entry times over subsequent arcs on a trip as if vehicles were traveling in free-flow, starting from the earliest departure $\TripReleaseTime{\Trip}$ at the first arc. Similarly, for the latest exit times, which we propagate backward over preceding arcs on a trip, starting from the latest arrival at the last arc $\TripDeadline{\Trip}$. In a second phase, we then refine and tighten these time windows by bounding the minimum and maximum congestion encountered by each trip over its arcs. We now describe both of these steps in further detail and refer to Algorithm~\ref{algo:compute_bounds} for a formal overview of the approach.

\begin{algorithm}[!b]
\footnotesize
\caption{Compute arc-specific time windows}\label{algo:compute_bounds}
\textbf{Input:} Set of trips $\SetTrips$
\begin{algorithmic}[1]
\State $\Queue,~\EarliestEntryVector,~\LatestExitVector,~\MinDelayVector \gets$ \textsc{initialize}($\SetTrips$)
\While{$\Queue$ is not empty}
\State $(\Arc,\Trip) \gets \Queue.\textsc{pop()}$
\Comment{Extract tuple with highest priority.}
\State $\EarliestEntryVector,~\MinDelayVector \gets$ \textsc{updateEarliestEntryTimes}($\Arc,\Trip$)
\State $\LatestExitVector \gets$ \textsc{updateLatestExitTimes}($\Arc,\Trip$)
\EndWhile
\State $\LowerBound \gets \textsc{getInitialLowerBound}(\MinDelayVector)$
\State \textbf{return} ($\EarliestEntryVector,~\LatestExitVector, \LowerBound$)
\end{algorithmic}
\end{algorithm}

\vspace{-.5cm}
\subparagraph{Initialization.}
For each trip, we set the earliest entry time $\EarliestEntry{\Arc}{\Trip}$ for its first arc equal to the earliest departure time of the trip. We then calculate the earliest entry times for the subsequent arc by adding the nominal travel time of the preceding arc, that is
$
\EarliestEntry{\Successor{\Arc}{\Trip}}{\Trip} = 
\ubar{e}_{\Arc}^{\Trip} + \ArcNominalTravelTime{\Arc}$
for all $\Arc \in \TripPath{\Trip}\setminus{\LastArc}$.
%$\forall \Trip \in \SetTrips, \Arc \in \TripPath{\Trip}\setminus{\FirstArc}$. 
For each trip, we set the latest exit time $\LatestExit{\Arc}{\Trip}$ for its last arc equal to the latest arrival time of the trip $\TripDeadline{\Trip}$. We then calculate the latest exit time for the preceding arc by subtracting the nominal travel time of the arc, that is
$
\LatestExit{\Arc}{\Trip} = 
\LatestExit{\Successor{\Arc}{\Trip}}{\Trip} - \ArcNominalTravelTime{\Successor{\Arc}{\Trip}}$
for all $\Arc \in \TripPath{\Trip}\setminus{\LastArc}$.
%$\forall \Trip \in \SetTrips, \Arc \in \TripPath{\Trip}\setminus{\FirstArc}$. 

We then construct a queue whose elements are the tuples indexed by $(\Arc,\Trip)$, i.e., arc and trip, and sort such queue based on the smallest earliest entry times. We then process, in order, one element of the queue at a time to update its earliest entry and latest exit time. 
First, we describe the method for updating time windows. Following that, we detail how employing a priority queue results in tighter bounds compared to an unsorted approach. %\textcolor{red}{[DP: Antonio mentions that this sorting is done so that after we have processed one element, by processing subsequent elements, we will never change the bounds of previous ones. This is unclear to me and needs to be explained/justified clearly.]}
%{[I added the rationale behind the priority queue at the end of this section.].}

\vspace{-.5cm}
\subparagraph{Updating earliest entry times.} 
%{[AC: I changed this discussion to reflect better what I do in my implementation. The version of Dario is in comment. I mainly disagreed on the point of determining the minimum congestion on the current arc by looking at the previous arc -- this cannot be done on the first arc, for example. What I do is looking at the congestion on the current arc and update the earliest entry on all the subsequent arcs.]}

% DARIO VERSION
% Given an element $(\Arc,\Trip)$, we now update, i.e., move forward in time, its earliest entry time, so as to tighten the time windows. 
% Towards this goal, we compute the minimum congestion that trips $\Trip$ inevitably incurs on the arc $a'$ which precedes $a$, i.e., for which $\Arc=\delta^\Trip(\Arc')$, and denote it with $\ubar{f}_{\Arc'}^\Trip$. This quantity can be determined, as follows, by leveraging the existing time windows, 
% \begin{align*}
% \ubar{f}_{\Arc'}^\Trip = \sum_{\substack{\SecondTrip \in \SetTripsOnArc{\Arc'} \\ \hspace{-6.5pt}\SecondTrip \neq \Trip}}\mathbb{I}(\ubar{e}_{\Arc'}^\Trip>\bar{e}_{\Arc'}^{\Trip'}   \wedge 
% \bar{e}^\Trip_{\Arc'}<\ubar{\ell}_{\Arc'}^{\Trip'}).
% \end{align*}
% Then, the earliest entry time over arc $a$ is given by the earliest exit possible from the preceding arc, which can be computed by summing its earliest entry time to the travel time on that arc due to the presence of $\ubar{f}_\Arc^\Trip$, that is,
% $\ubar{e}_\Arc^\Trip \gets \ubar{e}_{\Arc'}^\Trip + \tau_{\Arc'}^\Trip(\ubar{f}_{\Arc'}^\Trip)$.

% ANTONIO VERSION
Given an element $(\Arc,\Trip)$, we now update, i.e., move forward in time, the earliest entry times on arcs subsequent to $\Arc$ to tighten the time windows. 
Towards this goal, we compute the minimum conflicts that trips $\Trip$ inevitably incur on arc $\Arc$ and denote it with $\MinFlow{\Arc}{\Trip}$. This quantity can be determined by leveraging the existing time windows as
\begin{align*}
\MinFlow{\Arc}{\Trip} = \sum_{\SecondTrip \in \SetTripsOnArc{\Arc} \setminus \{ \Trip\}}\mathbb{I}(\EarliestEntry{\Arc}{\Trip} > \LatestEntry{\Arc}{\SecondTrip}   \wedge 
\LatestEntry{\Arc}{\Trip}<\EarliestExit{\Arc}{\SecondTrip}).
\end{align*}

We denote $\MinDelay{\Arc}{\Trip}$ as the delay that $\Trip$ incurs on arc $\Arc$ in presence of $\MinFlow{\Arc}{\Trip}$. As the earliest entries on arcs subsequent to $\Arc$ are surely delayed by $\MinDelay{\Arc}{\Trip}$, we update $\EarliestEntry{\Arc'}{\Trip} \gets \EarliestEntry{\Arc'}{\Trip} + \MinDelay{\Arc}{\Trip}$ for every arc $\Arc' \in \TripPath{\Trip}$ such that $\EarliestEntry{\Arc'}{\Trip} > \EarliestEntry{\Arc}{\Trip}$.

%\textcolor{red}{[DP: Antonio, are you doing this only for one arc in the route (extracted according to the priority queue), or are you doing this for all arcs in that trip? unclear]}
%\textcolor{blue}{[AC: I do this for every subsequent arc to the arc extracted from the priority queue: I updated the text to reflect this better. Note that I haven't used $\Successor{\Arc}{\Trip}$ because in the updates we need to maintain the reference to the minimum delay $ \MinDelay{\Arc}{\Trip}$ on the current arc $\Arc$]}
\vspace{-.5cm}
\subparagraph{Updating latest exit times.} To update the latest exit time $\LatestExit{\Arc}{\Trip}$, we determine the maximum number of conflicts, $\MaxFlow{\Arc}{\Trip}$, that a trip $\Trip$ can encounter on arc $\Arc$. 
A simple counting argument reveals that this is equal to the number of other trips $\SecondTrip \in \SetTripsOnArc{\Arc}$ whose time windows $[\EarliestEntry{\Arc}{\SecondTrip}, \LatestExit{\Arc}{\SecondTrip}]$ intersect $[\EarliestEntry{\Arc}{\Trip},\LatestEntry{\Arc}{\Trip}]$.
We then update the latest exit time from arc $\Arc$ based on the latest entry in such arc and the maximum number of conflicts computed, that is, $
\LatestExit{\Arc}{\Trip} \gets \min(\LatestExit{\Arc}{\Trip}, \LatestEntry{\Arc}{\Trip} + \tau_\Arc^\Trip(\MaxFlow{\Arc}{\Trip})).$

\vspace{-.5cm}
\subparagraph{Computing an initial lower bound.} 
We conclude Algorithm~\ref{algo:compute_bounds} by calculating an initial lower bound $\LowerBound$, which already allows evaluating the quality of a solution $\Solution$. After finalizing the computation of arc-dependent time windows, the initial lower bound results from summing up the minimum delays that trips unavoidably incur on each arc:
\begin{align*}
\LowerBound = \sum_{\Trip \in \SetTrips}\sum_{\Arc \in \TripPath{\Trip}} \MinDelay{\Arc}{\Trip}.
\end{align*}
%\textcolor{red}{DP: Unclear what Maintaining the queue below is doing; see also my previous comment}

\vspace{-.5cm}
\subparagraph{The priority queue.} 
 
The priority queue, organized by ascending earliest entry times $\EarliestEntry{\Arc}{\Trip}$, improves upon an unsorted approach to processing the sequence of the time window for two main reasons:
\begin{itemize}
    \item[i)]  Processing pair $(\Arc, \Trip)$ enables us to tighten both the earliest entry times for arcs after $\Arc$ in $\Trip$'s route, as well as the latest exit time from $\Arc$. In this context, proceeding in an ascending order aligns with this forward-looking update mechanism.
    \item[ii)] Ordering by $\EarliestEntry{\Arc}{\Trip}$ allows to tighten the maximum number of conflicts $\MaxFlow{\Arc}{\Trip}$. As seen above, $\Trip$ may conflict with $\SecondTrip$ on $\Arc$ if $[\EarliestEntry{\Arc}{\Trip}, \LatestEntry{\Arc}{\Trip}]$ overlaps with $[\EarliestEntry{\Arc}{\SecondTrip}, \LatestExit{\Arc}{\SecondTrip}]$. At the time the algorithm processes pair $(\Arc,\Trip)$, the forward-looking update logic highlighted in point i) ensures $\EarliestEntry{\Arc}{\Trip}$, $\EarliestEntry{\Arc}{\SecondTrip}$ and $\LatestEntry{\Arc}{\Trip}$ to be as tight as possible. Nevertheless, $\LatestExit{\Arc}{\SecondTrip}$ may or may not be as tight as possible, depending on the following case distinction:
    \begin{itemize}
            \item[1)] If $\EarliestEntry{\Arc}{\Trip} > \EarliestEntry{\Arc}{\SecondTrip}$, the algorithm processed $\SecondTrip$ before $\Trip$ on $\Arc$. Thus, $\LatestExit{\Arc}{\SecondTrip}$ is as tight as possible consequent to the update logic of point i).
        \item[2)] If $\EarliestEntry{\Arc}{\Trip} < \EarliestEntry{\Arc}{\SecondTrip}$, $\Trip$ is processed before $\SecondTrip$. Although $\LatestExit{\Arc}{\SecondTrip}$ is not as tight as possible, the determinant for overlap is $\EarliestEntry{\Arc}{\SecondTrip}$, which is as tight as possible.
    \end{itemize}
\end{itemize}
Note that changes to a trip's earliest entry time can disrupt the queue's ordering. To address this, we maintain the queue's priorities to ensure that tuples are processed in an ascending time sequence.

\vspace{-.5cm}
\subparagraph{Discussion.} 
Notice that our current method tends to overestimate the worst-case delay by focusing solely on the bounds related to the specific arc under consideration. We do not dynamically adjust these bounds based on changing scenarios across consecutive arcs. This conservative approach might not capture the interplay of delays across arcs. As the scalability of the approach largely hinges on the tightness of departure bounds, addressing this limitation through a dynamic bound adaptation process presents a promising direction for future research.

Still, one can iterate the bounding procedure multiple times to achieve more precise bounds. This can be achieved by populating $\Queue$ with newly computed earliest departure bounds. This iterative approach continues until no further modifications in the bounds are observed.

\subsubsection{Calculating conflicting sets}\label{sec:conflicting_sets}

After determining arc-specific time windows for all trips, we proceed by calculating conflicting sets $\ConflictingSet{\Arc}{i}$, i.e., sets containing pairs of trips that can incur delay by conflicting. To do so, we define for each arc $\Arc$ a set $\SetOfPairs{\Arc}$ that contains all pairs of distinct trips belonging to $\SetTripsOnArc{\Arc}$ that might potentially conflict:
\begin{align*}
    \SetOfPairs{\Arc} := \{(\Trip,\SecondTrip) \in \SetTripsOnArc{\Arc} \times \SetTripsOnArc{\Arc} \,|\, \Trip \neq \SecondTrip \wedge
    [
    \EarliestEntry{\Arc}{\Trip},
    \LatestEntry{\Arc}{\Trip}
    ] 
    \cap
    [
    \EarliestEntry{\Arc}{\SecondTrip},
    \LatestExit{\Arc}{\SecondTrip}
    ] 
    \neq \emptyset\}.
\end{align*}
By construction, $\SetOfPairs{\Arc}$ contains all tuples representing potential conflicts between trips that utilize~$\Arc$. Still, within $\SetOfPairs{\Arc}$, distinct subsets of tuples built on trips that do not appear in any other subset may emerge. 
This separation enables us to identify subsets of trips that may conflict with each other but not with other trips utilizing arc $\Arc$.
To isolate these dependencies, we split $\SetOfPairs{\Arc}$ into subsets of maximal size such that each subset contains only tuples formed from trips exclusive to that particular subset, ensuring that no trip is included in more than one subset. Technically, we partition $\SetOfPairs{\Arc}$ into $n$ sets $\SetOfPartition{\Arc}{i}$ of maximal size that satisfy the following conditions:
\begin{enumerate}
    \item[i)] The subsets $\SetOfPartition{\Arc}{i}$ are pairwise disjoint and their union constitutes the entire set $\SetOfPairs{\Arc}$, i.e., $\SetOfPairs{\Arc} := \bigcup_{i=1}^n \SetOfPartition{\Arc}{i}$ and $\SetOfPartition{\Arc}{i} \cap\SetOfPartition{\Arc}{j} = \varnothing$ for all $i \neq j$.
    \item[ii)] For any subset $\SetOfPartition{\Arc}{i}$ containing more than one pair, it holds that for each pair ${(\Trip,\SecondTrip) \in \SetOfPartition{\Arc}{i}}$, there exists another pair \( (r'',r''') \in \SetOfPartition{\Arc}{i} \) such that the intersection of their components is non-empty: \( \{\Trip,\SecondTrip\} \cap \{r'',r'''\} \neq \varnothing \). If no such pair \( (r'',r''') \) exists, then \( (\Trip,\SecondTrip) \) constitutes a singleton subset in the partition.
\end{enumerate}

With $\SetOfPartition{\Arc}{i}$ containing all tuples relative to a subset of trips, we can finally identify those trips that may incur a delay. According to Equation~\eqref{eq:travel_time_approximation}, this is the case when the number of conflicts that a trip in $\SetOfPartition{\Arc}{i}$ has, i.e., the number of trips traversing an arc at the same time, exceeds the first threshold $\ThFlow{\Arc}{1}$, beyond which it is not possible to travel on arc $\Arc$ at free-flow speed.
Then, we can tighten the maximum number of conflicts $\MaxFlow{\Arc}{\Trip}$ by counting the number of pairs in $\SetOfPartition{\Arc}{i}$ in which a trip $\Trip$ appears as the first element. Accordingly, we filter our partition of $\SetOfPartition{\Arc}{i}$ to derive conflicting sets $\ConflictingSet{\Arc}{i}$, each containing trips that may induce delay, as:
\begin{align*}
    \ConflictingSet{\Arc}{i} := \left\{(\Trip,\SecondTrip) \in \SetOfPartition{\Arc}{i} \,|\, \MaxFlow{\Arc}{\Trip} \geq \ThFlow{\Arc}{1} \right\}.
\end{align*}

\subsubsection{Multigraph representation}\label{sec:graph_processing}

\begin{figure}[!t]    
    \centering
    \subcaptionbox{Original graph $\Graph$\label{fig:multigraph_a}}{\includegraphics[width=0.3\textwidth]{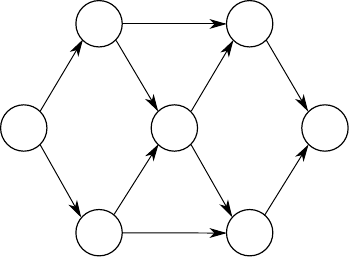}}
    \hfill
    \subcaptionbox{$\Graph$ expanded by conflicting arcs\label{fig:multigraph_b}}{\includegraphics[width=0.3\textwidth]{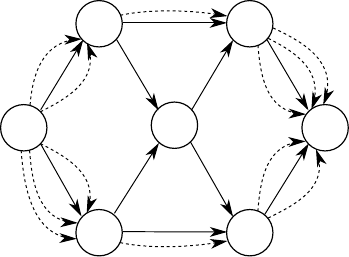}}
    \hfill
    \subcaptionbox{Non-conflicting arcs processing\label{fig:multigraph_c}}{\includegraphics[width=0.3\textwidth]{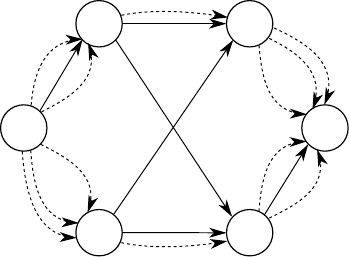}}
    \caption{Schematic representation of the graph transformation process upon identification of conflicting sets: (a) initial graph structure. (b) expansion with conflicting arcs (dashed lines). (c) {result of merging sequences} of non-conflicting arcs and removing unused arcs. The merging process involves removing the central node and replacing its incident arcs with non-conflicting arcs with the original arcs' cumulative length.}
    \label{fig:multigraph}
\end{figure}

After identifying all conflicting sets, we aim to modify our original graph $\Graph$ towards a multigraph representation. This allows us to maintain Constraints~\eqref{constr:original_alpha_constraint_first}-\eqref{constr:compute_number_of_conflicts} at a minimum. Figure~\ref{fig:multigraph} shows an example of our graph modification that consists of two steps:

%\todo{aren't we merging too many? What if one agent is trip is traveling to the center? Where has the arc from hours 9 to hours 6 gone?bit unclear.}
in the first step, we transition to a multigraph by adding a parallel arc for each conflicting set; in the second step, we sparsen our graph by removing arcs that are no longer used and merging sequences of conflict-free arcs. 
\vspace{-.5cm}
\subparagraph{Step 1:} \label{par:step_1_mg}
For each conflicting set $\ConflictingSet{\Arc}{i}$, we add a conflicting arc $\hat{\Arc}$ as a parallel arc to its original counterpart~$a$. This arc shares the characteristics of $a$, i.e., it has the same start point, endpoint, nominal travel time, and capacity. After adding a conflicting arc $\hat{\Arc}$, we modify the routes that belong to the trips contained in $\ConflictingSet{\Arc}{i}$ by replacing the original arc $a$ with the respective conflicting arc $\hat{\Arc}$, such that all trips that share a potential conflict traverse the same conflicting arc. 

Let the set of conflicting arcs be $\SetConflictingArcs\subset \SetArcs$. After Step 1, only trips that cannot induce delay on an arc $\Arc$ will still use arc $\Arc \in \SetArcs \setminus \SetConflictingArcs$. Accordingly, we will refer to trips that utilize only arcs $\Arc \in \SetArcs \setminus \SetConflictingArcs$ as non-conflicting trips and to the respective arcs as non-conflicting arcs. We remove trips from $\SetTrips$ that traverse only non-conflicting arcs.
\vspace{-.5cm}
\subparagraph{Step 2:} To sparsen our arc set, we aim to merge as many non-conflicting arcs as possible. We combine sequences of non-conflicting arcs in trip routes into single arcs, which we add to $\SetArcs$. Each new arc starts at the origin of the first arc in the sequence and ends at the destination of the last arc, and its length equals the sum of the lengths of all arcs in the sequence. %\todo{unclear if this is done correctly. If there are two trips, t1 and t2, through nodes 1-2-3, and if 1-2 is not conflicting for t1 but it is conflicting for t2, and the opposite for 2-3, it's not clear how this is handled}
 We update accordingly the routes $\TripPath{\Trip}$ in which such sequences appear.
Finally, we remove arcs from $\SetArcs$ that do not appear in any trip's route and update $\SetTripsOnArc{\Arc}$ for every arc in $\SetArcs$ accordingly.

\subsubsection{Tightened Big-M constraints}\label{sec:constr_redefinition} 
The remaining task involves updating the definition of Constraints~\eqref{constr:original_alpha_constraint_first}-\eqref{constr:original_gamma_constraint_second} as follows:
\begin{align*}
&\text{Constraints}~\eqref{constr:original_alpha_constraint_first}-\eqref{constr:original_gamma_constraint_second}, &\quad \forall \Arc \in \SetConflictingArcs, (\Trip,\SecondTrip) \in \SetOfPairs{\Arc}
\end{align*}
thereby applying these constraints exclusively to conflicting arcs, $\SetConflictingArcs$, and to each pair of trips that might conflict and potentially cause delays.

To minimize the $M_i$ constants in the remaining big-M constraints, we leverage the time windows obtained with Algorithm~\ref{algo:compute_bounds}. For any pair of trips $(\Trip, \SecondTrip)$ within $\SetOfPairs{\Arc}$, where $\Arc$ belongs to the set of conflicting arcs $\SetConflictingArcs$, we extract the following relationships from Constraints~$\eqref{constr:original_alpha_constraint_first}-\eqref{constr:original_beta_constraint_second}$: 
\begin{align*}
\BigM{1} &\geq \max(\TripDepartureOnArc{\Arc}{\Trip} - \TripDepartureOnArc{\Arc}{\SecondTrip} + \SmallConstant) = \max(\TripDepartureOnArc{\Arc}{\Trip}) - \min(\TripDepartureOnArc{\Arc}{\SecondTrip}) + \SmallConstant = \LatestEntry{\Arc}{\Trip} - \EarliestEntry{\Arc}{\SecondTrip} + \SmallConstant, &&\forall \Arc \in \SetConflictingArcs, (\Trip, \SecondTrip) \in \SetOfPairs{\Arc}\\
\BigM{2} &\geq \max(\TripDepartureOnArc{\Arc}{\SecondTrip} - \TripDepartureOnArc{\Arc}{\Trip}) = \max(\TripDepartureOnArc{\Arc}{\SecondTrip}) - \min(\TripDepartureOnArc{\Arc}{\Trip}) = \LatestEntry{\Arc}{\SecondTrip} - \EarliestEntry{\Arc}{\Trip}, &&\forall \Arc \in \SetConflictingArcs, (\Trip, \SecondTrip) \in \SetOfPairs{\Arc}\\
\BigM{3} &\geq \max(\TripDepartureOnArc{\Arc}{\SecondTrip} + \ArcNominalTravelTime{\Arc} + \DelayOnArc{\Arc}{\SecondTrip} - \TripDepartureOnArc{\Arc}{\Trip})\\
&= \max(\TripDepartureOnArc{\Arc}{\SecondTrip} + \ArcNominalTravelTime{\Arc} + \DelayOnArc{\Arc}{\SecondTrip}) - \min(\TripDepartureOnArc{\Arc}{\Trip}) = \LatestExit{\Arc}{\SecondTrip} - \EarliestEntry{\Arc}{\Trip}, &&\forall \Arc \in \SetConflictingArcs, (\Trip, \SecondTrip) \in \SetOfPairs{\Arc}\\
\BigM{4} &\geq \max(\TripDepartureOnArc{\Arc}{\Trip} + \ArcNominalTravelTime{\Arc} + \DelayOnArc{\Arc}{\Trip} - \TripDepartureOnArc{\Arc}{\SecondTrip} + \SmallConstant)\\
&= \max(\TripDepartureOnArc{\Arc}{\Trip} + \ArcNominalTravelTime{\Arc} + \DelayOnArc{\Arc}{\Trip}) - \min(\TripDepartureOnArc{\Arc}{\SecondTrip}) + \SmallConstant = \LatestExit{\Arc}{\Trip} - \EarliestEntry{\Arc}{\SecondTrip} + \SmallConstant. &&\forall \
\Arc \in \SetConflictingArcs, (\Trip, \SecondTrip) \in \SetOfPairs{\Arc}
\end{align*}

\subsection{Matheuristic}\label{subsec:matheuristic}

\begin{algorithm}[!b]
\footnotesize
\caption{Matheuristic}\label{algo:matheuristic}
\hspace{\algorithmicindent}\textbf{Input} earliest departure times $\ReleaseTimesVector$, initial lower bound $\LowerBound$, time limit $\TimeLimitMILP$
\begin{algorithmic}[1]
\State $\Solution \gets \textsc{constructSolution}(\ReleaseTimesVector)$
\State $\Solution \gets \textsc{localSearch}(\Solution)$
\While{time elapsed $\leq~\TimeLimitMILP$}
\State $\LowerBound, \Solution \gets \textsc{solveMILP}(\LowerBound, \Solution)$
\If{$\TotalDelay{\Solution} = \LowerBound$}
\State \textbf{return} $\Solution$
\EndIf
\State $\Solution \gets \textsc{localSearch}(\Solution)$
\EndWhile
\State \textbf{return} $\Solution$
\end{algorithmic}
\end{algorithm}

We develop a matheuristic incorporating our \gls{milp} into a local search scheme to find good solutions for large real-world instances. While this matheuristic is not guaranteed to find an optimal solution, it leverages the MILP from Section~\ref{subsec:milp_formulation} and can thus detect an optimal solution if encountered. Algorithm~\ref{algo:matheuristic} shows a high-level pseudo-code of our matheuristic. After creating an initial solution $\Solution$ that is agnostic to staggering, we apply a local search to improve it by staggering its trip departures. We then aim to further improve $\Solution$ by iterating between our \gls{milp} and our local search: whenever we find a new incumbent by solving the \gls{milp}, we stop and improve this incumbent with our local search. 
Once the local search terminates, we feed its final solution back to warmstart the \gls{milp}. Our matheuristic stops either when the \gls{milp} certifies optimality of $\Solution$, i.e., $\TotalDelay{\Solution} = \LowerBound$, or after a maximum time $\eta$ elapsed. Note that even if the \gls{milp} cannot certify optimality within the computational time limit, our matheuristic has the advantage of providing a lower bound to indicate the quality of $\Solution$. In the remainder of this section, we detail our construction routine to obtain an initial solution (Section~\ref{subsubsec:initial_schedule}) before we describe our local search (Section~\ref{subsubsec:local_search}).

\subsubsection{Constructing initial solutions}\label{subsubsec:initial_schedule}
To construct an initial solution $\Solution$, we set its departure time variables $\TripDepartureOnArc{\Arc}{\Trip}$, delay variables $\DelayOnArc{\Arc}{\Trip}$ and conflict counting variables $\Flow{\Arc}{\Trip}$ to infinity.
We then sort the trip-specific earliest departures $\TripReleaseTime{\Trip}$ of every trip $\Trip$ into a priority queue $\Queue$ and initialize an empty set $\SetOfArcArrivals{\Arc}$ for each arc $\Arc$ to track when each trip completes the traversal of $\Arc$. Then, we construct an initial solution iteratively: 
\begin{enumerate}
\item[i)] We start by extracting the departure time $\TripDepartureOnArc{\Arc}{\Trip}$ with the highest priority from $\Queue$. 
\item[ii)] We then set $\Flow{\Arc}{\Trip} = |\{ \Arrival \in \SetOfArcArrivals{\Arc} \,|\, \TripDepartureOnArc{\Arc}{\Trip} < \Arrival\}|$,
that is the number of trips that have not completed their traversal of arc $\Arc$ by time $\TripDepartureOnArc{\Arc}{\Trip}$, and update $\DelayOnArc{\Arc}{\Trip} = \max_{k\in K} \{0,\mu_\Arc^k\cdot(\Flow{\Arc}{\Trip}-\hat{f}_\Arc^k)\}$. 

\item[iii)] Finally, we add the time in which the arc traversal is completed, which is $\TripDepartureOnArc{\Arc}{\Trip} + \DelayOnArc{\Arc}{\Trip} + \ArcNominalTravelTime{\Arc}$, to $\SetOfArcArrivals{\Arc}$ by setting $\SetOfArcArrivals{\Arc} =  \SetOfArcArrivals{\Arc} \cup (\TripDepartureOnArc{\Arc}{\Trip} + \DelayOnArc{\Arc}{\Trip} + \ArcNominalTravelTime{\Arc})$. If $\Arc$ is not the last arc of route $\TripPath{\Trip}$, we insert $\TripDepartureOnArc{\Arc}{\Trip} + \DelayOnArc{\Arc}{\Trip} + \ArcNominalTravelTime{\Arc}$ into $\Queue$. 
\end{enumerate}
We repeat these steps until $\Queue$ is empty. 
% This procedure may generate solutions $\Solution$ for which a pair $(\Trip,\SecondTrip)$ that violates Constraints~\eqref{eq:min_distance} exists. To repair this infeasibility, we increase $\TripReleaseTime{\Trip}$ by $\SmallConstant$ and recompute the solution until the infeasibility is resolved.

Note that once we obtained a feasible solution, it is possible to compute values for the \gls{milp}'s binary variables $\VarAlpha{\Arc}{\Trip}{\SecondTrip},\VarBeta{\Arc}{\Trip}{\SecondTrip}$ and $\VarGamma{\Arc}{\Trip}{\SecondTrip}$ with a support function checking their activation conditions as introduced in Section~\ref{subsec:milp_formulation}.

\begin{figure}[b]    
    \centering
    \begin{subfigure}[b]{0.4\textwidth}
        \def\svgwidth{0.9\linewidth}
        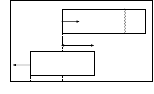
        \caption{Conflict metrics}
    \end{subfigure}
    \hfill
    \begin{subfigure}[b]{0.4\textwidth}
        \def\svgwidth{0.9\linewidth}
        %% Creator: Inkscape 1.1.1 (3bf5ae0d25, 2021-09-20), www.inkscape.org
%% Accompanies image file 'metrics_after.pdf' (pdf, eps, ps)
\begingroup%
  \makeatletter%
  \providecommand\color[2][]{%
    \errmessage{(Inkscape) Color is used for the text in Inkscape, but the package 'color.sty' is not loaded}%
    \renewcommand\color[2][]{}%
  }%
  \providecommand\transparent[1]{%
    \errmessage{(Inkscape) Transparency is used (non-zero) for the text in Inkscape, but the package 'transparent.sty' is not loaded}%
    \renewcommand\transparent[1]{}%
  }%
  \providecommand\rotatebox[2]{#2}%
  \newcommand*\fsize{\dimexpr\f@size pt\relax}%
  \newcommand*\lineheight[1]{\fontsize{\fsize}{#1\fsize}\selectfont}%
  \ifx\svgwidth\undefined%
    \setlength{\unitlength}{73.61256697bp}%
    \ifx\svgscale\undefined%
      \relax%
    \else%
      \setlength{\unitlength}{\unitlength * \real{\svgscale}}%
    \fi%
  \else%
    \setlength{\unitlength}{\svgwidth}%
  \fi%
  \global\let\svgwidth\undefined%
  \global\let\svgscale\undefined%
  \makeatother%
  \begin{picture}(1,0.6118412)%
    \lineheight{1}%
    \setlength\tabcolsep{0pt}%
    \put(0,0){\includegraphics[width=\unitlength,page=1]{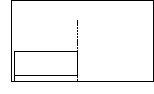}}%
    \put(0.27326415,0.18423771){\makebox(0,0)[lt]{\lineheight{1.25}\smash{\begin{tabular}[t]{l}$\tau_a^{\text{N}}$\end{tabular}}}}%
    \put(0,0){\includegraphics[width=\unitlength,page=2]{metrics_after.pdf}}%
    \put(0.6869586,0.45885437){\makebox(0,0)[lt]{\lineheight{1.25}\smash{\begin{tabular}[t]{l}$\tau_a^{\text{N}}$\end{tabular}}}}%
    \put(0.08492794,0.01465923){\makebox(0,0)[lt]{\lineheight{1.25}\smash{\begin{tabular}[t]{l}$s_a^r - \phi^{\text{f}}$\end{tabular}}}}%
    \put(0.50058592,0.01467803){\makebox(0,0)[lt]{\lineheight{1.25}\smash{\begin{tabular}[t]{l}$s_a^{r'}+ \rho^{\text{f}}$\end{tabular}}}}%
    \put(-0.00358189,0.1666773){\makebox(0,0)[lt]{\lineheight{1.25}\smash{\begin{tabular}[t]{l}$r$\end{tabular}}}}%
    \put(-0.00358189,0.44992782){\makebox(0,0)[lt]{\lineheight{1.25}\smash{\begin{tabular}[t]{l}$r'$\end{tabular}}}}%
  \end{picture}%
\endgroup%

        \caption{Backward and forward shifts}
    \end{subfigure}
    \caption{Figure (a) illustrates the quantities computed in the local search for resolving a conflict $(\Trip,\SecondTrip,\Arc)$, where $\SecondTrip$ incurs delay $\DelayOnArc{\Arc}{\Trip}$. It illustrates the time overlap $\Overlap$, the feasible backward shift $\BackwardShift{f}$, and forward shift $\ForwardShift{f}$. Figure (b) depicts trips $\Trip$ and $\SecondTrip$ post the application of these shifts, resulting in conflict resolution and the consequent elimination of the arc delay.}
    \label{fig:ls_metrics}
\end{figure}

\subsubsection{Local search}\label{subsubsec:local_search}

Our local search aims to iteratively remove conflicts based on a priority queue, aiming to first resolve conflicts that induce the largest delays. Let the tuple $(\Trip,\SecondTrip,\Arc)$ index the conflict that $\Trip$ has with $\SecondTrip$ on arc $\Arc$. Given a conflict $(\Trip,\SecondTrip,\Arc)$, we resolve it by either shifting $\Trip$ backward or $\SecondTrip$ forward in time without violating the trips' time windows. In this context, we compute additional quantities for each conflict, utilizing the relations shown in Figure~\ref{fig:ls_metrics}. Specifically, consider a conflict $(\Trip,\SecondTrip,\Arc)$ with $\TripDepartureOnArc{\Arc}{\Trip} < \TripDepartureOnArc{\Arc}{\SecondTrip}$. We then compute the time overlap $\Overlap$ between $\Trip$ and $\SecondTrip$, which denotes how much staggering we need to apply to resolve the conflict. Formally:
$$\Overlap = \TripDepartureOnArc{\Arc}{\Trip} + \ArcNominalTravelTime{\Arc} + \DelayOnArc{\Arc}{\Trip} - \TripDepartureOnArc{\Arc}{\SecondTrip} + \SmallConstant.$$
Clearly, $\Overlap = \BackwardShift{} = \ForwardShift{} $ indicates the required backward shift $\BackwardShift{}$ for $\Trip$ and the required forward shift $\ForwardShift{}$ for $\SecondTrip$ to resolve the conflict by staggering just one of the two trips.
However, shifting a trip by $\Overlap$ may cause time window violations. Therefore, we limit $\BackwardShift{}$ and $\ForwardShift{}$ to their feasible counterparts, $\BackwardShift{f}$ and $\ForwardShift{f}$. These values represent the maximum backward and forward time shifts applicable without violating the time window of each trip.
Formally, $\BackwardShift{f}$ reads:
$$\BackwardShift{f} = \min(\TripDepartureOnArc{}{\Trip} - \TripReleaseTime{\Trip}, \Overlap),$$
while $\ForwardShift{f}$ reads:
$$\ForwardShift{f} = \min(\MaxStaggering{\SecondTrip} - (\TripDepartureOnArc{}{\SecondTrip} - \TripReleaseTime{\SecondTrip}), \Overlap),$$
with $\TripDepartureOnArc{}{\Trip}$ denoting the starting time of a trip $\Trip$.

\begin{algorithm}[!b]
\footnotesize
\caption{Local Search}\label{alg:local_search}
\hspace{\algorithmicindent}\textbf{Input:} solution $\Solution$
\begin{algorithmic}[1]
\State $\ConflictsList \gets \textsc{identifyConflicts}(\Solution)$ \label{line:get_conflicts_1}
\While{$\ConflictsList$ is not empty} \label{line:empty_q}
    \State $(\Trip, \SecondTrip, \Arc) \gets \ConflictsList.\textsc{pop}()$ \label{line:get_single_conflict}
    \State $\ForwardShift{f}, \BackwardShift{f},\Overlap \gets \textsc{analyzeConflict}( \Trip, \SecondTrip,\Arc)$ \label{line:analyze_conflict}
    \If {$\ForwardShift{f} + \BackwardShift{f} \geq \Overlap$} \label{line:check_compensation}
        \State $\bar{\Solution} \gets \textsc{resolveConflict}(\Solution, \Trip, \SecondTrip, \ForwardShift{f}, \Overlap)$\label{line:resolve_conflict}
        \If {$\TotalDelay{\bar{\Solution}} < \TotalDelay{\Solution}$}\label{line:check_delay_improvement}
            \State $\Solution \gets \bar{\Solution}$ \label{line:update_solution}
            \State $\ConflictsList \gets \textsc{identifyConflicts}(\bar{\Solution})$\label{line:update_priority_queue}
        \EndIf
    \EndIf
\EndWhile
\State \textbf{return} $\Solution$
\end{algorithmic}
\end{algorithm}

Based on these quantities, Algorithm~\ref{alg:local_search} provides the pseudocode for our local search: given an initial solution $\Solution$, we sort conflicts into a priority queue $\ConflictsList$, thus ordering conflicts based on their induced delay in decreasing order (l.\ref{line:get_conflicts_1}). We then process $\ConflictsList$ until it is empty~(l.\ref{line:empty_q}). In each iteration, we pop the first conflict from $\ConflictsList$~(l. \ref{line:get_single_conflict}) and compute $\ForwardShift{f}$, $\BackwardShift{f}$, and $\Overlap$ as defined above~(l. \ref{line:analyze_conflict}). Only if the sum of $\ForwardShift{f}$ and $\BackwardShift{f}$ allows to compensate $\Overlap$ (l.\ref{line:check_compensation}), we resolve the conflict~(l.\ref{line:resolve_conflict}). If this leads to an improved solution (l.\ref{line:check_delay_improvement}), we update our solution (l.\ref{line:update_solution}) as well as our priority queue (l.\ref{line:update_priority_queue}) since new conflicts may arise. 
If the conflict cannot be resolved, we directly proceed with the next item in $\ConflictsList$ without updating $\Solution$ or $\ConflictsList$. In the remainder of this section, we elaborate on constructing and evaluating a neighborhood to create a new solution that resolves a conflict.

While our local search follows a rather simple greedy repair selection based on conflict severity, its complexity lies in efficiently evaluating a neighborhood, i.e., computing the impact of resolving one conflict on all other trips' departure times. In the best case, resolving a conflict induces changes only in arc entry times for the conflicting trips, while in the worst case, resolving a conflict can lead to a very large cascade of changes in other trip departures. Algorithm~\ref{algo:resolve_conflict} shows how we efficiently update our solution when resolving a conflict: we first apply a minimum staggering to the departures of $\Trip$ and $\SecondTrip$. When staggering a trip departure time to resolve a conflict, we insert a corresponding tuple $(\Trip,\Arc)$, consisting of the trip's index $\Trip$ and its route's first arc $\Arc$ into a priority queue $\Queue$, which we use to maintain trips that require recalculating arc entry times due to the applied staggering. We sort the tuples in $\Queue$ in ascending order based on their start times~$\TripDepartureOnArc{\Arc}{\Trip}$~(l.\ref{line:stagger_departure}). After resolving the initial conflict, we start processing trips in $\Queue$ (l.\ref{line:processing}) by checking whether the trip's updated departure time $\TripDepartureOnArc{\Arc}{\Trip}$ creates a change in the arc entry times of another trip $\SecondTrip$ (l.\ref{line:check_new_conflicts}). In this case, we insert $\SecondTrip$ into $\Queue$. If inserting $\SecondTrip$ into $\Queue$ invalidates a further propagation of arc entry times for trip $\Trip$ -- i.e., if $\TripDepartureOnArc{\Arc}{\SecondTrip} < \TripDepartureOnArc{\Arc}{\Trip}$, which disrupts the order of our priority queue -- we reinsert $\TripDepartureOnArc{\Arc}{\Trip}$ into $\Queue$ and restart evaluating $\Queue$~(l.~\ref{line:trip_flag_is_false}\&\ref{line:reinsert_new_conflict}). If $\TripDepartureOnArc{\Arc}{\Trip}$ does not induce changes in other trips, we recursively propagate the change in trip $\Trip$'s arc entry times on subsequent arcs until a change arises~(l.\ref{line:stag_dep_else}--l.\ref{line:insert_new_trips_2}). 

\begin{algorithm}[!t]
\footnotesize
\caption{\textsc{resolveConflict}}\label{algo:resolve_conflict}
\hspace{\algorithmicindent}\textbf{Input:} solution $\Solution$, trips $\Trip$ and $\SecondTrip$, feasible forward shift $\ForwardShift{f}$, time overlap $\Overlap$ 
\begin{algorithmic}[1]
\State $\Queue \gets \textsc{staggerDepartures}(\Trip,\SecondTrip,\ForwardShift{f},\Overlap)$ \label{line:stagger_departure}
\While{$\Queue$ is not empty} 
    \State $(\Trip,\Arc) \gets \Queue.\textsc{pop}()$ \label{line:processing}
    \State $\SecondTrip, \DelayOnArc{\Arc}{\Trip} \gets \textsc{updateSolution}(\Trip,\Arc)$ \label{line:check_new_conflicts}
    \If{$\SecondTrip \neq \textsc{false}$}\label{line:trip_flag_is_false}
        \State $\Queue \gets \textsc{insertTrips}(\Trip,\SecondTrip,\Arc)$ \label{line:reinsert_new_conflict}
    \Else \label{line:stag_dep_else}
        \While{$\SecondTrip = \textsc{false}$}
            \State $\TripDepartureOnArc{\Successor{\Arc}{\Trip}}{\Trip} \gets \TripDepartureOnArc{\Arc}{\Trip} + \ArcNominalTravelTime{\Arc} + \DelayOnArc{\Arc}{\Trip}$
            \State $\Arc \gets \Successor{\Arc}{\Trip}$
            \State $\SecondTrip,\DelayOnArc{\Arc}{\Trip} \gets \textsc{updateSolution}(\Trip,\Arc)$
            \If{$\SecondTrip \neq \textsc{false}$}
                \State $\Queue \gets \textsc{insertTrips}(\Trip,\SecondTrip,\Arc)$ \label{line:insert_new_trips_2} 
            \EndIf            
        \EndWhile
    \EndIf
\EndWhile
\State \textbf{return} $\Solution$
\end{algorithmic}
\end{algorithm}

Our overall update ends once $\Queue$ is empty and all arc entry times have been correctly adjusted. To conclude this section, we detail how we stagger departure times, check for solution changes, and maintain the priority queue.

$\textsc{staggerDepartures}(\Trip,\SecondTrip,\ForwardShift{f},\Overlap):$ To resolve a conflict, we need to shift either $r$ or $r'$ or both in time by $\Overlap$. Clearly, it is desirable to shift only one of these trips in time to keep the number of newly arising conflicts and entry time changes at a minimum. Accordingly, we proceed as follows to resolve a conflict: 
\begin{enumerate}
\item If $\ForwardShift{f} \geq \Overlap$ we adjust $\TripDepartureOnArc{\Arc}{\Trip} \gets \TripDepartureOnArc{\Arc}{\Trip} +  \Overlap$ and push $(\Trip,\Arc)$ into $\Queue$.
\item If a conflict cannot be resolved by shifting only trip $r$ in time, i.e., $\ForwardShift{f} < \Overlap$, we stagger both trips such that $\TripDepartureOnArc{\Arc}{\Trip} \gets \TripDepartureOnArc{\Arc}{\Trip} + \ForwardShift{f}$ and $\TripDepartureOnArc{\Arc}{\SecondTrip} \gets \TripDepartureOnArc{\Arc}{\SecondTrip}- (\Overlap - \ForwardShift{f})$ and push $(\Trip,\Arc)$ as well as $(\SecondTrip,\Arc)$ into $\Queue$.
\end{enumerate}

$\textsc{updateSolution}(\Trip,\Arc)$: 
After staggering trips to resolve a conflict, we need to update the arc entry times of these trips as well as the arc entry times of all other trips that might be affected by the change. To do so, we first compute the number of conflicts $\Flow{\Arc}{\Trip}$ that arise for the respective trip and arc as:
$${\Flow{\Arc}{\Trip}} = \sum_{\SecondTrip \in \SetTripsOnArc{\Arc} \setminus \{\Trip\}}\mathbb{I}({\TripDepartureOnArc{\Arc}{\SecondTrip}} \leq {\TripDepartureOnArc{\Arc}{\Trip}} < {\TripDepartureOnArc{\Arc}{\SecondTrip} + \ArcNominalTravelTime{\Arc} + \DelayOnArc{\Arc}{\SecondTrip}}),$$
which subsequently allows us to update the corresponding delay as:
$${\DelayOnArc{\Arc}{\Trip}} = \max_{k\in K} \{0,\mu_\Arc^k\cdot(\Flow{\Arc}{\Trip}-\hat{f}_\Arc^k)\}.$$
We can then check whether a new conflict arises or an existing conflict gets resolved between $r$ and some other trip $r'$ due to the updated entry time of $r$.
%Now, we check if the modified arc entry time and delay introduce changes in the arc entry times of other trips $\SecondTrip$. This occurs in one of these two cases:
%\begin{itemize}
%    \item [Case 1:] A new conflict arises between trips $\Trip$ and $\SecondTrip$ on arc $\Arc$;
%    \item [Case 2:] An existing conflict between trips $\Trip$ and $\SecondTrip$ on arc $\Arc$ resolves.
%\end{itemize}
If either is the case, we return the index of trip $\SecondTrip$; otherwise, we return a $\textsc{false}$ flag to indicate that no new conflict was induced or resolved by staggering $r$.

$\textsc{insertTrips}(\Trip,\SecondTrip,\Arc)$: 
If we detect a new conflict between $r$ and $r'$, we add $(\SecondTrip,\Arc)$ to $\Queue$. Here, we need to account for the following two cases: if $\TripDepartureOnArc{\Arc}{\Trip} < \TripDepartureOnArc{\Arc}{\SecondTrip}$, changes induced in $r'$ do not affect $r$. It suffices to push $(\SecondTrip,\Arc)$ into $\Queue$. However, if $\TripDepartureOnArc{\Arc}{\SecondTrip} < \TripDepartureOnArc{\Arc}{\Trip}$, changes in $r'$ affect $r$, which invalidates preceding updates made to $r$. In this case, we additionally push $(\Trip,\Arc)$ back to $\Queue$ to ensure valid arc entry time updates upon the algorithm's termination.

\subsection{Online algorithm}\label{subsec:online_algorithm}

\begin{algorithm}[b]
\footnotesize
\caption{Online algorithm}\label{algo:online_algorithm}
\hspace{\algorithmicindent}\textbf{Input} Set of trips $\SetTrips$, number of epochs $\NumberOfEpochs$, time limit $\TimeLimitMILP$
\begin{algorithmic}[1]
\State $\TransferTrips \gets \varnothing;\, \DummyTrips \gets \varnothing$
\For{$ n \in \{1, \ldots ,\NumberOfEpochs\}$} \label{line:online_iter}
\State $\SetTripsEpoch{n} \gets \textsc{getTripsEpoch}(\SetTrips)$\label{line:online_get_trips_epoch}
\State $\ReleaseTimesVector^n,\LowerBound \gets \textsc{preprocess}(\SetTripsEpoch{n}, \TransferTrips, \DummyTrips)$ \label{line:online_preprocess}
\State $\Solution^{n} \gets \textsc{matheuristic}(\ReleaseTimesVector^n , \LowerBound, \TimeLimitMILP)$ \label{line:online_math}
\If{$n < N$}
\State $\TransferTrips \gets \textsc{getTransferTrips}(\Solution^{n})$ \label{line:online_transfer}
\State $\DummyTrips \gets \textsc{getDummyTrips}(\Solution^{n},\TransferTrips)$ \label{line:online_artificial}
\EndIf
\EndFor
% \State $\Solution \get \textsc{getCompleteSolution}(\textbf{\})$
\State $\boldsymbol{\Solution} \gets (\Solution^n)_{n \in \{1, \ldots, N\}}$  \Comment{Vector of epoch solutions}
\State $\Solution \gets \textsc{getSolution}(\boldsymbol{\Solution})$ \Comment{Solution for every trip in $\SetTrips$}
\State \textbf{return} $\Solution$
\vspace{0.4em}
\end{algorithmic}
\end{algorithm}

While we discussed our matheuristic for a static problem setting, we can straightforwardly apply it to an online problem setting where trips enter the system over time, e.g., following a point process. In this context, we focus on batched decision-making, which allows us to apply our matheuristic in a rolling-horizon fashion.

Specifically, we divide the static time horizon into $n \in \{1, \ldots,\NumberOfEpochs\}$ distinct epochs, each spanning a fixed time interval of size $\EpochLength$.
Before a new epoch $n$ starts, we then take a staggering decision for all requests that entered the system after the last decision taken, collected in a set $\SetTripsEpoch{n}$,
formally:
\begin{align*}
    \SetTripsEpoch{n} := \{\Trip \in \SetTrips \,|\, (n-1) \cdot \EpochLength \leq \TripReleaseTime{\Trip} < n \cdot \EpochLength\}. &&\forall n \in \{1,\ldots,N\} 
\end{align*}
We can then compute a solution $\Solution^n$ for the respective batch of requests by running our matheuristic with a suitable time limit $\eta$.

Algorithm~\ref{algo:online_algorithm} shows how we apply our matheuristic in such a rolling horizon setting: we iterate over all epochs (l.\ref{line:online_iter}), each time focusing on the respective request set $\SetTripsEpoch{n}$ (l.\ref{line:online_get_trips_epoch}), and calculating the relative earliest departure times vector $\textbf{e}^n$, and a lower bound for epoch $n$ as described in Sections~\ref{sec:arc_dependent_time_windows}-\ref{sec:graph_processing} (l.\ref{line:online_preprocess}). 

For each epoch, we compute a solution $\Solution^n$ by running our matheuristic (l.\ref{line:online_math}). Afterward, we can fix this solution for all trips in $\SetTripsEpoch{n}$. Still, we need to identify trips that span multiple epochs to identify all conflicts in future epochs during preprocessing. To do so, we collect such trips in a set $\TransferTrips$, formally:
\begin{align*}
\TransferTrips := \{\Trip \in \SetTripsEpoch{n} \,|\, \TripDepartureOnArc{\Arc'}{\Trip} > n \cdot \EpochLength\},
\end{align*}
with $\Arc'$ being the last arc of route $\TripPath{\Trip}$. Note that after the respective adjustments, trips in $\TransferTrips$ can also conflict with trips being in $\SetTripsEpoch{n}$ but not in $\TransferTrips$. Maintaining those trips by incorporating them directly in $\TransferTrips$ poses a computational overhead in the matheuristic, as it necessitates the inclusion of additional conflicting trips to accurately calculate their arc entry times.

To mitigate this overhead, we identify all conflicts $(\Trip, \SecondTrip, \Arc)$ in $\Solution^{n}$, where $\Trip$ is part of $\TransferTrips$, $\SecondTrip$ is not, and $\Arc$ is included in $\TripPath{\Trip}$. 
For each identified conflict, we generate an artificial trip $\SecondTrip$ with $\TripPath{\SecondTrip}$ consisting solely of the arc $\Arc$, $\TripReleaseTime{\SecondTrip}$ being the arc entry time of the conflicting trip $\TripDepartureOnArc{\Arc}{\SecondTrip}$, $\TripDeadline{\SecondTrip}$ matching the time in which the conflicting trip completes the arc traversal, and $\MaxStaggering{\SecondTrip}$ being zero to preclude any staggering time decisions by the algorithm. We maintain these dummy trips in a separate set $\DummyTrips$ and take them into account in future epochs.
\section{Design of Experiments}\label{sec:case_study} 
We used Python 3.9.10 and Gurobi 10.0.1 to build and solve the \gls{milp} model. To implement our matheuristic, we used C\texttt{++} to implement the local search and integrated it into the Python framework via pybind11 \citep{pybind11}. We performed all experiments on an Intel(R) i9-9900 CPU, 3.1 GHz, with 56 GB of RAM.

To create a realistic case study, we use the New York taxi data set \citep{NTLC2015} and focus on the area of Manhattan in New York City.
%Note: I've used 4 CPUs, to which I've allocated 14GB of memory
 \vspace{-.5cm}
\subparagraph{Road network.}
We extract the Manhattan road network from the OpenStreetMap dataset to obtain a network with 7782 arcs and 3213 nodes \citep{OpenStreetMap}. Studying staggered routing on this raw network remains challenging as it requires checking conflicts for a large number of arcs. 
% For a given arc $\Arc$, the worst-case scenario occurs when every trip in $\SetTripsOnArc{\Arc}$ potentially conflicts with every other trip using the same arc. In this situation, the number of constraints necessary to verify conflicts between pairs of trips grows as $\BigO{\frac{|\SetTripsOnArc{\Arc}|!}{(|\SetTripsOnArc{\Arc}|-2)!}}$, which represents the number of ways to permute two elements within the set $\SetTripsOnArc{\Arc}$.
To decrease the number of arcs while preserving the length of the shortest paths within the network, we use contraction hierarchies \citep[cf.][]{GeisbergerSandersEtAl2008}, a well-known reduction technique for online routing, to contract the road network graph. We employ iterative node contractions that merge arcs while preserving the shortest routes between the remaining nodes.
The efficiency of this contraction hinges on the order in which nodes are processed. We first contract nodes with the smallest \textit{edge difference} - a metric capturing the difference between the number of shortcuts added at the time of node contraction and the number of incident arcs. Here, we adopt the \textit{lazy update} approach of \citet{GeisbergerSandersEtAl2008} to reduce edge difference recalculations and enhance efficiency.

%By integrating shortcuts into the \gls{osm} network, trips can be routed on routes made of fewer but lengthier arcs, reducing the complexity of our instances. However, adding shortcuts inevitably alters the spatial distribution of conflicts. To better approximate the spatial distribution of conflicts intrinsic to the OSM road network, we cap the shortcut length for trip routing at 2 kilometers. We further discuss the implications of this approximation at the end of this section.

\vspace{-.5cm}
\subparagraph{Trip data.} We obtain trips $\SetTrips$ from the NYC taxi dataset, collected in January 2015 \citep{NTLC2015}, which contains taxi trips' origins, destinations, and departure times.
We determine trip routes by selecting the shortest route connecting the nearest network nodes to the origin and destination with the minimum number of arcs. Trip start times $\TripReleaseTime{\Trip}$ correspond to the taxi ride's start times provided in the dataset.
We assign trip deadlines $\TripDeadline{\Trip}$ as the time of arrival at the destination in the uncontrolled scenario plus a relative quota equal to 25\% of the route's nominal travel time plus a fixed tolerance of 30 seconds for each trip. We calculate the maximum allowable departure time shift $\MaxStaggering{\Trip}$ as a percentage $\TimePercentage$ of the total nominal travel time required to traverse the route $\TripPath{\Trip}$. Unless specified differently, we set this percentage to 10\%.
\vspace{-.5cm}
\subparagraph{Parameterization.} We set the nominal travel times $\ArcNominalTravelTime{\Arc}$ assuming vehicles traveling at 20 kilometers per hour, which approximates the average traffic speed recorded in New York City downtown in 2022 \citep{Pishue2023}. As we focus on the impact and control of a large \gls{amod} fleet, we assume a baseload of traffic that is out of the operator's control and adjust arc capacities $\ArcCapacity{\Arc}$ accordingly to account for potential congestion induced by the \gls{amod} fleet. Specifically, we consider two scenarios: first, a \gls{lc} scenario, which reflects a setting in which the exogenous traffic in the system is uncongested, and the AMoD fleet may induce congestion in case of inefficient operations. Second, a \gls{hc} scenario in which the exogenous traffic already encounters congestion that might be amplified by the AMoD fleet operations. To create these scenarios, we calibrate the residual capacities to accommodate one \gls{amod} trip every 15 seconds for the \gls{lc} scenario and 30 seconds for the \gls{hc} scenario. This is done by setting, respectively, arc capacities as $\ArcCapacity{\Arc} = \ArcNominalTravelTime{\Arc} / 15$ and $\ArcCapacity{\Arc} = \ArcNominalTravelTime{\Arc} / 30$ and rounding the resulting value to its nearest integer.
To calibrate the arc travel time function, we considered a single non-flat piece, a maximum flow at which it is possible to travel the arc under free-flow conditions equal to $\ThFlow{\Arc}{} = \ArcCapacity{\Arc}$, and a slope equal to $\PWLSlope{\Arc}{} = \frac{0.5 \cdot \ArcNominalTravelTime{\Arc}}{\ArcCapacity{\Arc}}$.

\begin{table}[!b]
\caption{Minimum, maximum, and average values of the total delay in the uncontrolled solution $\TotalDelayUncSolution$, the number of conflicting arcs $|\SetConflictingArcs|$, and the maximum number of trips traversing the same conflicting arc $|\SetTripsOnArc{\hat{\Arc}}|^{\text{MAX}}$ for low congestion (LC) and high congestion (HC) instances.}
\label{table:overview_instances}
\centering
\begin{tabularx}{\textwidth}{l|XXX|XXX}
\toprule
 & \multicolumn{3}{c}{LC} & \multicolumn{3}{c}{HC} \\
\cmidrule(lr){2-4} \cmidrule(lr){5-7}
 & Min & Max & Avg & Min & Max & Avg \\
\midrule
$\TotalDelayUncSolution$ [min]             & 54   & 133  &  82  & 348  & 1005 & 577 \\
$|\hat{\mathcal{A}}|$                      & 1400 & 2489 & 1971 & 2496 & 3660 & 3165 \\
$|\SetTripsOnArc{\hat{\Arc}}|^{\text{MAX}}$ & 52   & 161  & 102  & 76   & 161  & 110 \\
\bottomrule
\end{tabularx}
\end{table}
\vspace{-.5cm}
\subparagraph{Instances.}
To create our instances, we sampled 5000 trips between four and five p.m. over 31 days, sorting them by start time.
To gauge their complexity, we compute the uncontrolled solution for each instance both in the \gls{lc} and \gls{hc} scenarios and calculate its corresponding total delay $\TotalDelayUncSolution$, the number of conflicting arcs $|\SetConflictingArcs|$, and the maximum number of trips traversing the same conflicting arc $|\SetTripsOnArc{\hat{\Arc}}|^{\text{MAX}}$. Table~\ref{table:overview_instances} summarizes the characteristics of these instances, detailing their minimum, maximum, and average values. The total fleet delay in the \gls{hc} scenario is roughly ten times greater than in the \gls{lc} scenario. Additionally, there is a consistent rise in the number of conflicting arcs between scenarios, with an increase of approximately 1000 arcs. However, the difference in the peak number of trips per conflicting arc is relatively modest between scenarios, averaging an increase of eight trips.

\section{Results}\label{sec:results}
In the following, we discuss the findings of our numerical studies. First, we analyze our matheuristic performance for both an idealized offline setting as well as a myopic online setting in Section~\ref{sec:algo_performance}. We then proceed with a fine-grained analysis in Section~\ref{sec:structural_analysis} before discussing the trade-off between congestion-related delay reduction and trip departure time shifts in Section~\ref{sec:sensitivity_analyses}.

\begin{figure}[!b]
    \centering
    % First row of subfigures
    \begin{subfigure}{\textwidth}
        \centering
        % This file was created with tikzplotlib v0.10.1.
\begin{tikzpicture}

\definecolor{darkslategray38}{RGB}{38,38,38}
\definecolor{lightgray204}{RGB}{204,204,204}

\begin{axis}[
height=\DelayReductionHeight,
width=\DelayReductionWidth,
yticklabels={}, % Hide y tick labels
ytick style={draw=none}, % 
axis line style={darkslategray38},
legend style={
  fill opacity=0.8,
  draw opacity=1,
  text opacity=1,
  at={(0.03,0.97)},
  anchor=north west,
  draw=none
},
tick align=outside,
tick pos=left,
x grid style={lightgray204},
xlabel=\textcolor{darkslategray38}{\(\displaystyle \Theta[\%]\)},
xmajorgrids,
xmin=-1, xmax=101,
xtick style={color=darkslategray38},
xtick={0,20,40,60,80,100},
xticklabels={
  \(\displaystyle {0}\),
  \(\displaystyle {20}\),
  \(\displaystyle {40}\),
  \(\displaystyle {60}\),
  \(\displaystyle {80}\),
  \(\displaystyle {100}\)
},
y dir=reverse,
y grid style={lightgray204},
ymin=-0.5, ymax=1.5,
ytick style={color=darkslategray38},
ytick={0,1},
yticklabels={OFF,ON}
]
\path [draw=black, fill=white]
(axis cs:85.7673943969413,-0.4)
--(axis cs:85.7673943969413,0.4)
--(axis cs:100,0.4)
--(axis cs:100,-0.4)
--(axis cs:85.7673943969413,-0.4)
--cycle;
\path [draw=black, fill=white]
(axis cs:89.9678711302725,0.6)
--(axis cs:89.9678711302725,1.4)
--(axis cs:91.8862645426898,1.4)
--(axis cs:91.8862645426898,0.6)
--(axis cs:89.9678711302725,0.6)
--cycle;
\addplot [black, forget plot]
table {%
85.7673943969413 0
80.1721875412009 0
};
\addplot [black, forget plot]
table {%
100 0
100 0
};
\addplot [black, forget plot]
table {%
80.1721875412009 -0.2
80.1721875412009 0.2
};
\addplot [black, forget plot]
table {%
100 -0.2
100 0.2
};
\addplot [black, forget plot]
table {%
89.9678711302725 1
87.3923931840809 1
};
\addplot [black, forget plot]
table {%
91.8862645426898 1
94.3145425165934 1
};
\addplot [black, forget plot]
table {%
87.3923931840809 0.8
87.3923931840809 1.2
};
\addplot [black, forget plot]
table {%
94.3145425165934 0.8
94.3145425165934 1.2
};
\addplot [black, forget plot]
table {%
99.6306369376103 -0.4
99.6306369376103 0.4
};
\addplot [black, forget plot]
table {%
90.9335974319942 0.6
90.9335974319942 1.4
};
\addplot [draw=black, fill=white, forget plot, mark=*, only marks]
table{%
x  y
85.9501357759841 -0.229903242967598
85.777862670751 -0.257157573668744
100 0.127690936634137
100 -0.270302351287626
85.6135553241588 -0.299183490994118
99.8289892586735 0.00824509372820353
99.8305621132975 0.26693548946663
99.8172382545536 0.0543641539766408
87.2709125623164 0.169609775580074
96.3537321478898 -0.0816426895632003
99.4380246249146 -0.290183199806861
85.7569261231317 -0.231753505856298
100 -0.204948148992856
83.5171191228929 0.258491915493171
84.2611553729852 0.111859654017154
99.8530813657235 -0.142154288952707
100 0.162129845937416
100 -0.00736017437762188
100 -0.229768909182323
100 0.0108074242263506
86.3311479526958 -0.119419883122961
99.6306369376103 -0.238809457841289
98.9768023301197 0.0548098323867848
100 -0.169898740783202
84.7785665969419 -0.0077549487718287
84.3265490387304 0.234474712841051
80.1721875412009 0.1471836027463
100 0.262983407076413
99.8293938099938 0.0663130002943371
85.1654571438188 0.0351109648977892
99.8180700689071 -0.00715680736296748
82.4649349385438 0.931865688756932
92.1429887087138 0.93878800068193
89.5204419078277 0.712049664831507
90.9335974319942 0.857067691430218
87.3923931840809 0.794745633954409
89.9349721583992 1.14703426051976
91.6295403766658 1.29431417062016
93.7714109687106 1.0402133610088
90.3908180762624 1.17088636085357
91.3854940461695 0.855181606396342
92.6422615021059 0.826482443303831
91.0479239905113 0.723387510482467
90.7466008202 1.28817749218449
91.4700614795043 0.853547149332815
84.8971267286748 1.06605834341164
84.9870282954121 1.05673557457333
94.3145425165934 0.923316678003352
90.3504611722708 1.25635953589919
94.015612449064 1.07930281473784
93.4436884305983 1.20406880906535
90.8391413336395 1.0964405251101
91.2183622588506 1.14452001332879
91.1244315364511 1.21158480941391
90.8558233613643 1.03708877734473
93.2017900339511 1.27500266740493
87.7902579322175 0.743274526693473
88.784582628964 1.21731636704954
91.0847178754412 0.921824321981312
90.389479524641 0.884040053417282
90.0007701021458 1.2868365328272
93.1657876145315 0.891451524181201
};
\end{axis}

\end{tikzpicture}
        % This file was created with tikzplotlib v0.10.1.
\begin{tikzpicture}

\definecolor{darkslategray38}{RGB}{38,38,38}
\definecolor{lightgray204}{RGB}{204,204,204}

\begin{axis}[
height=\DelayReductionHeight,
width=\DelayReductionWidth,
yticklabels={}, % Hide y tick labels
ytick style={draw=none}, % 
axis line style={darkslategray38},
legend style={fill opacity=0.8, draw opacity=1, text opacity=1, draw=none},
tick align=outside,
tick pos=left,
x grid style={lightgray204},
xlabel=\textcolor{darkslategray38}{\(\displaystyle \Theta[\mathrm{min}]\)},
xmajorgrids,
xmin=46.7890617603938, xmax=120.197213443846,
xtick style={color=darkslategray38},
xtick={40,60,80,100,120,140},
xticklabels={
  \(\displaystyle {40}\),
  \(\displaystyle {60}\),
  \(\displaystyle {80}\),
  \(\displaystyle {100}\),
  \(\displaystyle {120}\),
  \(\displaystyle {140}\)
},
y dir=reverse,
y grid style={lightgray204},
ymin=-0.5, ymax=1.5,
ytick style={color=darkslategray38},
ytick={0,1},
yticklabels={OFF,ON}
]
\path [draw=black, fill=white]
(axis cs:67.0780217069725,-0.4)
--(axis cs:67.0780217069725,0.4)
--(axis cs:81.9141729700261,0.4)
--(axis cs:81.9141729700261,-0.4)
--(axis cs:67.0780217069725,-0.4)
--cycle;
\path [draw=black, fill=white]
(axis cs:61.4865788844662,0.6)
--(axis cs:61.4865788844662,1.4)
--(axis cs:82.7108163684664,1.4)
--(axis cs:82.7108163684664,0.6)
--(axis cs:61.4865788844662,0.6)
--cycle;
\addplot [black, forget plot]
table {%
67.0780217069725 0
53.6427839800571 0
};
\addplot [black, forget plot]
table {%
81.9141729700261 0
95.5117317774352 0
};
\addplot [black, forget plot]
table {%
53.6427839800571 -0.2
53.6427839800571 0.2
};
\addplot [black, forget plot]
table {%
95.5117317774352 -0.2
95.5117317774352 0.2
};
\addplot [black, forget plot]
table {%
61.4865788844662 1
50.1257959278234 1
};
\addplot [black, forget plot]
table {%
82.7108163684664 1
105.551036973033 1
};
\addplot [black, forget plot]
table {%
50.1257959278234 0.8
50.1257959278234 1.2
};
\addplot [black, forget plot]
table {%
105.551036973033 0.8
105.551036973033 1.2
};
\addplot [black, forget plot]
table {%
74.8929610122652 -0.4
74.8929610122652 0.4
};
\addplot [black, forget plot]
table {%
74.2964354788497 0.6
74.2964354788497 1.4
};
\addplot [draw=black, fill=white, forget plot, mark=*, only marks]
table{%
x  y
110.011921623282 -0.23948751357714
81.3517935059996 -0.253201664852342
89.1278106077748 0.219605509496781
79.1195629792297 -0.043618052742975
73.988116577691 0.0360648027028904
67.8542830307707 -0.147905682935119
67.3084409365509 -0.243100871812598
64.9350647299996 0.15027167013304
79.4998585142027 -0.190272446243032
95.5117317774352 0.229136662187397
84.9383034860589 -0.0869080184019093
78.4691434834197 -0.0415712412506489
74.8929610122652 -0.246838835249589
68.0714344326788 0.0658162344813291
82.4573245321324 0.0151805391362962
80.7049405072304 -0.282503099293593
78.0446262075766 -0.252105797324312
82.231384892645 -0.243385357938995
54.6486426145496 0.128003583358525
53.6427839800571 0.108122511056601
70.9310049456654 0.0247781300566324
64.3973909298601 0.0810885986115172
91.4790259922751 0.136532484656242
66.7627280938757 -0.0737333535508351
73.487343317016 -0.115180678579856
112.249823254886 -0.128712560064456
81.5969610474072 -0.184006580753954
69.5989777305069 -0.166074973679935
66.8476024773941 0.161781108232772
63.0175974335035 -0.0727962286486113
65.5634828605047 -0.152085372951882
105.551036973033 0.819383191064503
87.3884841270703 1.12412781761327
79.7876099188518 1.00011444962453
71.9462648894859 0.879390953014528
75.5254065833963 1.26613847678123
61.1292681666643 0.992732550090013
61.779092253199 1.28215839043281
61.0020147576947 0.955851439835481
82.3419514825703 0.783753367377728
90.5869093092239 0.734281917666046
79.1334758788956 0.811079882491158
83.3105025385443 0.943324875444863
67.9628163722283 0.8159785904795
74.5535568988327 0.842453674681577
83.0796812543624 0.861892495523591
68.6896485181676 1.13801913339825
73.6074321664612 0.817993047565913
74.2964354788497 0.868719676640712
51.378256049169 1.19388487618433
50.1257959278234 0.911443855004777
74.6348419544581 0.769492577856633
58.960021886072 0.8786122687111
84.221494782698 1.00515751205593
60.6578263081997 1.10408328218118
80.7887207452793 1.1444276611258
116.860479276416 0.756837469843618
90.3624106135668 1.03936048705564
63.3940325100234 1.15032882921398
60.526461844516 0.836698189306551
66.5954541807326 0.884307894820417
61.1940655157333 0.943106079593448
};
\end{axis}

\end{tikzpicture}
        \vspace{-5pt}
        \caption{LC Scenario}
        \label{fig:performance_results_a}
    \end{subfigure}
    \begin{subfigure}{\textwidth}
    \vspace{10 pt}
        \centering
        % This file was created with tikzplotlib v0.10.1.
\begin{tikzpicture}

\definecolor{darkgray176}{RGB}{176,176,176}

\begin{axis}[
height=\DelayReductionHeight,
width=\DelayReductionWidth,
legend style={fill opacity=0.8, draw opacity=1, text opacity=1, draw=none},
tick align=outside,
tick pos=left,
x grid style={darkgray176},
xlabel={\(\displaystyle \Theta[\%]\)},
xmajorgrids,
xmin=-1, xmax=101,
xtick style={color=black},
xtick={0,20,40,60,80,100},
xticklabels={
  \(\displaystyle {0}\),
  \(\displaystyle {20}\),
  \(\displaystyle {40}\),
  \(\displaystyle {60}\),
  \(\displaystyle {80}\),
  \(\displaystyle {100}\)
},
y dir=reverse,
y grid style={darkgray176},
ymin=-0.5, ymax=1.5,
ytick style={color=black},
ytick={0,1},
yticklabels={OFF,ON}
]
\path [draw=black, fill=white]
(axis cs:62.2427571347263,-0.4)
--(axis cs:62.2427571347263,0.4)
--(axis cs:69.9664992733457,0.4)
--(axis cs:69.9664992733457,-0.4)
--(axis cs:62.2427571347263,-0.4)
--cycle;
\path [draw=black, fill=white]
(axis cs:49.779896853799,0.6)
--(axis cs:49.779896853799,1.4)
--(axis cs:61.5192273694175,1.4)
--(axis cs:61.5192273694175,0.6)
--(axis cs:49.779896853799,0.6)
--cycle;
\addplot [black, forget plot]
table {%
62.2427571347263 0
56.3094423690441 0
};
\addplot [black, forget plot]
table {%
69.9664992733457 0
76.1696361216737 0
};
\addplot [black, forget plot]
table {%
56.3094423690441 -0.2
56.3094423690441 0.2
};
\addplot [black, forget plot]
table {%
76.1696361216737 -0.2
76.1696361216737 0.2
};
\addplot [black, forget plot]
table {%
49.779896853799 1
40.6719412339434 1
};
\addplot [black, forget plot]
table {%
61.5192273694175 1
70.530587412668 1
};
\addplot [black, forget plot]
table {%
40.6719412339434 0.8
40.6719412339434 1.2
};
\addplot [black, forget plot]
table {%
70.530587412668 0.8
70.530587412668 1.2
};
\addplot [black, forget plot]
table {%
67.6612801974365 -0.4
67.6612801974365 0.4
};
\addplot [black, forget plot]
table {%
57.6288007428668 0.6
57.6288007428668 1.4
};
\addplot [draw=black, fill=white, forget plot, mark=*, only marks]
table{%
x  y
56.3094423690441 -0.017263321805695
69.2728442304206 -0.172998860759978
60.4823616734851 -0.179368905914045
68.2581248763516 0.288729561470671
60.74930152853 -0.0187123779325285
72.8253883228207 -0.244438864553317
69.4596745853088 0.204671689151699
73.6026509695226 -0.186590984413287
58.4041116669603 -0.0261291232599909
64.7790761155489 0.163583721980165
65.5576328205764 0.249772463521065
59.4373857322522 -0.20490869947846
68.5148872264387 -0.0619744853758455
67.4712756430865 0.189211011698029
58.2829492564864 0.162730487628929
68.7469885656863 -0.156824327774193
76.1660000495808 0.0272454831421036
63.9692541324862 -0.0476601567749552
72.3378413260858 0.114474071982536
75.1526522813941 -0.282439146226429
64.9988239546038 -0.0310217541922664
76.1696361216737 -0.0495741685864277
63.7362127409226 0.020682088872701
68.8927588955511 0.0522252968668382
70.1983512098386 0.00196252441174594
59.6177668308448 0.20465043066546
48.0924876155559 -0.101239790503166
69.7346473368528 -0.0483819431226807
66.737303602896 -0.047775290451413
67.6612801974365 -0.0416447116379061
71.8479678912227 -0.118157784498932
40.6719412339434 0.736847881961262
54.9021113110007 1.19277999946018
46.9959893083155 1.01339009288637
60.1073282751523 0.984179529454323
49.2782996206671 1.27420274960179
70.530587412668 1.23976270641501
59.3932033961761 1.04596110998289
63.6334104627313 0.847367131993124
45.2696789351077 1.19218683041518
50.4030633995016 0.973888325092167
55.3681611872254 1.10502119342624
46.006748958026 1.25636596697597
55.9828430152788 0.702005068049107
60.6584793895397 1.17373360403072
45.8426594004774 0.720513080054781
59.9147158184651 1.11696963758708
67.3337702633902 0.974551510406028
53.9420502608636 0.990432234958518
63.2606448764328 1.09139169901506
67.6444799197027 1.08075135450013
59.3020652862078 0.815705325821767
67.6125432965706 0.993225187320823
50.281494086931 0.872813199516511
59.9157676427074 0.81795807197907
57.6288007428668 0.702022871059728
43.6762703253294 0.784485356422466
42.7722184270685 1.07854341611863
57.0710529525453 1.04411819422571
62.3799753492952 0.802933359570285
57.9955355907884 0.98595077437421
63.8001616295366 0.82502845453305
};
\end{axis}

\end{tikzpicture}
        % This file was created with tikzplotlib v0.10.1.
\begin{tikzpicture}

\definecolor{darkslategray38}{RGB}{38,38,38}
\definecolor{lightgray204}{RGB}{204,204,204}

\begin{axis}[
height=\DelayReductionHeight,
width=\DelayReductionWidth,
axis line style={darkslategray38},
legend style={fill opacity=0.8, draw opacity=1, text opacity=1, draw=none},
tick align=outside,
tick pos=left,
x grid style={lightgray204},
xlabel=\textcolor{darkslategray38}{\(\displaystyle \Theta[\mathrm{min}]\)},
xmajorgrids,
xmin=202.653457377597, xmax=583.427102132624,
xtick style={color=darkslategray38},
xtick={200,300,400,500,600},
xticklabels={
  \(\displaystyle {200}\),
  \(\displaystyle {300}\),
  \(\displaystyle {400}\),
  \(\displaystyle {500}\),
  \(\displaystyle {600}\)
},
y dir=reverse,
y grid style={lightgray204},
ymin=-0.5, ymax=1.5,
ytick style={color=darkslategray38},
ytick={0,1},
yticklabels={OFF,ON}
]
\path [draw=black, fill=white]
(axis cs:332.674312185974,-0.4)
--(axis cs:332.674312185974,0.4)
--(axis cs:399.844174597185,0.4)
--(axis cs:399.844174597185,-0.4)
--(axis cs:332.674312185974,-0.4)
--cycle;
\path [draw=black, fill=white]
(axis cs:287.328484855396,0.6)
--(axis cs:287.328484855396,1.4)
--(axis cs:334.408442455697,1.4)
--(axis cs:334.408442455697,0.6)
--(axis cs:287.328484855396,0.6)
--cycle;
\addplot [black, forget plot]
table {%
332.674312185974 0
251.523348971111 0
};
\addplot [black, forget plot]
table {%
399.844174597185 0
463.419032261057 0
};
\addplot [black, forget plot]
table {%
251.523348971111 -0.2
251.523348971111 0.2
};
\addplot [black, forget plot]
table {%
463.419032261057 -0.2
463.419032261057 0.2
};
\addplot [black, forget plot]
table {%
287.328484855396 1
219.961350321007 1
};
\addplot [black, forget plot]
table {%
334.408442455697 1
359.201265186681 1
};
\addplot [black, forget plot]
table {%
219.961350321007 0.8
219.961350321007 1.2
};
\addplot [black, forget plot]
table {%
359.201265186681 0.8
359.201265186681 1.2
};
\addplot [black, forget plot]
table {%
358.868637039548 -0.4
358.868637039548 0.4
};
\addplot [black, forget plot]
table {%
306.305563481916 0.6
306.305563481916 1.4
};
\addplot [draw=black, fill=white, forget plot, mark=*, only marks]
table{%
x  y
566.119209189214 -0.189454337642454
441.756046005305 -0.138573095916339
454.514868711695 -0.220562454502083
338.1740729844 -0.200516279901318
370.766829317541 0.199862826558792
334.797650493258 0.0991103556326704
346.045414936215 0.25128764864024
333.044118170874 -0.0786081845432483
463.419032261057 0.129415517664341
395.555954500679 0.229024410402122
369.814416699693 0.140435761116481
458.604249764889 0.0126042073750326
347.097100967308 -0.0217939076931833
359.944921443741 -0.197737093317333
435.781739704283 -0.0325365798559398
358.868637039548 0.148708076913475
353.73321580037 0.164183023967698
404.132394693691 -0.0998918939517686
251.523348971111 -0.193874296673582
302.83651176093 -0.0590987603947135
331.094157714241 -0.110374276600007
292.77320713121 0.0175264406519683
388.269221357244 0.199336317740983
330.722752926818 -0.120551038348848
349.633145223485 -0.161578671402083
554.453880808882 -0.0600693023979963
368.835057973195 0.0231829991438926
373.469490104344 0.00694204795176401
301.401662562333 0.297467947796142
332.304506201074 0.298978327082322
310.507234084474 -0.255105803374692
408.90419508413 1.0755467449479
350.113235273226 1.02672611026624
353.167027864413 0.747258187355812
297.792241668522 0.868011077633931
300.756690937984 0.813447576727719
324.247841274758 0.889532242293293
295.894644429674 1.03308981121472
287.93437185491 0.928286459917985
359.201265186681 1.15322699185117
307.77271070037 0.704452985423358
312.335015042849 0.7049442846831
354.97675966203 1.0941103767602
283.609640198209 1.07407892599343
323.600100793438 0.9625262077805
342.765665105522 0.810512297169991
312.762971193133 0.943227033983149
312.714217258721 1.00641180441294
340.784494711491 1.19767172312937
219.961350321007 1.12336150102232
272.581442130109 0.934830556478286
302.075732483484 0.856195473318405
259.882312049042 1.10027009749979
306.305563481916 1.17492572065921
287.628307186278 0.716826208090316
287.028662524513 0.734122123127069
406.195650532279 1.05631920733926
328.032390199902 1.20466165071618
305.648595926018 1.08070008867357
281.722923550346 0.898671604456434
284.833183175477 1.08646145494503
275.726820161742 0.898986685615317
};
\end{axis}

\end{tikzpicture}
       \vspace{-5pt}
        \caption{HC Scenario}
        \label{fig:performance_results_b}
    \end{subfigure}
    \caption{Distribution of the relative and absolute delay reduction $\Theta$ over all instances for a full information offline scenario (OFF) and a rolling-horizon online scenario (ON).}
    \label{fig:performance_results}
\end{figure}
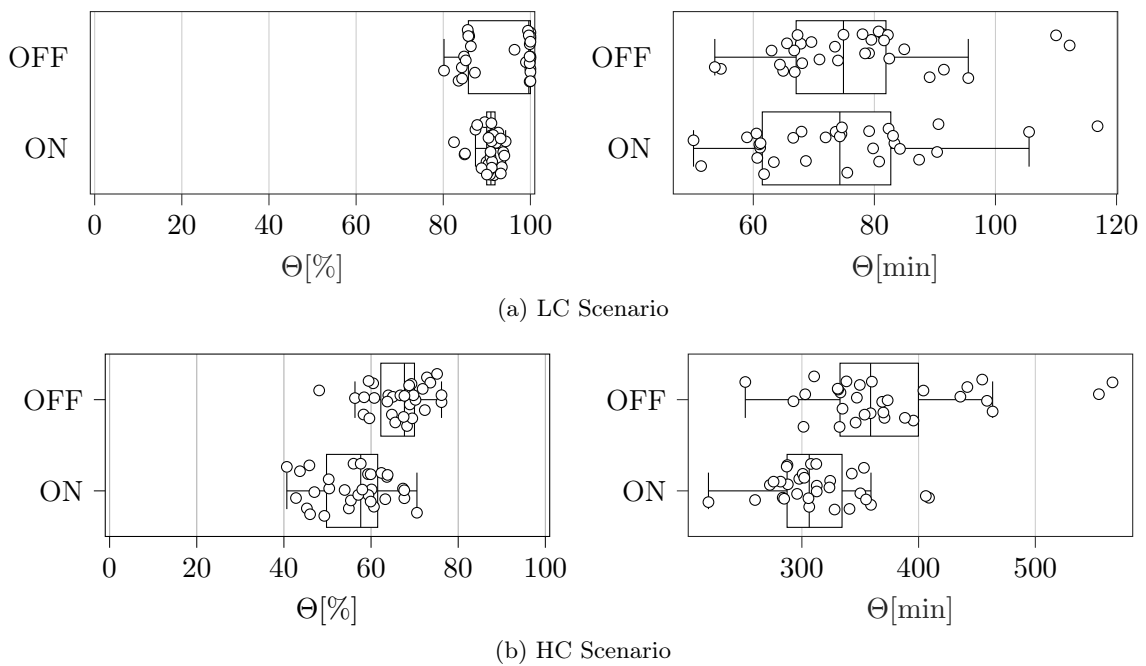

\subsection{Algorithmic performance}\label{sec:algo_performance}

Figure~\ref{fig:performance_results} shows the distribution of the delay reduction $\Theta$ obtained when using our matheuristic compared to the uncontrolled setting for both the \gls{lc} and \gls{hc} scenario. We report the delay reduction obtained when using our matheuristic in a full-information offline scenario (OFF) to obtain an upper bound on the staggered routing performance. Additionally, we report the delay reduction obtained when using our matheuristic in a rolling-horizon fashion in an online scenario (ON) to quantify the benefit of staggering in a real-world setting. For further technical discussions, we refer the readers to Appendix~\ref{sec:appendix_optimality_gaps}. We used a time limit of two hours to compute the offline solutions while limiting the computation of the rolling horizon solutions to six minutes, which equals the length of an epoch.

Focusing on the offline delay reduction, we observe a median delay reduction of $\LCOfflineMedianReduction$\% and a minimum reduction of $\LCOfflineMinReduction$\% for the \gls{lc} scenario, which shows that our matheuristic allows us to mitigate at minimum $\LCOfflineMinReduction$\% of the fleet-induced congestion in an initially uncongested system. For the \gls{hc} scenario, we observe a delay reduction within $\HCOfflineMinReduction$\% and $\HCOfflineMaxReduction$\%, with a median reduction of $\HCOfflineMedianReduction$\%, which shows that our matheuristic effectively reduces fleet-induced congestion even in a setting where the system is already congested by exogenous traffic. 

\begin{result}[offline performance]
In a full information setting, staggered routing allows for an average delay reduction of $\LCOfflineMeanReduction$\% and $\HCOfflineMeanReduction$\%, respectively, for an \gls{lc} and \gls{hc} scenario.
\end{result}

Focusing on the improvement potential in an online decision-making setting, we observe an average delay improvement of $\LCOnlineMeanReduction$\% and $\HCOnlineMeanReduction$\% for the \gls{lc} and \gls{hc} scenario when applying our matheuristic in a naive rolling horizon fashion. The observed performance difference of ten percentage points between our naive online algorithm and the full information bound remains surprisingly small in both \gls{lc} and \gls{hc} scenarios and points to the fact that even a temporally local staggering of trips allows for the reduction of a significant share of congestion-related delay. At the same time, it points to the potential of developing a more sophisticated prescriptive online algorithm that might obtain a performance close to the full information bound.

\begin{result}[online performance]
In a rolling horizon setting, staggered routing allows for an average delay reduction of $\LCOnlineMeanReduction$\% and $\HCOnlineMeanReduction$\%, respectively, for an \gls{lc} and \gls{hc} scenario.
\end{result}
\subsection{Fine-grained analysis.}\label{sec:structural_analysis}

\begin{figure}[!b]
    \centering
    \begin{minipage}{\textwidth}
        \begin{subfigure}{\textwidth}
            \centering
            \vspace{-10pt}% This file was created with tikzplotlib v0.10.1.
\begin{tikzpicture}

\definecolor{darkslategray38}{RGB}{38,38,38}
\definecolor{gray}{RGB}{128,128,128}
\definecolor{lightgray}{RGB}{211,211,211}
\definecolor{lightgray204}{RGB}{204,204,204}

\begin{axis}[
width=\columnwidth,
height=4.5cm,
axis line style={darkslategray38},
legend cell align={left},
legend style={
  fill opacity=0.8,
  draw opacity=1,
  text opacity=1,
  at={(0.867,1)},
  anchor=north west,
  draw=lightgray204
},
tick align=outside,
tick pos=left,
x grid style={lightgray204},
xlabel=\textcolor{darkslategray38}{\(\displaystyle \mathcal{E}_a\)[min]},
xmin=-2.424, xmax=32.326,
xtick style={color=darkslategray38},
xtick={0,2,4,6,8,10,12,14,16,18,20,22,24,26,28,30},
xticklabels={
  \(\displaystyle {0}\),
  \(\displaystyle {2}\),
  \(\displaystyle {4}\),
  \(\displaystyle {6}\),
  \(\displaystyle {8}\),
  \(\displaystyle {10}\),
  \(\displaystyle {12}\),
  \(\displaystyle {14}\),
  \(\displaystyle {16}\),
  \(\displaystyle {18}\),
  \(\displaystyle {20}\),
  \(\displaystyle {22}\),
  \(\displaystyle {24}\),
  \(\displaystyle {26}\),
  \(\displaystyle {28}\),
  \(\displaystyle {30}\)
},
y grid style={lightgray204},
ylabel=\textcolor{darkslategray38}{Observations},
ymajorgrids,
ymin=0.640158488872091, ymax=11693.9791163119,
yminorgrids,
ymode=log,
ytick style={color=darkslategray38},
ytick={0.01,0.1,1,10,100,1000,10000,100000,1000000},
yticklabels={
  \(\displaystyle {10^{-2}}\),
  \(\displaystyle {10^{-1}}\),
  \(\displaystyle {10^{0}}\),
  \(\displaystyle {10^{1}}\),
  \(\displaystyle {10^{2}}\),
  \(\displaystyle {10^{3}}\),
  \(\displaystyle {10^{4}}\),
  \(\displaystyle {10^{5}}\),
  \(\displaystyle {10^{6}}\)
}
]
\begin{pgfonlayer}{foreground}
  \draw[black,thick] (current axis.south west) rectangle (current axis.north east);
\end{pgfonlayer}
\draw[draw=black,fill=gray,postaction={pattern=north east lines}] (axis cs:-0.749,0.1) rectangle (axis cs:-0.249,0.1);
\addlegendimage{area legend,draw=black,fill=gray,postaction={pattern=north east lines}}
\addlegendentry{UNC}

\draw[draw=black,fill=gray,postaction={pattern=north east lines}] (axis cs:1.251,0.1) rectangle (axis cs:1.751,7486);
\draw[draw=black,fill=gray,postaction={pattern=north east lines}] (axis cs:3.251,0.1) rectangle (axis cs:3.751,100.1);
\draw[draw=black,fill=gray,postaction={pattern=north east lines}] (axis cs:5.251,0.1) rectangle (axis cs:5.751,31);
\draw[draw=black,fill=gray,postaction={pattern=north east lines}] (axis cs:7.251,0.1) rectangle (axis cs:7.751,6);
\draw[draw=black,fill=gray,postaction={pattern=north east lines}] (axis cs:9.251,0.1) rectangle (axis cs:9.751,2);
\draw[draw=black,fill=gray,postaction={pattern=north east lines}] (axis cs:11.251,0.1) rectangle (axis cs:11.751,3);
\draw[draw=black,fill=gray,postaction={pattern=north east lines}] (axis cs:13.251,0.1) rectangle (axis cs:13.751,0.1);
\draw[draw=black,fill=gray,postaction={pattern=north east lines}] (axis cs:15.251,0.1) rectangle (axis cs:15.751,1);
\draw[draw=black,fill=gray,postaction={pattern=north east lines}] (axis cs:17.251,0.1) rectangle (axis cs:17.751,0.1);
\draw[draw=black,fill=gray,postaction={pattern=north east lines}] (axis cs:19.251,0.1) rectangle (axis cs:19.751,1);
\draw[draw=black,fill=gray,postaction={pattern=north east lines}] (axis cs:21.251,0.1) rectangle (axis cs:21.751,0.1);
\draw[draw=black,fill=gray,postaction={pattern=north east lines}] (axis cs:23.251,0.1) rectangle (axis cs:23.751,0.1);
\draw[draw=black,fill=gray,postaction={pattern=north east lines}] (axis cs:25.251,0.1) rectangle (axis cs:25.751,0.1);
\draw[draw=black,fill=gray,postaction={pattern=north east lines}] (axis cs:27.251,0.1) rectangle (axis cs:27.751,0.1);
\draw[draw=black,fill=gray,postaction={pattern=north east lines}] (axis cs:29.251,0.1) rectangle (axis cs:29.751,0.1);
\draw[draw=black,fill=gray,postaction={pattern=north east lines}] (axis cs:31.251,0.1) rectangle (axis cs:31.751,0.1);
\draw[draw=black,fill=white,postaction={pattern=crosshatch dots}] (axis cs:-0.249,0.1) rectangle (axis cs:0.251,6841);
\addlegendimage{area legend,draw=black,fill=white,postaction={pattern=crosshatch dots}}
\addlegendentry{OFF}

\draw[draw=black,fill=white,postaction={pattern=crosshatch dots}] (axis cs:1.751,0.1) rectangle (axis cs:2.251,784);
\draw[draw=black,fill=white,postaction={pattern=crosshatch dots}] (axis cs:3.751,0.1) rectangle (axis cs:4.251,4);
\draw[draw=black,fill=white,postaction={pattern=crosshatch dots}] (axis cs:5.751,0.1) rectangle (axis cs:6.251,1);
\draw[draw=black,fill=white,postaction={pattern=crosshatch dots}] (axis cs:7.751,0.1) rectangle (axis cs:8.251,0.1);
\draw[draw=black,fill=white,postaction={pattern=crosshatch dots}] (axis cs:9.751,0.1) rectangle (axis cs:10.251,0.1);
\draw[draw=black,fill=white,postaction={pattern=crosshatch dots}] (axis cs:11.751,0.1) rectangle (axis cs:12.251,0.1);
\draw[draw=black,fill=white,postaction={pattern=crosshatch dots}] (axis cs:13.751,0.1) rectangle (axis cs:14.251,0.1);
\draw[draw=black,fill=white,postaction={pattern=crosshatch dots}] (axis cs:15.751,0.1) rectangle (axis cs:16.251,0.1);
\draw[draw=black,fill=white,postaction={pattern=crosshatch dots}] (axis cs:17.751,0.1) rectangle (axis cs:18.251,0.1);
\draw[draw=black,fill=white,postaction={pattern=crosshatch dots}] (axis cs:19.751,0.1) rectangle (axis cs:20.251,0.1);
\draw[draw=black,fill=white,postaction={pattern=crosshatch dots}] (axis cs:21.751,0.1) rectangle (axis cs:22.251,0.1);
\draw[draw=black,fill=white,postaction={pattern=crosshatch dots}] (axis cs:23.751,0.1) rectangle (axis cs:24.251,0.1);
\draw[draw=black,fill=white,postaction={pattern=crosshatch dots}] (axis cs:25.751,0.1) rectangle (axis cs:26.251,0.1);
\draw[draw=black,fill=white,postaction={pattern=crosshatch dots}] (axis cs:27.751,0.1) rectangle (axis cs:28.251,0.1);
\draw[draw=black,fill=white,postaction={pattern=crosshatch dots}] (axis cs:29.751,0.1) rectangle (axis cs:30.251,0.1);
\draw[draw=black,fill=white,postaction={pattern=crosshatch dots}] (axis cs:31.751,0.1) rectangle (axis cs:32.251,0.1);
\draw[draw=black,fill=lightgray,postaction={pattern=north west lines}] (axis cs:0.251,0.1) rectangle (axis cs:0.751,6483);
\addlegendimage{area legend,draw=black,fill=lightgray,postaction={pattern=north west lines}}
\addlegendentry{ON}

\draw[draw=black,fill=lightgray,postaction={pattern=north west lines}] (axis cs:2.251,0.1) rectangle (axis cs:2.751,1142);
\draw[draw=black,fill=lightgray,postaction={pattern=north west lines}] (axis cs:4.251,0.1) rectangle (axis cs:4.751,3);
\draw[draw=black,fill=lightgray,postaction={pattern=north west lines}] (axis cs:6.251,0.1) rectangle (axis cs:6.751,2);
\draw[draw=black,fill=lightgray,postaction={pattern=north west lines}] (axis cs:8.251,0.1) rectangle (axis cs:8.751,0.1);
\draw[draw=black,fill=lightgray,postaction={pattern=north west lines}] (axis cs:10.251,0.1) rectangle (axis cs:10.751,0.1);
\draw[draw=black,fill=lightgray,postaction={pattern=north west lines}] (axis cs:12.251,0.1) rectangle (axis cs:12.751,0.1);
\draw[draw=black,fill=lightgray,postaction={pattern=north west lines}] (axis cs:14.251,0.1) rectangle (axis cs:14.751,0.1);
\draw[draw=black,fill=lightgray,postaction={pattern=north west lines}] (axis cs:16.251,0.1) rectangle (axis cs:16.751,0.1);
\draw[draw=black,fill=lightgray,postaction={pattern=north west lines}] (axis cs:18.251,0.1) rectangle (axis cs:18.751,0.1);
\draw[draw=black,fill=lightgray,postaction={pattern=north west lines}] (axis cs:20.251,0.1) rectangle (axis cs:20.751,0.1);
\draw[draw=black,fill=lightgray,postaction={pattern=north west lines}] (axis cs:22.251,0.1) rectangle (axis cs:22.751,0.1);
\draw[draw=black,fill=lightgray,postaction={pattern=north west lines}] (axis cs:24.251,0.1) rectangle (axis cs:24.751,0.1);
\draw[draw=black,fill=lightgray,postaction={pattern=north west lines}] (axis cs:26.251,0.1) rectangle (axis cs:26.751,0.1);
\draw[draw=black,fill=lightgray,postaction={pattern=north west lines}] (axis cs:28.251,0.1) rectangle (axis cs:28.751,0.1);
\draw[draw=black,fill=lightgray,postaction={pattern=north west lines}] (axis cs:30.251,0.1) rectangle (axis cs:30.751,0.1);
\draw[draw=black,fill=lightgray,postaction={pattern=north west lines}] (axis cs:32.251,0.1) rectangle (axis cs:32.751,0.1);
\end{axis}

\end{tikzpicture}
            \label{fig:barplot_a}
            \vspace{-20 pt}
            \caption{LC Scenario}
        \end{subfigure}
        \begin{subfigure}{\textwidth}
        \vspace{1 pt}
            \centering
            % This file was created with tikzplotlib v0.10.1.
\begin{tikzpicture}

\definecolor{darkslategray38}{RGB}{38,38,38}
\definecolor{gray}{RGB}{128,128,128}
\definecolor{lightgray}{RGB}{211,211,211}
\definecolor{lightgray204}{RGB}{204,204,204}

\begin{axis}[
width=\columnwidth,
height=4.5cm,
axis line style={darkslategray38},
legend cell align={left},
legend style={
  fill opacity=0.8,
  draw opacity=1,
  text opacity=1,
  at={(0.867,1)},
  anchor=north west,
  draw=lightgray204
},
tick align=outside,
tick pos=left,
x grid style={lightgray204},
xlabel=\textcolor{darkslategray38}{\(\displaystyle \mathcal{E}_a\)[min]},
xmin=-2.324, xmax=32.326,
xtick style={color=darkslategray38},
xtick={0,2,4,6,8,10,12,14,16,18,20,22,24,26,28,30},
xticklabels={
  \(\displaystyle {0}\),
  \(\displaystyle {2}\),
  \(\displaystyle {4}\),
  \(\displaystyle {6}\),
  \(\displaystyle {8}\),
  \(\displaystyle {10}\),
  \(\displaystyle {12}\),
  \(\displaystyle {14}\),
  \(\displaystyle {16}\),
  \(\displaystyle {18}\),
  \(\displaystyle {20}\),
  \(\displaystyle {22}\),
  \(\displaystyle {24}\),
  \(\displaystyle {26}\),
  \(\displaystyle {28}\),
  \(\displaystyle {30}\)
},
y grid style={lightgray204},
ylabel=\textcolor{darkslategray38}{Observations},
ymajorgrids,
ymin=1.27660160426135, ymax=24859.752560285,
yminorgrids,
ymode=log,
ytick style={color=darkslategray38},
ytick={0.1,1,10,100,1000,10000,100000,1000000},
yticklabels={
  \(\displaystyle {10^{-1}}\),
  \(\displaystyle {10^{0}}\),
  \(\displaystyle {10^{1}}\),
  \(\displaystyle {10^{2}}\),
  \(\displaystyle {10^{3}}\),
  \(\displaystyle {10^{4}}\),
  \(\displaystyle {10^{5}}\),
  \(\displaystyle {10^{6}}\)
}
]
\begin{pgfonlayer}{foreground}
  \draw[black,thick] (current axis.south west) rectangle (current axis.north east);
\end{pgfonlayer}
\draw[draw=black,fill=gray,postaction={pattern=north east lines}] (axis cs:-0.749,0.1) rectangle (axis cs:-0.249,0.1);
\addlegendimage{area legend,draw=black,fill=gray,postaction={pattern=north east lines}}
\addlegendentry{UNC}

\draw[draw=black,fill=gray,postaction={pattern=north east lines}] (axis cs:1.251,0.1) rectangle (axis cs:1.751,15868);
\draw[draw=black,fill=gray,postaction={pattern=north east lines}] (axis cs:3.251,0.1) rectangle (axis cs:3.751,743);
\draw[draw=black,fill=gray,postaction={pattern=north east lines}] (axis cs:5.251,0.1) rectangle (axis cs:5.751,281);
\draw[draw=black,fill=gray,postaction={pattern=north east lines}] (axis cs:7.251,0.1) rectangle (axis cs:7.751,162);
\draw[draw=black,fill=gray,postaction={pattern=north east lines}] (axis cs:9.251,0.1) rectangle (axis cs:9.751,68);
\draw[draw=black,fill=gray,postaction={pattern=north east lines}] (axis cs:11.251,0.1) rectangle (axis cs:11.751,60.1);
\draw[draw=black,fill=gray,postaction={pattern=north east lines}] (axis cs:13.251,0.1) rectangle (axis cs:13.751,32);
\draw[draw=black,fill=gray,postaction={pattern=north east lines}] (axis cs:15.251,0.1) rectangle (axis cs:15.751,25);
\draw[draw=black,fill=gray,postaction={pattern=north east lines}] (axis cs:17.251,0.1) rectangle (axis cs:17.751,14);
\draw[draw=black,fill=gray,postaction={pattern=north east lines}] (axis cs:19.251,0.1) rectangle (axis cs:19.751,19);
\draw[draw=black,fill=gray,postaction={pattern=north east lines}] (axis cs:21.251,0.1) rectangle (axis cs:21.751,13);
\draw[draw=black,fill=gray,postaction={pattern=north east lines}] (axis cs:23.251,0.1) rectangle (axis cs:23.751,8);
\draw[draw=black,fill=gray,postaction={pattern=north east lines}] (axis cs:25.251,0.1) rectangle (axis cs:25.751,16);
\draw[draw=black,fill=gray,postaction={pattern=north east lines}] (axis cs:27.251,0.1) rectangle (axis cs:27.751,4);
\draw[draw=black,fill=gray,postaction={pattern=north east lines}] (axis cs:29.251,0.1) rectangle (axis cs:29.751,5);
\draw[draw=black,fill=white,postaction={pattern=crosshatch dots}] (axis cs:-0.249,0.1) rectangle (axis cs:0.251,10479);
\addlegendimage{area legend,draw=black,fill=white,postaction={pattern=crosshatch dots}}
\addlegendentry{OFF}

\draw[draw=black,fill=white,postaction={pattern=crosshatch dots}] (axis cs:1.751,0.1) rectangle (axis cs:2.251,6439);
\draw[draw=black,fill=white,postaction={pattern=crosshatch dots}] (axis cs:3.751,0.1) rectangle (axis cs:4.251,242);
\draw[draw=black,fill=white,postaction={pattern=crosshatch dots}] (axis cs:5.751,0.1) rectangle (axis cs:6.251,68);
\draw[draw=black,fill=white,postaction={pattern=crosshatch dots}] (axis cs:7.751,0.1) rectangle (axis cs:8.251,35);
\draw[draw=black,fill=white,postaction={pattern=crosshatch dots}] (axis cs:9.751,0.1) rectangle (axis cs:10.251,30.1);
\draw[draw=black,fill=white,postaction={pattern=crosshatch dots}] (axis cs:11.751,0.1) rectangle (axis cs:12.251,14);
\draw[draw=black,fill=white,postaction={pattern=crosshatch dots}] (axis cs:13.751,0.1) rectangle (axis cs:14.251,7);
\draw[draw=black,fill=white,postaction={pattern=crosshatch dots}] (axis cs:15.751,0.1) rectangle (axis cs:16.251,2);
\draw[draw=black,fill=white,postaction={pattern=crosshatch dots}] (axis cs:17.751,0.1) rectangle (axis cs:18.251,2);
\draw[draw=black,fill=white,postaction={pattern=crosshatch dots}] (axis cs:19.751,0.1) rectangle (axis cs:20.251,0.1);
\draw[draw=black,fill=white,postaction={pattern=crosshatch dots}] (axis cs:21.751,0.1) rectangle (axis cs:22.251,0.1);
\draw[draw=black,fill=white,postaction={pattern=crosshatch dots}] (axis cs:23.751,0.1) rectangle (axis cs:24.251,0.1);
\draw[draw=black,fill=white,postaction={pattern=crosshatch dots}] (axis cs:25.751,0.1) rectangle (axis cs:26.251,0.1);
\draw[draw=black,fill=white,postaction={pattern=crosshatch dots}] (axis cs:27.751,0.1) rectangle (axis cs:28.251,0.1);
\draw[draw=black,fill=white,postaction={pattern=crosshatch dots}] (axis cs:29.751,0.1) rectangle (axis cs:30.251,0.1);
\draw[draw=black,fill=lightgray,postaction={pattern=north west lines}] (axis cs:0.251,0.1) rectangle (axis cs:0.751,9524);
\addlegendimage{area legend,draw=black,fill=lightgray,postaction={pattern=north west lines}}
\addlegendentry{ON}

\draw[draw=black,fill=lightgray,postaction={pattern=north west lines}] (axis cs:2.251,0.1) rectangle (axis cs:2.751,7200.1);
\draw[draw=black,fill=lightgray,postaction={pattern=north west lines}] (axis cs:4.251,0.1) rectangle (axis cs:4.751,359);
\draw[draw=black,fill=lightgray,postaction={pattern=north west lines}] (axis cs:6.251,0.1) rectangle (axis cs:6.751,105);
\draw[draw=black,fill=lightgray,postaction={pattern=north west lines}] (axis cs:8.251,0.1) rectangle (axis cs:8.751,43);
\draw[draw=black,fill=lightgray,postaction={pattern=north west lines}] (axis cs:10.251,0.1) rectangle (axis cs:10.751,32);
\draw[draw=black,fill=lightgray,postaction={pattern=north west lines}] (axis cs:12.251,0.1) rectangle (axis cs:12.751,19);
\draw[draw=black,fill=lightgray,postaction={pattern=north west lines}] (axis cs:14.251,0.1) rectangle (axis cs:14.751,21);
\draw[draw=black,fill=lightgray,postaction={pattern=north west lines}] (axis cs:16.251,0.1) rectangle (axis cs:16.751,7);
\draw[draw=black,fill=lightgray,postaction={pattern=north west lines}] (axis cs:18.251,0.1) rectangle (axis cs:18.751,4);
\draw[draw=black,fill=lightgray,postaction={pattern=north west lines}] (axis cs:20.251,0.1) rectangle (axis cs:20.751,2);
\draw[draw=black,fill=lightgray,postaction={pattern=north west lines}] (axis cs:22.251,0.1) rectangle (axis cs:22.751,2);
\draw[draw=black,fill=lightgray,postaction={pattern=north west lines}] (axis cs:24.251,0.1) rectangle (axis cs:24.751,0.1);
\draw[draw=black,fill=lightgray,postaction={pattern=north west lines}] (axis cs:26.251,0.1) rectangle (axis cs:26.751,0.1);
\draw[draw=black,fill=lightgray,postaction={pattern=north west lines}] (axis cs:28.251,0.1) rectangle (axis cs:28.751,0.1);
\draw[draw=black,fill=lightgray,postaction={pattern=north west lines}] (axis cs:30.251,0.1) rectangle (axis cs:30.751,0.1);
\end{axis}

\end{tikzpicture}
            \vspace{-20 pt}
            \caption{HC Scenario}
            \label{fig:barplot_b}
        \end{subfigure}
    \end{minipage}
    \vspace{-10pt}
    \caption{Distribution of the arc-based cumulative delay ($\mathcal{E}_a$) over all instances for the uncontrolled (UNC), offline full information (OFF), and online rolling horizon (ON) settings.\label{fig:lc_barplot}}
\end{figure}
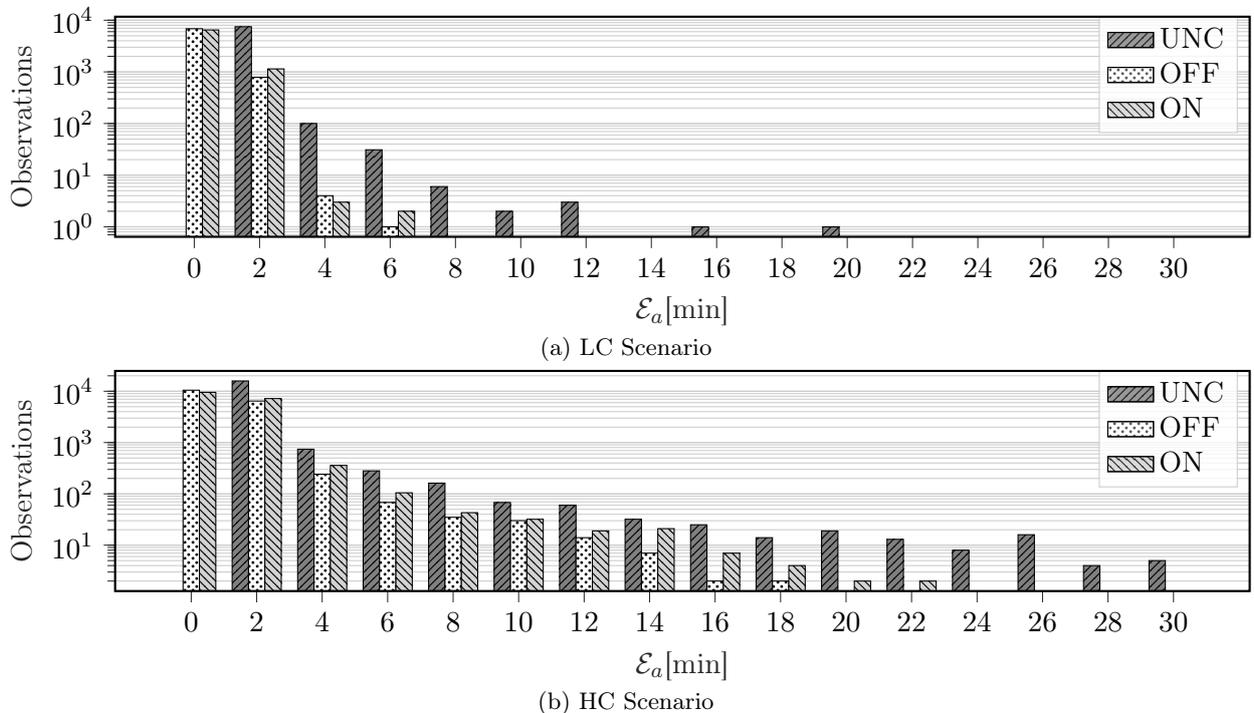

\begin{figure}[!t]
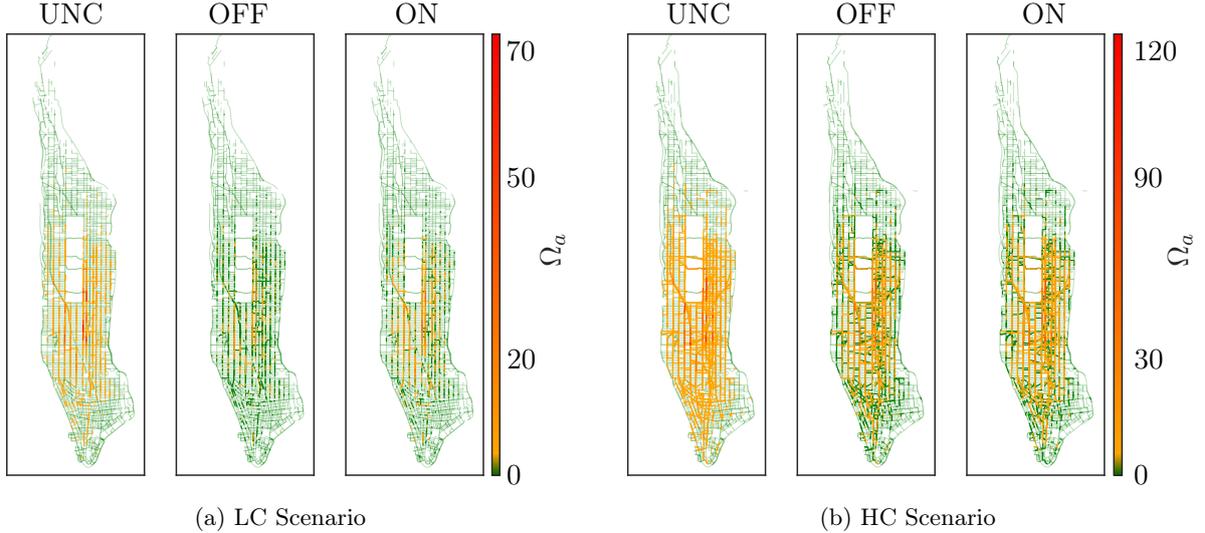

    \centering
    % First row of subfigures for LC Scenario
    \begin{subfigure}{0.5\textwidth}
        \centering
        %% Creator: Inkscape 1.1.1 (3bf5ae0d25, 2021-09-20), www.inkscape.org
%% Accompanies image file 'network_unc.pdf' (pdf, eps, ps)
\begingroup%
  \makeatletter%
  \providecommand\color[2][]{%
    \errmessage{(Inkscape) Color is used for the text in Inkscape, but the package 'color.sty' is not loaded}%
    \renewcommand\color[2][]{}%
  }%
  \providecommand\transparent[1]{%
    \errmessage{(Inkscape) Transparency is used (non-zero) for the text in Inkscape, but the package 'transparent.sty' is not loaded}%
    \renewcommand\transparent[1]{}%
  }%
  \providecommand\rotatebox[2]{#2}%
  \newcommand*\fsize{\dimexpr\f@size pt\relax}%
  \newcommand*\lineheight[1]{\fontsize{\fsize}{#1\fsize}\selectfont}%
  \ifx\svgwidth\undefined%
    \setlength{\unitlength}{90bp}%
    \ifx\svgscale\undefined%
      \relax%
    \else%
      \setlength{\unitlength}{\unitlength * \real{\svgscale}}%
    \fi%
  \else%
    \setlength{\unitlength}{\svgwidth}%
  \fi%
  \global\let\svgwidth\undefined%
  \global\let\svgscale\undefined%
  \makeatother%
  \begin{picture}(1,2.1428572)%
    \lineheight{1}%
    \setlength\tabcolsep{0pt}%
    \put(0,0){\includegraphics[width=\unitlength,page=1]{network_unc.pdf}}%
    \put(0.22166885,1.99712258){\makebox(0,0)[lt]{\lineheight{1.25}\smash{\begin{tabular}[t]{l}UNC\end{tabular}}}}%
  \end{picture}%
\endgroup%

        \hspace{-1.2cm}
        %% Creator: Inkscape 1.1.1 (3bf5ae0d25, 2021-09-20), www.inkscape.org
%% Accompanies image file 'network_off.pdf' (pdf, eps, ps)
\begingroup%
  \makeatletter%
  \providecommand\color[2][]{%
    \errmessage{(Inkscape) Color is used for the text in Inkscape, but the package 'color.sty' is not loaded}%
    \renewcommand\color[2][]{}%
  }%
  \providecommand\transparent[1]{%
    \errmessage{(Inkscape) Transparency is used (non-zero) for the text in Inkscape, but the package 'transparent.sty' is not loaded}%
    \renewcommand\transparent[1]{}%
  }%
  \providecommand\rotatebox[2]{#2}%
  \newcommand*\fsize{\dimexpr\f@size pt\relax}%
  \newcommand*\lineheight[1]{\fontsize{\fsize}{#1\fsize}\selectfont}%
  \ifx\svgwidth\undefined%
    \setlength{\unitlength}{90bp}%
    \ifx\svgscale\undefined%
      \relax%
    \else%
      \setlength{\unitlength}{\unitlength * \real{\svgscale}}%
    \fi%
  \else%
    \setlength{\unitlength}{\svgwidth}%
  \fi%
  \global\let\svgwidth\undefined%
  \global\let\svgscale\undefined%
  \makeatother%
  \begin{picture}(1,2.1428572)%
    \lineheight{1}%
    \setlength\tabcolsep{0pt}%
    \put(0,0){\includegraphics[width=\unitlength,page=1]{network_off.pdf}}%
\put(0.22166885,1.99712258){\makebox(0,0)[lt]{\lineheight{1.25}\smash{\begin{tabular}[t]{l}OFF\end{tabular}}}}%
  \end{picture}%
\endgroup%

        \hspace{-1.2cm}
        %% Creator: Inkscape 1.1.1 (3bf5ae0d25, 2021-09-20), www.inkscape.org
%% Accompanies image file 'network_on.pdf' (pdf, eps, ps)
\begingroup%
  \makeatletter%
  \providecommand\color[2][]{%
    \errmessage{(Inkscape) Color is used for the text in Inkscape, but the package 'color.sty' is not loaded}%
    \renewcommand\color[2][]{}%
  }%
  \providecommand\transparent[1]{%
    \errmessage{(Inkscape) Transparency is used (non-zero) for the text in Inkscape, but the package 'transparent.sty' is not loaded}%
    \renewcommand\transparent[1]{}%
  }%
  \providecommand\rotatebox[2]{#2}%
  \newcommand*\fsize{\dimexpr\f@size pt\relax}%
  \newcommand*\lineheight[1]{\fontsize{\fsize}{#1\fsize}\selectfont}%
  \ifx\svgwidth\undefined%
    \setlength{\unitlength}{90bp}%
    \ifx\svgscale\undefined%
      \relax%
    \else%
      \setlength{\unitlength}{\unitlength * \real{\svgscale}}%
    \fi%
  \else%
    \setlength{\unitlength}{\svgwidth}%
  \fi%
  \global\let\svgwidth\undefined%
  \global\let\svgscale\undefined%
  \makeatother%
  \begin{picture}(1,2.1428572)%
    \lineheight{1}%
    \setlength\tabcolsep{0pt}%
    \put(0,0){\includegraphics[width=\unitlength,page=1]{network_on.pdf}}%
    \put(0.97975575,0.98219445){\color[rgb]{0.14901961,0.14901961,0.14901961}\rotatebox{90}{\makebox(0,0)[lt]{\lineheight{1.25}\smash{\begin{tabular}[t]{l}$\Omega_a$\end{tabular}}}}}%
    \put(0.75563504,0.07506142){\makebox(0,0)[lt]{\lineheight{1.25}\smash{\begin{tabular}[t]{l}0\end{tabular}}}}%
    \put(0.75679757,0.55597165){\makebox(0,0)[lt]{\lineheight{1.25}\smash{\begin{tabular}[t]{l}20\end{tabular}}}}%
    %\put(0.75726261,0.93471193){\makebox(0,0)[lt]{\lineheight{1.25}\smash{\begin{tabular}[t]{l}40\end{tabular}}}}%
    \put(0.75509245,1.31453698){\makebox(0,0)[lt]{\lineheight{1.25}\smash{\begin{tabular}[t]{l}50\end{tabular}}}}%
    \put(0.75369735,1.84436223){\makebox(0,0)[lt]{\lineheight{1.25}\smash{\begin{tabular}[t]{l}70\end{tabular}}}}%
    \put(0,0){\includegraphics[width=\unitlength,page=2]{network_on.pdf}}%
    \put(0.29166885,1.99712258){\makebox(0,0)[lt]{\lineheight{1.25}\smash{\begin{tabular}[t]{l}ON\end{tabular}}}}%
  \end{picture}%
\endgroup%

        \vspace{-5pt}
        \caption{LC Scenario}
        \label{fig:heatmap_a} % Only one label is needed if the caption is for all three
    \end{subfigure}%
    \begin{subfigure}{0.5\textwidth}
        \centering
        %% Creator: Inkscape 1.1.1 (3bf5ae0d25, 2021-09-20), www.inkscape.org
%% Accompanies image file 'network_unc.pdf' (pdf, eps, ps)
\begingroup%
  \makeatletter%
  \providecommand\color[2][]{%
    \errmessage{(Inkscape) Color is used for the text in Inkscape, but the package 'color.sty' is not loaded}%
    \renewcommand\color[2][]{}%
  }%
  \providecommand\transparent[1]{%
    \errmessage{(Inkscape) Transparency is used (non-zero) for the text in Inkscape, but the package 'transparent.sty' is not loaded}%
    \renewcommand\transparent[1]{}%
  }%
  \providecommand\rotatebox[2]{#2}%
  \newcommand*\fsize{\dimexpr\f@size pt\relax}%
  \newcommand*\lineheight[1]{\fontsize{\fsize}{#1\fsize}\selectfont}%
  \ifx\svgwidth\undefined%
    \setlength{\unitlength}{90bp}%
    \ifx\svgscale\undefined%
      \relax%
    \else%
      \setlength{\unitlength}{\unitlength * \real{\svgscale}}%
    \fi%
  \else%
    \setlength{\unitlength}{\svgwidth}%
  \fi%
  \global\let\svgwidth\undefined%
  \global\let\svgscale\undefined%
  \makeatother%
  \begin{picture}(1,2.1428572)%
    \lineheight{1}%
    \setlength\tabcolsep{0pt}%
    \put(0,0){\includegraphics[width=\unitlength,page=1]{network_unc.pdf}}%
    \put(0.22166885,1.99712258){\makebox(0,0)[lt]{\lineheight{1.25}\smash{\begin{tabular}[t]{l}UNC\end{tabular}}}}%
  \end{picture}%
\endgroup%

        \hspace{-1.2cm}
        %% Creator: Inkscape 1.1.1 (3bf5ae0d25, 2021-09-20), www.inkscape.org
%% Accompanies image file 'network_off.pdf' (pdf, eps, ps)
\begingroup%
  \makeatletter%
  \providecommand\color[2][]{%
    \errmessage{(Inkscape) Color is used for the text in Inkscape, but the package 'color.sty' is not loaded}%
    \renewcommand\color[2][]{}%
  }%
  \providecommand\transparent[1]{%
    \errmessage{(Inkscape) Transparency is used (non-zero) for the text in Inkscape, but the package 'transparent.sty' is not loaded}%
    \renewcommand\transparent[1]{}%
  }%
  \providecommand\rotatebox[2]{#2}%
  \newcommand*\fsize{\dimexpr\f@size pt\relax}%
  \newcommand*\lineheight[1]{\fontsize{\fsize}{#1\fsize}\selectfont}%
  \ifx\svgwidth\undefined%
    \setlength{\unitlength}{90bp}%
    \ifx\svgscale\undefined%
      \relax%
    \else%
      \setlength{\unitlength}{\unitlength * \real{\svgscale}}%
    \fi%
  \else%
    \setlength{\unitlength}{\svgwidth}%
  \fi%
  \global\let\svgwidth\undefined%
  \global\let\svgscale\undefined%
  \makeatother%
  \begin{picture}(1,2.1428572)%
    \lineheight{1}%
    \setlength\tabcolsep{0pt}%
    \put(0,0){\includegraphics[width=\unitlength,page=1]{network_off.pdf}}%
\put(0.22166885,1.99712258){\makebox(0,0)[lt]{\lineheight{1.25}\smash{\begin{tabular}[t]{l}OFF\end{tabular}}}}%
  \end{picture}%
\endgroup%

        \hspace{-1.2cm}
        %% Creator: Inkscape 1.1.1 (3bf5ae0d25, 2021-09-20), www.inkscape.org
%% Accompanies image file 'network_on.pdf' (pdf, eps, ps)
\begingroup%
  \makeatletter%
  \providecommand\color[2][]{%
    \errmessage{(Inkscape) Color is used for the text in Inkscape, but the package 'color.sty' is not loaded}%
    \renewcommand\color[2][]{}%
  }%
  \providecommand\transparent[1]{%
    \errmessage{(Inkscape) Transparency is used (non-zero) for the text in Inkscape, but the package 'transparent.sty' is not loaded}%
    \renewcommand\transparent[1]{}%
  }%
  \providecommand\rotatebox[2]{#2}%
  \newcommand*\fsize{\dimexpr\f@size pt\relax}%
  \newcommand*\lineheight[1]{\fontsize{\fsize}{#1\fsize}\selectfont}%
  \ifx\svgwidth\undefined%
    \setlength{\unitlength}{90bp}%
    \ifx\svgscale\undefined%
      \relax%
    \else%
      \setlength{\unitlength}{\unitlength * \real{\svgscale}}%
    \fi%
  \else%
    \setlength{\unitlength}{\svgwidth}%
  \fi%
  \global\let\svgwidth\undefined%
  \global\let\svgscale\undefined%
  \makeatother%
  \begin{picture}(1,2.1428572)%
    \lineheight{1}%
    \setlength\tabcolsep{0pt}%
    \put(0,0){\includegraphics[width=\unitlength,page=1]{network_on.pdf}}%
    \put(0.97975575,0.98219445){\color[rgb]{0.14901961,0.14901961,0.14901961}\rotatebox{90}{\makebox(0,0)[lt]{\lineheight{1.25}\smash{\begin{tabular}[t]{l}$\Omega_a$\end{tabular}}}}}%
    \put(0.75563504,0.07506142){\makebox(0,0)[lt]{\lineheight{1.25}\smash{\begin{tabular}[t]{l}0\end{tabular}}}}%
    \put(0.75679757,0.55597165){\makebox(0,0)[lt]{\lineheight{1.25}\smash{\begin{tabular}[t]{l}30\end{tabular}}}}%
    %\put(0.75726261,0.93471193){\makebox(0,0)[lt]{\lineheight{1.25}\smash{\begin{tabular}[t]{l}40\end{tabular}}}}%
    \put(0.75509245,1.31453698){\makebox(0,0)[lt]{\lineheight{1.25}\smash{\begin{tabular}[t]{l}90\end{tabular}}}}%
    \put(0.75369735,1.84436223){\makebox(0,0)[lt]{\lineheight{1.25}\smash{\begin{tabular}[t]{l}120\end{tabular}}}}%
    \put(0,0){\includegraphics[width=\unitlength,page=2]{network_on.pdf}}%
    \put(0.29166885,1.99712258){\makebox(0,0)[lt]{\lineheight{1.25}\smash{\begin{tabular}[t]{l}ON\end{tabular}}}}%
  \end{picture}%
\endgroup%

        \vspace{-5pt}
        \caption{HC Scenario}
        \label{fig:heatmap_d} % Only one label is needed if the caption is for all three
    \end{subfigure}
    \caption{Spatial distribution of the number of trips per arc that experience congestion ($\Omega_a$) in the uncontrolled (UNC), offline (OFF), and online (ON) settings for both the LC and HC scenarios.}
    \label{fig:heatmap}
\end{figure}

To deepen our understanding of the impact of staggered routing, Figure~\ref{fig:lc_barplot} shows a distribution over arcs that exhibit a certain cumulative delay within the \gls{lc} or \gls{hc} scenario. We detail this distribution for the arc-based delay resulting in an uncontrolled setting (UNC) in which no staggering is applied, when solving the offline full information setting (OFF), and when solving an online rolling horizon setting (ON). Note that we exclude arcs that do not exhibit any congestion in either of these settings to focus the distributions on congestion effects. Figure~\ref{fig:heatmap} complements this analysis by showing the spatial distribution of trips incurring delay on certain arcs.

Clearly, the \gls{hc} scenario exhibits a significantly higher cumulative delay on arcs as the existing exogenous congestion in the system reinforces the congestion induced by the AMoD fleet. However, we observe in both the \gls{lc} and the \gls{hc} scenarios that staggering trip departures allow us to mitigate a significant share of the congestion induced by the AMoD fleet. Specifically, our matheuristic removes congestion completely for approximately 6500 and 10000 arcs in the \gls{lc} and \gls{hc} scenarios, independent of whether we apply it in an offline or online setting. Interestingly, the differences between the offline and online settings arise for the remaining mildly congested arcs for which the online algorithm does not succeed in anticipating temporal interdependencies.  

Figure~\ref{fig:heatmap} shows the spatial distribution of the number of trips per arc that experience delay in the respective settings. As can be seen, only a few local congestion effects remain in the \gls{lc} scenario when comparing the offline and online solutions after staggering to the uncontrolled setting. Contrarily, some main roads exhibit a decreased but still significant congestion in the \gls{hc} scenario.

\begin{result}[congestion reduction]
Staggered routing increases the number of uncongested arcs, respectively, for \gls{lc} and \gls{hc} scenarios by 6841 and 10479 when comparing the full information solutions and 6483 and 9524 when comparing the naive online solutions against the uncontrolled solutions.
\end{result}

\FloatBarrier

\subsection{Sensitivity analyses}\label{sec:sensitivity_analyses}
Clearly, the improvement potential of staggered routing depends on the amount of staggering allowed for each trip and constitutes a trade-off between passenger convenience, i.e., modifying a customer's departure time up to a certain threshold and obtaining a system optimum from a total congestion perspective. For this analysis, we focus on the full information setting for one single instance to discuss the maximum improvement potential of the respective trade-off. To analyze this trade-off, we ignore our basic setting and vary the maximum amount a trip's departure can be staggered within $\TimePercentage\in\{0\%,\ldots, 25\%\}$ of its nominal travel time with a step width of 2.5 percentage points. Here, $\TimePercentage=0\%$ equals the uncontrolled setting discussed in the preceding subsections.

\begin{figure}[!t]
    \centering
    \begin{minipage}{\textwidth}
        \centering
        \begin{subfigure}{\textwidth}
            \centering
            % This file was created with tikzplotlib v0.10.1.
\begin{tikzpicture}
\footnotesize
\definecolor{darkgray176}{RGB}{176,176,176}
\definecolor{lightgray}{RGB}{211,211,211}

\begin{axis}[
width = \TotalDelayBarplotWidth,
height = 5cm,
tick align=outside,
tick pos=left,
x grid style={darkgray176},
xlabel={\(\displaystyle \TimePercentage[\%]\)},
xmin=-2.35, xmax=27.35,
xtick style={color=black},
xtick distance=2.5,
xticklabels={-2.5,0,2.5,5,7.5,10,12.5,15,17.5,20,22.5,25},
y grid style={darkgray176},
ylabel={\(\displaystyle \TotalDelay{\Solution}\) [min]},
ymajorgrids,
ymin=0, ymax=57.3101527939095,
ytick style={color=black}
]
\draw[draw=black,fill=lightgray] (axis cs:-1,0) rectangle (axis cs:1,54.6486426145496);
\draw[draw=black,fill=lightgray] (axis cs:1.5,0) rectangle (axis cs:3.5,14.2038409736591);
\draw[draw=black,fill=lightgray] (axis cs:4,0) rectangle (axis cs:6,3.0588722208307);
\draw[draw=black,fill=lightgray] (axis cs:6.5,0) rectangle (axis cs:8.5,0.718511449075512);
\draw[draw=black,fill=lightgray] (axis cs:9,0) rectangle (axis cs:11,0);
\draw[draw=black,fill=lightgray] (axis cs:11.5,0) rectangle (axis cs:13.5,0);
\draw[draw=black,fill=lightgray] (axis cs:14,0) rectangle (axis cs:16,0);
\draw[draw=black,fill=lightgray] (axis cs:16.5,0) rectangle (axis cs:18.5,0);
\draw[draw=black,fill=lightgray] (axis cs:19,0) rectangle (axis cs:21,0);
\draw[draw=black,fill=lightgray] (axis cs:21.5,0) rectangle (axis cs:23.5,0);
\draw[draw=black,fill=lightgray] (axis cs:24,0) rectangle (axis cs:26,0);
\end{axis}

\end{tikzpicture}
            \hfill
            % This file was created with tikzplotlib v0.10.1.
\begin{tikzpicture}
\footnotesize
\definecolor{darkgray176}{RGB}{176,176,176}

\begin{axis}[
width = \ShiftAppliedWidth,
height = 5cm,
tick align=outside,
tick pos=left,
x grid style={darkgray176},
xlabel={\(\displaystyle \TimePercentage[\%]\)},
xmin=-0.5, xmax=10.5,
xtick style={color=black},
xtick={0,1,2,3,4,5,6,7,8,9,10},
xticklabels={0,2.5,5,7.5,10,12.5,15,17.5,20,22.5,25},
y grid style={darkgray176},
ylabel={\(\displaystyle \sigma^{r}\)[min]},
ymajorgrids,
ymin=-0.40208669303515, ymax=4.5,
ytick style={color=black}
]
\path [draw=black, fill=white, semithick]
(axis cs:-0.4,0)
--(axis cs:0.4,0)
--(axis cs:0.4,0)
--(axis cs:-0.4,0)
--(axis cs:-0.4,0)
--cycle;
\path [draw=black, fill=white, semithick]
(axis cs:0.6,0.0119018505426842)
--(axis cs:1.4,0.0119018505426842)
--(axis cs:1.4,0.185962426033957)
--(axis cs:0.6,0.185962426033957)
--(axis cs:0.6,0.0119018505426842)
--cycle;
\path [draw=black, fill=white, semithick]
(axis cs:1.6,0.063201442002142)
--(axis cs:2.4,0.063201442002142)
--(axis cs:2.4,0.33805983764759)
--(axis cs:1.6,0.33805983764759)
--(axis cs:1.6,0.063201442002142)
--cycle;
\path [draw=black, fill=white, semithick]
(axis cs:2.6,0.112839732038129)
--(axis cs:3.4,0.112839732038129)
--(axis cs:3.4,0.435496390174617)
--(axis cs:2.6,0.435496390174617)
--(axis cs:2.6,0.112839732038129)
--cycle;
\path [draw=black, fill=white, semithick]
(axis cs:3.6,0.129303740622193)
--(axis cs:4.4,0.129303740622193)
--(axis cs:4.4,0.506609743800724)
--(axis cs:3.6,0.506609743800724)
--(axis cs:3.6,0.129303740622193)
--cycle;
\path [draw=black, fill=white, semithick]
(axis cs:4.6,0.142418064069352)
--(axis cs:5.4,0.142418064069352)
--(axis cs:5.4,0.620827499679744)
--(axis cs:4.6,0.620827499679744)
--(axis cs:4.6,0.142418064069352)
--cycle;
\path [draw=black, fill=white, semithick]
(axis cs:5.6,0.142281048273183)
--(axis cs:6.4,0.142281048273183)
--(axis cs:6.4,0.561488572803648)
--(axis cs:5.6,0.561488572803648)
--(axis cs:5.6,0.142281048273183)
--cycle;
\path [draw=black, fill=white, semithick]
(axis cs:6.6,0.142281048273198)
--(axis cs:7.4,0.142281048273198)
--(axis cs:7.4,0.715678425000001)
--(axis cs:6.6,0.715678425000001)
--(axis cs:6.6,0.142281048273198)
--cycle;
\path [draw=black, fill=white, semithick]
(axis cs:7.6,0.144352814892651)
--(axis cs:8.4,0.144352814892651)
--(axis cs:8.4,0.703884306425313)
--(axis cs:7.6,0.703884306425313)
--(axis cs:7.6,0.144352814892651)
--cycle;
\path [draw=black, fill=white, semithick]
(axis cs:8.6,0.215607577265431)
--(axis cs:9.4,0.215607577265431)
--(axis cs:9.4,0.739902777234128)
--(axis cs:8.6,0.739902777234128)
--(axis cs:8.6,0.215607577265431)
--cycle;
\path [draw=black, fill=white, semithick]
(axis cs:9.6,0.203079095758936)
--(axis cs:10.4,0.203079095758936)
--(axis cs:10.4,0.819956211437339)
--(axis cs:9.6,0.819956211437339)
--(axis cs:9.6,0.203079095758936)
--cycle;
\addplot [semithick, black]
table {%
0 0
0 0
};
\addplot [semithick, black]
table {%
0 0
0 0
};
\addplot [semithick, black]
table {%
-0.2 0
0.2 0
};
\addplot [semithick, black]
table {%
-0.2 0
0.2 0
};
\addplot [semithick, black]
table {%
1 0.0119018505426842
1 -7.5791225147744e-15
};
\addplot [semithick, black]
table {%
1 0.185962426033957
1 0.446158697068176
};
\addplot [semithick, black]
table {%
0.8 -7.5791225147744e-15
1.2 -7.5791225147744e-15
};
\addplot [semithick, black]
table {%
0.8 0.446158697068176
1.2 0.446158697068176
};
\addplot [black, mark=x, mark size=2.5, mark options={solid}, only marks]
table {%
1 0.56803741549795
1 0.579414632147874
1 0.478936918080418
1 0.528010620157382
1 0.588342834557716
};
\addplot [semithick, black]
table {%
2 0.063201442002142
2 -9.473903143468e-16
};
\addplot [semithick, black]
table {%
2 0.33805983764759
2 0.696140572483777
};
\addplot [semithick, black]
table {%
1.8 -9.473903143468e-16
2.2 -9.473903143468e-16
};
\addplot [semithick, black]
table {%
1.8 0.696140572483777
2.2 0.696140572483777
};
\addplot [black, mark=x, mark size=2.5, mark options={solid}, only marks]
table {%
2 0.792184321834704
2 0.979150751875862
2 0.750555643348192
2 0.953745434551733
2 0.92075216049459
};
\addplot [semithick, black]
table {%
3 0.112839732038129
3 0
};
\addplot [semithick, black]
table {%
3 0.435496390174617
3 0.888070690757485
};
\addplot [semithick, black]
table {%
2.8 0
3.2 0
};
\addplot [semithick, black]
table {%
2.8 0.888070690757485
3.2 0.888070690757485
};
\addplot [black, mark=x, mark size=2.5, mark options={solid}, only marks]
table {%
3 1.56032883107058
3 0.957654097940171
3 1.19004606533945
3 0.937550205543626
3 1.10689942032078
3 1.02281479320058
};
\addplot [semithick, black]
table {%
4 0.129303740622193
4 -1.8947806286936e-15
};
\addplot [semithick, black]
table {%
4 0.506609743800724
4 1.0675948305869
};
\addplot [semithick, black]
table {%
3.8 -1.8947806286936e-15
4.2 -1.8947806286936e-15
};
\addplot [semithick, black]
table {%
3.8 1.0675948305869
4.2 1.0675948305869
};
\addplot [black, mark=x, mark size=2.5, mark options={solid}, only marks]
table {%
4 1.07662045049945
4 1.29162075270866
4 2.08043844142744
};
\addplot [semithick, black]
table {%
5 0.142418064069352
5 0
};
\addplot [semithick, black]
table {%
5 0.620827499679744
5 1.24045751731142
};
\addplot [semithick, black]
table {%
4.8 0
5.2 0
};
\addplot [semithick, black]
table {%
4.8 1.24045751731142
5.2 1.24045751731142
};
\addplot [black, mark=x, mark size=2.5, mark options={solid}, only marks]
table {%
5 1.42798702136943
5 1.50117402678664
5 2.6005480517843
};
\addplot [semithick, black]
table {%
6 0.142281048273183
6 0
};
\addplot [semithick, black]
table {%
6 0.561488572803648
6 1.16992417717729
};
\addplot [semithick, black]
table {%
5.8 0
6.2 0
};
\addplot [semithick, black]
table {%
5.8 1.16992417717729
6.2 1.16992417717729
};
\addplot [black, mark=x, mark size=2.5, mark options={solid}, only marks]
table {%
6 1.5116258231738
6 1.58746581698958
6 1.25785121711343
6 1.2383393019798
6 3.12065766214116
6 1.21016852241953
6 1.23209899999999
6 1.21055168209984
6 1.42979720324894
6 2.06948361444072
6 1.20322125
};
\addplot [semithick, black]
table {%
7 0.142281048273198
7 -3.7895612573872e-15
};
\addplot [semithick, black]
table {%
7 0.715678425000001
7 1.52664837921954
};
\addplot [semithick, black]
table {%
6.8 -3.7895612573872e-15
7.2 -3.7895612573872e-15
};
\addplot [semithick, black]
table {%
6.8 1.52664837921954
7.2 1.52664837921954
};
\addplot [black, mark=x, mark size=2.5, mark options={solid}, only marks]
table {%
7 2.3308625465704
7 1.71628553006487
7 1.66809673712376
7 2.41439755018084
7 1.69429004588154
};
\addplot [semithick, black]
table {%
8 0.144352814892651
8 -9.473903143468e-16
};
\addplot [semithick, black]
table {%
8 0.703884306425313
8 1.54196658544831
};
\addplot [semithick, black]
table {%
7.8 -9.473903143468e-16
8.2 -9.473903143468e-16
};
\addplot [semithick, black]
table {%
7.8 1.54196658544831
8.2 1.54196658544831
};
\addplot [black, mark=x, mark size=2.5, mark options={solid}, only marks]
table {%
8 1.96181733333333
8 1.90639627099858
8 1.74474100482233
8 1.604295
};
\addplot [semithick, black]
table {%
9 0.215607577265431
9 -2.65269288017104e-14
};
\addplot [semithick, black]
table {%
9 0.739902777234128
9 1.43167208364601
};
\addplot [semithick, black]
table {%
8.8 -2.65269288017104e-14
9.2 -2.65269288017104e-14
};
\addplot [semithick, black]
table {%
8.8 1.43167208364601
9.2 1.43167208364601
};
\addplot [black, mark=x, mark size=2.5, mark options={solid}, only marks]
table {%
9 1.53742605607043
9 1.68403439886366
9 2.38142189382051
9 1.81582752314976
9 2.14469580487341
9 1.96283363042512
9 1.804831875
9 1.53274801324027
};
\addplot [semithick, black]
table {%
10 0.203079095758936
10 -3.7895612573872e-15
};
\addplot [semithick, black]
table {%
10 0.819956211437339
10 1.72319069188468
};
\addplot [semithick, black]
table {%
9.8 -3.7895612573872e-15
10.2 -3.7895612573872e-15
};
\addplot [semithick, black]
table {%
9.8 1.72319069188468
10.2 1.72319069188468
};
\addplot [black, mark=x, mark size=2.5, mark options={solid}, only marks]
table {%
10 1.87372136473546
10 2.22528616037064
10 1.77959394340295
10 2.000582829612
10 1.78854388203985
10 1.77626233333332
10 2.00536875
};
\addplot [semithick, black]
table {%
-0.4 0
0.4 0
};
\addplot [semithick, black]
table {%
0.6 0.0962156430305825
1.4 0.0962156430305825
};
\addplot [semithick, black]
table {%
1.6 0.196051333333336
2.4 0.196051333333336
};
\addplot [semithick, black]
table {%
2.6 0.243402851191714
3.4 0.243402851191714
};
\addplot [semithick, black]
table {%
3.6 0.30506170000001
4.4 0.30506170000001
};
\addplot [semithick, black]
table {%
4.6 0.389867740277676
5.4 0.389867740277676
};
\addplot [semithick, black]
table {%
5.6 0.338761999999785
6.4 0.338761999999785
};
\addplot [semithick, black]
table {%
6.6 0.430875952882017
7.4 0.430875952882017
};
\addplot [semithick, black]
table {%
7.6 0.40066687340216
8.4 0.40066687340216
};
\addplot [semithick, black]
table {%
8.6 0.486385632046627
9.4 0.486385632046627
};
\addplot [semithick, black]
table {%
9.6 0.482396999999999
10.4 0.482396999999999
};
\end{axis}

\end{tikzpicture}
            \caption{LC Scenario}
            \label{fig:sensitivity_analysis_a}
        \end{subfigure}
        \begin{subfigure}{\textwidth}
        \vspace{10pt}
            \centering
            % This file was created with tikzplotlib v0.10.1.
\begin{tikzpicture}
\footnotesize
\definecolor{darkgray176}{RGB}{176,176,176}
\definecolor{lightgray}{RGB}{211,211,211}

\begin{axis}[
width = \TotalDelayBarplotWidth,
height = 5cm,
tick align= outside,
tick pos=left,
x grid style={darkgray176},
xlabel={\(\displaystyle \TimePercentage[\%]\)},
xmin=-2.35, xmax=27.35,
xtick style={color=black},
xtick distance=2.5,
xticklabels={-2.5,0,2.5,5,7.5,10,12.5,15,17.5,20,22.5,25},
y grid style={darkgray176},
ylabel={\(\displaystyle \TotalDelay{\Solution}\) [min]},
ymajorgrids,
ymin=0, ymax=511.461983713751,
ytick style={color=black}
]
\draw[draw=black,fill=lightgray] (axis cs:-1,0) rectangle (axis cs:1,347.706462288929);
\draw[draw=black,fill=lightgray] (axis cs:1.5,0) rectangle (axis cs:3.5,227.779798747265);
\draw[draw=black,fill=lightgray] (axis cs:4,0) rectangle (axis cs:6,161.438166817742);
\draw[draw=black,fill=lightgray] (axis cs:6.5,0) rectangle (axis cs:8.5,120.476884200582);
\draw[draw=black,fill=lightgray] (axis cs:9,0) rectangle (axis cs:11,96.183113317817);
\draw[draw=black,fill=lightgray] (axis cs:11.5,0) rectangle (axis cs:13.5,84.554813652697);
\draw[draw=black,fill=lightgray] (axis cs:14,0) rectangle (axis cs:16,73.4799333977004);
\draw[draw=black,fill=lightgray] (axis cs:16.5,0) rectangle (axis cs:18.5,70.5500602072657);
\draw[draw=black,fill=lightgray] (axis cs:19,0) rectangle (axis cs:21,67.3224565806686);
\draw[draw=black,fill=lightgray] (axis cs:21.5,0) rectangle (axis cs:23.5,67.0194401650671);
\draw[draw=black,fill=lightgray] (axis cs:24,0) rectangle (axis cs:26,65.7795371918359);
\end{axis}

\end{tikzpicture}
            \import{figures/hc_sensitivity_analysis/}{shift_applied}
            \caption{HC Scenario}
            \label{fig:sensitivity_analysis_b}
        \end{subfigure}
    \end{minipage}
    \caption{Trade off between realized trip departure times shifts $\StaggeringApplied{\Trip}$, the total system delay $\TotalDelay{\Solution}$, and the maximum possible shift of a trip's departure relative to its length $\TimePercentage$ for both the LC and HC scenario.
    \label{fig:sensitivity_analysis}}
\end{figure}

Figure~\ref{fig:sensitivity_analysis} highlights the resulting trade-off by showing the dependency of the total accumulated delay $\TotalDelay{\Solution}$ and the distribution of realized trip departure shifts $\StaggeringApplied{\Trip}$ in dependence of the maximum possible trip departure shift $\TimePercentage$, which indicates the maximum possible departure shift relative to each trip's length. We report results for both the \gls{lc} and \gls{hc} scenarios. As can be seen, we observe diminishing returns for delay reduction when increasing $\TimePercentage$ in both scenarios but at different scales. While it is possible to entirely mitigate the average cumulative delay in the  \gls{lc} scenario with any $\TimePercentage \geq 7.5\%$, the \gls{hc} scenario shows a plateau of an average cumulative delay of around 65 minutes for $\TimePercentage \geq 20.0\%$. Remarkably, mitigating the entire delay in the \gls{lc} scenario with maximum realized departure time shifts below $\LCMaximumStaggeringApplied$ minutes is possible when choosing $\TimePercentage = 7.5\%$. Choosing a higher $\TimePercentage$ solely leads to higher departure time shifts. In the \gls{hc} scenario, we observe increasing realized departure time shifts when increasing $\TimePercentage$ while yielding a decreasing overall delay in the system. Surprisingly, the realized departure time shifts remain mostly below $\HCMaximumStaggeringApplied$ minutes even for high $\TimePercentage$; only two outliers exhibit a shift of six minutes. 

\begin{result}[\gls{lc} scenario trade off]
    In the \gls{lc} scenario, a maximum trip departure shift of 7.5\% relative to each trip's length allows to mitigate delays induced by the AMoD fleet entirely at the price of shifting any trip departure by at most two minutes.
\end{result}

\begin{result}[\gls{hc} scenario trade off]
    In the \gls{hc} scenario, we observe a decreasing cumulative delay at the price of increasing realized trip departure times when increasing $\TimePercentage$. Still, the realized trip departure times remain, except from two outliers at most six minutes even for high $\TimePercentage$.
\end{result}

In conclusion, our results show that there is a trade-off between allowing higher maximum and realized departure time shifts and reducing delay in the system. However, a significant delay reduction can be obtained by accepting reasonable departure shifts between two and six minutes for both a system that comprises solely fleet-induced congestion (\gls{lc} scenario) and a system that is already congested and incurs additional congestion induced by the AMoD fleet (\gls{hc} scenario). 
In this context, it's important to recognize that while adjusting the departure time of a trip might not be preferable for a passenger, it nevertheless does not compromise the passenger's primary objective of reaching their destination punctually. This is because, independently of the $\TimePercentage$ parameterization, alterations to the departure time are only permissible as long as the trip finishes within the designated time window.

\FloatBarrier
\section{Conclusion} \label{sec:conclusion}
We studied the impact of staggered routing, i.e., delaying the departure of trips as long as a trip's prefixed arrival time will not be exceeded to reduce local congestion bottlenecks when operating an AMoD fleet. We formalized the underlying planning problem in a full information offline setting as a \gls{milp}. We leveraged this \gls{milp} to devise an efficient matheuristic that can be used to solve large-scale full information offline instances, while at the same time being amenable to be used in a rolling-horizon online setting. We applied this matheuristic to a real-world case study for the Manhattan area in New York City. We compared the impact of staggered routing for both a full-information bound computed in an idealized offline setting and for the solution obtained in a rolling-horizon online setting against an uncontrolled case in which staggering is not allowed. 

Our results show that in low-congestion scenarios, i.e., scenarios in which only the AMoD fleet induces congestion but the system without the fleet would remain uncongested, staggering trip departures allow to mitigate on average $\LCOfflineMeanReduction$\% of the induced congestion in a full information setting. While not reaching this full information bound, our naive rolling horizon approach still allows us to reduce $\LCOnlineMeanReduction$\% of the induced congestion in the LC scenario. In high-congestion scenarios, i.e., scenarios in which the system already encounters exogenous congestion and the AMoD fleet induces further congestion, we observe an average reduction of $\HCOfflineMeanReduction$\% as the full information bound and an average reduction of $\HCOnlineMeanReduction$\% in our online setting. Surprisingly, we show that these reductions can be reached by shifting trip departures at a maximum of $\LCMaximumStaggeringApplied$ minutes in the LC scenario and by $\HCMaximumStaggeringApplied$ minutes in the HC scenario, still preserving each trip's desired arrival time.

Our work presents a first step towards leveraging staggered routing in AMoD systems to reduce congestion. Specifically, we anticipate two fruitful avenues for future research. First, focusing on efficient solution methodologies that solve large-scale instances with high solution quality. Our paper represents a first step in this direction. Still, one may improve upon our results by finding better lower bounds and further enhancing algorithmic components tailored to the respective problem setting. Second, developing prescriptive online algorithms to apply staggered routing in practice. Our myopic rolling-horizon implementation already shows a significant improvement potential, giving hope to achieve even better performance with a more sophisticated online algorithm incorporating a learning-based prescriptive element.

%\onehalfspacing
% Acknowledgments
\noindent
{\section*{Acknowledgements}}
\noindent
The work of Antonio Coppola, Gerhard Hiermann, and Maximilian Schiffer was supported by the Deutsche Forschungsgemeinschaft (DFG) - Project number 449261765 (BalSAM). The work of Dario Paccagnan was supported by the EPSRC grant EP/Y001001/1, funded by the International Science Partnerships Fund (ISPF) and UKRI.

%
%% Finally the bib file
\singlespacing{
%\footnotesize
\bibliographystyle{model5-names}%\biboptions{authoryear}
\bibliography{references.bib}} % if more than one, comma separated
\newpage
%% and the appendices
\onehalfspacing
\begin{appendices}
	\normalsize
	\section{Hardness Proof}\label{sec:appendix_hardness}

To prove the hardness of the staggered routing problem as defined in Section~\ref{sec:problem_setting}, we utilize the hardness result for the \gls{j3uptnwt} with makespan (the maximum completion time of a job) minimization as objective, which is NP-hard \citep{SriskandarajahLadet1986}. 

\begin{theorem}\label{thr:np_hardness}
The staggered routing problem is NP-hard.
\end{theorem}

\proof{\emph{Proof}}
We begin by defining instances for \gls{j3uptnwt} and for the staggered routing problem:

\begin{definition}
Let $I$ be an instance of the \gls{j3uptnwt} with $m = 3$ machines and jobs $J_i$, with $i = 1, \ldots, n$, where each job comprises $N(J_i)$ operations that each require one unit of processing time and must be executed continuously from start to finish. The operations are allocated to specific machines, and the order in which machines are visited may vary among different jobs. Each machine can process a maximum of one job at a time. The \gls{j3uptnwt} seeks for a job schedule $S$, i.e., the starting times $s(J_i)$ of the first operation of each job $J_i$, that minimizes the schedule's makespan $D$. % NOTE: notation taken from SriskandarajahLadet1986
\end{definition}

\begin{definition}
Let $I'$ be an instance of the staggered routing problem with $n$ trips as defined in Section~\ref{sec:problem_setting}. Trips must proceed along their routes without interruption within designated time windows. The travel time for a trip along an arc is influenced by the number of trips on the arc at the specific time the trip enters the arc. If the number of trips on the arc falls below a predefined arc-specific threshold capacity, the trip proceeds at the nominal travel time designated for that arc. Conversely, exceeding this threshold triggers a delay, as detailed in Equation~\eqref{eq:travel_time_approximation}. The objective of the staggered routing problem is to find the trip departures $\TripDepartureOnArc{\Arc}{\Trip}$ that minimize the total travel time.
\end{definition}
Given these definitions, we prove the hardness of the staggered routing problem by contradiction.
\begin{assumption}\label{ass:poly_algo}
It exists an algorithm $\PolyAlgo$ that finds a feasible solution for every instance of the staggered routing problem in polynomial time.
\end{assumption}

\textit{Step 1:} Consider an instance $I'$ that comprises $n$ trips traversing $m=3$ arcs, where at least two trips have different routes. % NOTE: we need this assumption; otherwise it degenerates into a special instance of the flow-job problem, which can be solved in polynomial time \citep{SriskandarajahLadet1986}. 
The earliest possible departure time for trips is set to zero, with no restrictions on the staggering of departures.
All trips must finish by time $D$. Each arc has a nominal travel time and a threshold capacity of one; exceeding this threshold on any arc introduces an infinite delay.

\textit{Step 2:} Instance $I'$ can be transformed in polynomial time to an instance $I$ as follows: the route of each trip corresponds to a job to be executed without interruptions from start to finish. The act of traversing the arcs within a trip's route mirrors the operations within a job. These operations are executed by a dedicated machine (an arc), each with a unit processing time (equivalent to nominal travel time), and concurrent processing of operations on the same machine is prohibited (otherwise, one incurs infinite delay). 
The latest arrival times translate to the schedule's makespan under consideration. Finally, the departure times from the arcs on a route correspond to the scheduling times of the job operations.

\textit{Step 3:} Observe that selecting the job's starting times in instance $I$ to decide whether a schedule $S$ of makespan of $D$ exists is equivalent to finding feasible trip departure times in instance $I'$, such that we obtain a polynomial-time transformation from the respective staggered routing problem instance instance to an equivalent \gls{j3uptnwt} instance. 
Then, if $\PolyAlgo$ can find a feasible solution in polynomial time for $I'$, it can also find a feasible solution for $I$ within polynomial time. However, finding a solution for $I$ in polynomial time is impossible unless P = NP,  contradicting Assumption~\ref{ass:poly_algo}. Accordingly, there exists at least one instance of the staggered routing problem that cannot be solved in polynomial time.\hfill\Halmos %This contradiction effectively demonstrates the NP-hardness of the staggered routing problem, thereby completing the proof.\hfill\Halmos
\endproof

\section{Table of Variables in the MILP}\label{sec:appendix_vars_table}
\begin{table}[H]
\centering
\footnotesize
\caption{Table of Variables in the MILP.}
\label{table:basicVariables}
\begin{tabular}{rl}
\toprule
\textbf{Variable} & \textbf{Description} \\
\midrule
\(\SetTrips\) & Set of all trips $\Trip\in\SetTrips$\\
\(\SetArcs\) & Set of all arcs $\Arc\in \SetArcs$\\
\(\SetTripsOnArc{\Arc}\) & Subset of $\SetTrips$, set of trips traversing arc $\Arc$ \\
\(\TripPath{\Trip}\) & Route of trip $\Trip$ \\
\(\Successor{\Arc}{\Trip}\) & Arc succeeding $\Arc$ on route $\TripPath{\Trip}$ \\
\(\FirstArc\) & First arc of $\TripPath{\Trip}$ \\
\(\LastArc\) & Last arc of $\TripPath{\Trip}$ \\
\(\DelayOnArc{\Arc}{\Trip}\) & Delay on arc \(\Arc\) for trip \(\Trip\) \\
\(\PWLSlope{\Arc}{k}\) & Slope of the $k$-th segment of arc $\Arc$ travel time function. \\
\(\ArcNominalTravelTime{\Arc}\) & Nominal travel time for arc \(\Arc\) \\
\(\Flow{\Arc}{\Trip}\) & Number of vehicles encountered by trip \(\Trip\) on arc \(\Arc\) \\
\(\ThFlow{\Arc}{k}\) & Flow value at which the $k$-th segment of arc $\Arc$ travel time function begins. \\
\(\TripDepartureOnArc{\Arc}{\Trip}\) & Departure time of trip \(\Trip\) on arc \(\Arc\) \\
\(\TripEarliestDepartureTime{}{\Trip}\) & Earliest possible departure time for trip $\Trip$\\
\(\TripLatestExitTime{}{\Trip}\) & Latest possible arrival time for trip $\Trip$ \\
\(\MaxStaggering{\Trip}\) & Maximum staggering applicable to trip's $\Trip$ departure \\
\(\VarAlpha{\Arc}{\Trip}{\SecondTrip},\VarBeta{\Arc}{\Trip}{\SecondTrip},\VarGamma{\Arc}{\Trip}{\SecondTrip}\) & Logical variables checking concurrent presence of trips $(\Trip,\SecondTrip)$ on arc $\Arc$  \\
\(\SmallConstant\) & Small constant \\
\(M_i\) & Large constant \\
\bottomrule
\end{tabular}
\end{table}

\section{Optimality Gaps}\label{sec:appendix_optimality_gaps}

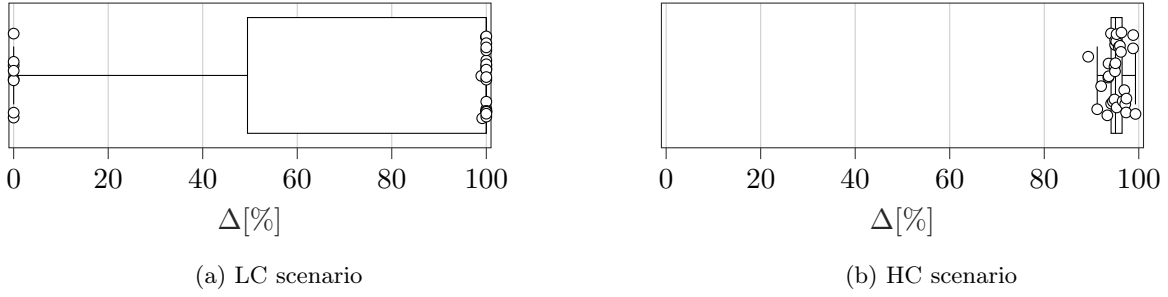
\begin{figure}[!t]
    \centering
    \begin{subfigure}{0.48\textwidth}
        % This file was created with tikzplotlib v0.10.1.
\begin{tikzpicture}

\definecolor{darkslategray38}{RGB}{38,38,38}
\definecolor{lightgray204}{RGB}{204,204,204}

\begin{axis}[
height=\GapPlotsHeight,
width=\textwidth,
yticklabels={}, % Hide y tick labels
ytick style={draw=none}, % 
axis line style={darkslategray38},
legend style={
  fill opacity=0.8,
  draw opacity=1,
  text opacity=1,
  at={(0.03,0.97)},
  anchor=north west,
  draw=none
},
tick align=outside,
tick pos=left,
x grid style={lightgray204},
xlabel=\textcolor{darkslategray38}{\(\displaystyle {\Delta}[\%]\)},
xmajorgrids,
xmin=-1, xmax=101,
xtick style={color=darkslategray38},
xtick={0,20,40,60,80,100},
xticklabels={
  \(\displaystyle {0}\),
  \(\displaystyle {20}\),
  \(\displaystyle {40}\),
  \(\displaystyle {60}\),
  \(\displaystyle {80}\),
  \(\displaystyle {100}\)
},
y dir=reverse,
y grid style={lightgray204},
ymin=-0.5, ymax=0.5,
ytick style={color=darkslategray38}
]
\path [draw=black, fill=white]
(axis cs:49.475,-0.4)
--(axis cs:49.475,0.4)
--(axis cs:100,0.4)
--(axis cs:100,-0.4)
--(axis cs:49.475,-0.4)
--cycle;
\addplot [black, forget plot]
table {%
49.475 0
0 0
};
\addplot [black, forget plot]
table {%
100 0
100 0
};
\addplot [black, forget plot]
table {%
0 -0.2
0 0.2
};
\addplot [black, forget plot]
table {%
100 -0.2
100 0.2
};
\addplot [black, forget plot]
table {%
100 -0.4
100 0.4
};
\addplot [draw=black, fill=white, forget plot, mark=*, only marks]
table{%
x  y
0 -0.0682597459167238
0 0.0270062144749302
0 0.29146482095007
0 0.0340047670023027
0 -0.0940401416583268
0 -0.289124652850632
0 0.257242172074642
0 -0.0325627036754814
98.95 0.00313646407221058
99.14 0.296060144872274
99.91 -0.2650867878665
100 -0.0501904748438652
100 0.0293666109620138
100 0.253137321943829
100 -0.070142826888813
100 0.233144868295986
100 -0.269867775090899
100 0.181424514695306
100 -0.043053557516937
100 -0.103535197704046
100 -0.173561391238547
100 0.246415559843751
100 0.245758264887113
100 0.250842543480347
100 -0.07728809844831
100 -0.222851922070893
100 -0.0401185068519828
100 -0.193764794427396
100 0.0112507340931477
100 0.284828215762876
100 0.264414378134761
};
\end{axis}

\end{tikzpicture}
        \caption{LC scenario}
        \label{fig:sub1}
    \end{subfigure}%
    \hfill
    \begin{subfigure}{0.48\textwidth}
        % This file was created with tikzplotlib v0.10.1.
\begin{tikzpicture}

\definecolor{darkslategray38}{RGB}{38,38,38}
\definecolor{lightgray204}{RGB}{204,204,204}

\begin{axis}[
height=\GapPlotsHeight,
width=\textwidth,
yticklabels={}, % Hide y tick labels
ytick style={draw=none}, % 
axis line style={darkslategray38},
legend style={fill opacity=0.8, draw opacity=1, text opacity=1, draw=none},
tick align=outside,
tick pos=left,
x grid style={lightgray204},
xlabel=\textcolor{darkslategray38}{\(\displaystyle {\Delta}[\%]\)},
xmajorgrids,
xmin=-1, xmax=101,
xtick style={color=darkslategray38},
xtick={0,20,40,60,80,100},
xticklabels={
  \(\displaystyle {0}\),
  \(\displaystyle {20}\),
  \(\displaystyle {40}\),
  \(\displaystyle {60}\),
  \(\displaystyle {80}\),
  \(\displaystyle {100}\)
},
y dir=reverse,
y grid style={lightgray204},
ymin=-0.5, ymax=0.5,
ytick style={color=darkslategray38}
]
\path [draw=black, fill=white]
(axis cs:94.14,-0.4)
--(axis cs:94.14,0.4)
--(axis cs:96.47,0.4)
--(axis cs:96.47,-0.4)
--(axis cs:94.14,-0.4)
--cycle;
\addplot [black, forget plot]
table {%
94.14 0
91.21 0
};
\addplot [black, forget plot]
table {%
96.47 0
99.3 0
};
\addplot [black, forget plot]
table {%
91.21 -0.2
91.21 0.2
};
\addplot [black, forget plot]
table {%
99.3 -0.2
99.3 0.2
};
\addplot [black, forget plot]
table {%
95.07 -0.4
95.07 0.4
};
\addplot [draw=black, fill=white, forget plot, mark=*, only marks]
table{%
x  y
89.3 -0.129370549496766
91.21 0.23355658902812
92.04 0.0740788004704418
93.35 0.275613022752581
93.56 0.0143409932076782
93.61 -0.0822968319859496
93.65 0.00558076900909121
94.12 -0.289913940067664
94.16 0.194879284274169
94.53 0.180322183586076
94.86 0.166121496094634
94.87 -0.060365782013143
94.97 -0.0284006215818802
94.99 -0.211339480649444
95.01 -0.243114191949793
95.07 -0.0843717648675194
95.24 -0.235050182669906
95.31 -0.285201830377954
95.32 0.221400748641856
95.76 -0.196459418341494
96.01 -0.20471483498712
96.23 -0.162907775131594
96.34 -0.296541815230086
96.6 0.177424166008049
96.93 0.101470201558422
97.14 0.19881399940502
97.32 0.255970943564185
97.38 0.159976821649983
98.79 -0.186913209639694
98.82 -0.278834744803288
99.3 0.266165780513443
};
\end{axis}

\end{tikzpicture}
        \caption{HC scenario}
        \label{fig:sub2}
    \end{subfigure}
    \caption{Distribution of the optimality gap $\Delta$ [\%] over all instances.}
    \label{fig:optimality_gap}
\end{figure}

Figures~\ref{fig:optimality_gap}--\ref{fig:UBLBdist} show the distribution of the optimality gap (Figure~\ref{fig:optimality_gap}), the absolute value of the lower bound (Figure~\ref{fig:LBabs}), and the distance between the upper bound and the lower bound (Figure~\ref{fig:UBLBdist}) for both the \gls{lc} scenario and the \gls{hc} scenario, obtained when running our matheuristic with a time limit of two hours.

\begin{figure}[!t]
    \centering
    \begin{subfigure}{0.48\textwidth}
        % This file was created with tikzplotlib v0.10.1.
\begin{tikzpicture}

\definecolor{darkslategray38}{RGB}{38,38,38}
\definecolor{lightgray204}{RGB}{204,204,204}

\begin{axis}[
height=\GapPlotsHeight,
width=\textwidth,
yticklabels={}, % Hide y tick labels
ytick style={draw=none}, % 
axis line style={darkslategray38},
legend style={fill opacity=0.8, draw opacity=1, text opacity=1, draw=none},
tick align=outside,
tick pos=left,
x grid style={lightgray204},
xlabel=\textcolor{darkslategray38}{LB [min]},
xmajorgrids,
xmin=-0.0119371000000477, xmax=0.250679100000002,
xtick style={color=darkslategray38},
xtick={-0.05,0,0.05,0.1,0.15,0.2,0.25,0.3},
xticklabels={
  \(\displaystyle {\ensuremath{-}0.05}\),
  \(\displaystyle {0.00}\),
  \(\displaystyle {0.05}\),
  \(\displaystyle {0.10}\),
  \(\displaystyle {0.15}\),
  \(\displaystyle {0.20}\),
  \(\displaystyle {0.25}\),
  \(\displaystyle {0.30}\)
},
y dir=reverse,
y grid style={lightgray204},
ymin=-0.5, ymax=0.5,
ytick style={color=darkslategray38}
]
\path [draw=black, fill=white]
(axis cs:0,-0.4)
--(axis cs:0,0.4)
--(axis cs:1.87900975182285e-06,0.4)
--(axis cs:1.87900975182285e-06,-0.4)
--(axis cs:0,-0.4)
--cycle;
\addplot [black, forget plot]
table {%
0 0
-4.547473508865e-14 0
};
\addplot [black, forget plot]
table {%
1.87900975182285e-06 0
3.75801948090833e-06 0
};
\addplot [black, forget plot]
table {%
-4.547473508865e-14 -0.2
-4.547473508865e-14 0.2
};
\addplot [black, forget plot]
table {%
3.75801948090833e-06 -0.2
3.75801948090833e-06 0.2
};
% \addplot [black, mark=*, mark size=2.5, mark options={solid,fill=white}, only marks, forget plot]
% table {%
% 0.11657547 0
% 0.114239538712533 0
% 0.118893749999983 0
% 0.011263713606255 0
% 0.118745053813317 0
% 0.118430235 0
% 0.238742 0
% };
\addplot [black, forget plot]
table {%
7.579122514775e-15 -0.4
7.579122514775e-15 0.4
};
\addplot [draw=black, fill=white, forget plot, mark=*, only marks]
table{%
x  y
1.89478062869333e-14 -0.193261049427295
0.11657547 0.00210945406867813
0 -0.164044620309006
0 -0.0445846257298809
2.27373675443167e-14 -0.0736920921905953
-1.13686837721617e-14 -0.260864450995199
0.114239538712533 -0.296989401725486
0.118893749999983 0.257443803207393
-7.579122514775e-15 0.19490524098909
7.579122514775e-15 0.113011075663665
3.75801948090833e-06 0.131753310505594
0.011263713606255 -0.272318789550387
0 0.152558748952627
0 -0.253179449095448
7.579122514775e-15 0.0289918283408006
0.118745053813317 -0.0322413317643125
-1.51582450295483e-14 -0.270994109932752
7.579122514775e-15 -0.0257259285427293
0 0.0840278540931328
-4.547473508865e-14 0.271019753780499
0.118430235 -0.117409863965808
0.238742 -0.190829007570314
7.579122514775e-15 -0.181497043801634
0 -0.234571082311969
1.51582450295483e-14 -0.10193433855429
1.51582450295483e-14 -0.15737165720434
-3.78956125738667e-15 0.164799684690198
0 -0.190843600342549
7.579122514775e-15 0.235779918638494
0 0.0193329826400819
7.579122514775e-15 0.181755255814825
};
\end{axis}

\end{tikzpicture}
        \caption{LC scenario}
        \label{fig:sub3}
    \end{subfigure}%
    \hfill
    \begin{subfigure}{0.48\textwidth}
        % This file was created with tikzplotlib v0.10.1.
\begin{tikzpicture}

\definecolor{darkslategray38}{RGB}{38,38,38}
\definecolor{lightgray204}{RGB}{204,204,204}

\begin{axis}[
height=\GapPlotsHeight,
width=\textwidth,
yticklabels={}, % Hide y tick labels
ytick style={draw=none}, % 
axis line style={darkslategray38},
legend style={fill opacity=0.8, draw opacity=1, text opacity=1, draw=none},
tick align=outside,
tick pos=left,
x grid style={lightgray204},
xlabel=\textcolor{darkslategray38}{LB [min]},
xmajorgrids,
xmin=1.46567467776942, xmax=26.3968648312856,
xtick style={color=darkslategray38},
xtick={0,5,10,15,20,25,30},
xticklabels={
  \(\displaystyle {0}\),
  \(\displaystyle {5}\),
  \(\displaystyle {10}\),
  \(\displaystyle {15}\),
  \(\displaystyle {20}\),
  \(\displaystyle {25}\),
  \(\displaystyle {30}\)
},
y dir=reverse,
y grid style={lightgray204},
ymin=-0.5, ymax=0.5,
ytick style={color=darkslategray38}
]
\path [draw=black, fill=white]
(axis cs:6.04517921419667,-0.4)
--(axis cs:6.04517921419667,0.4)
--(axis cs:11.1659090054408,0.4)
--(axis cs:11.1659090054408,-0.4)
--(axis cs:6.04517921419667,-0.4)
--cycle;
\addplot [black, forget plot]
table {%
6.04517921419667 0
2.59891059383833 0
};
\addplot [black, forget plot]
table {%
11.1659090054408 0
16.07943136825 0
};
\addplot [black, forget plot]
table {%
2.59891059383833 -0.2
2.59891059383833 0.2
};
\addplot [black, forget plot]
table {%
16.07943136825 -0.2
16.07943136825 0.2
};
% \addplot [black, mark=*, mark size=2.5, mark options={solid,fill=white}, only marks, forget plot]
% table {%
% 25.2636289152167 0
% };
\addplot [black, forget plot]
table {%
7.77254340247 -0.4
7.77254340247 0.4
};
\addplot [draw=black, fill=white, forget plot, mark=*, only marks]
table{%
x  y
5.17194177581833 -0.0348070535466954
5.14148937825833 -0.111222418127275
12.6060929403767 -0.120723203629399
9.250131357755 0.156807932284954
13.111336319145 -0.0736020752997513
10.9760938578767 -0.108295100144787
7.79892670547 -0.0279721539662155
6.00229285417667 0.0527837857562342
9.44937919779667 -0.0844364472354632
2.59891059383833 -0.269707845976282
11.355724153005 0.159410796105303
11.8113460891717 -0.211227660597114
8.20476506957167 0.225368912893754
8.68916397596667 0.171051113444921
14.63169683623 -0.276350845803705
12.9941555817117 0.127100671317579
4.056591228515 0.290384884821023
6.09134419298667 -0.206186764005705
6.19362271682833 -0.22578232300165
4.99345492556667 0.182820302409043
8.34083894708 -0.148764731513823
6.08806557421667 -0.116759701707412
7.51626724612667 -0.241427104886843
7.366107983965 -0.14968764411508
4.55738602371 0.283142127525273
2.622254014085 0.252467845583114
25.2636289152167 0.113754321723533
6.46520601027167 0.200460232048286
16.07943136825 0.120473067461745
7.55371060065 -0.0936883374100053
7.77254340247 0.0021688723454979
};
\end{axis}

\end{tikzpicture}
        \caption{HC scenario}
        \label{fig:sub4}
    \end{subfigure}
    \caption{Distribution of the absolute value of the lower bound LB [min] over all instances.}
    \label{fig:LBabs}
\end{figure}
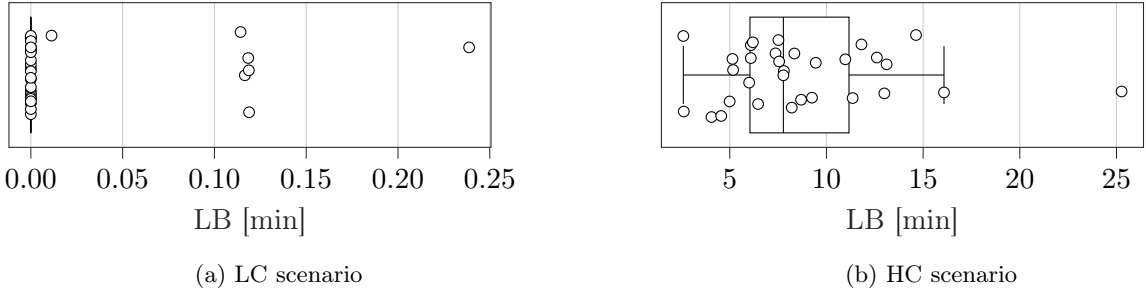
\begin{figure}[!t]
    \centering
    \begin{subfigure}{0.48\textwidth}
        % This file was created with tikzplotlib v0.10.1.
\begin{tikzpicture}

\definecolor{darkslategray38}{RGB}{38,38,38}
\definecolor{lightgray204}{RGB}{204,204,204}

\begin{axis}[
height=\GapPlotsHeight,
width=\textwidth,
yticklabels={}, % Hide y tick labels
ytick style={draw=none}, % 
axis line style={darkslategray38},
legend style={fill opacity=0.8, draw opacity=1, text opacity=1, draw=none},
tick align=outside,
tick pos=left,
x grid style={lightgray204},
xlabel=\textcolor{darkslategray38}{\(\displaystyle {\Delta}[\mathrm{min}]\)},
xmajorgrids,
xmin=-1.043, xmax=21.903,
xtick style={color=darkslategray38},
xtick={-5,0,5,10,15,20,25},
xticklabels={
  \(\displaystyle {\ensuremath{-}5}\),
  \(\displaystyle {0}\),
  \(\displaystyle {5}\),
  \(\displaystyle {10}\),
  \(\displaystyle {15}\),
  \(\displaystyle {20}\),
  \(\displaystyle {25}\)
},
y dir=reverse,
y grid style={lightgray204},
ymin=-0.5, ymax=0.5,
ytick style={color=darkslategray38}
]
\path [draw=black, fill=white]
(axis cs:0,-0.4)
--(axis cs:0,0.4)
--(axis cs:12.725,0.4)
--(axis cs:12.725,-0.4)
--(axis cs:0,-0.4)
--cycle;
\addplot [black, forget plot]
table {%
0 0
0 0
};
\addplot [black, forget plot]
table {%
12.725 0
20.86 0
};
\addplot [black, forget plot]
table {%
0 -0.2
0 0.2
};
\addplot [black, forget plot]
table {%
20.86 -0.2
20.86 0.2
};
\addplot [black, forget plot]
table {%
0.12 -0.4
0.12 0.4
};
\addplot [draw=black, fill=white, forget plot, mark=*, only marks]
table{%
x  y
0 -0.0664442390642232
0 -0.10228197569165
0 0.295127912313528
0 -0.0754775656666171
-0 0.187187507723614
0 -0.0583854517330993
0 -0.196290319283215
0 -0.176002150950401
-0 0.246453875497066
0 0.0227576018445104
0 -0.0133384458771214
0 -0.182519590059365
0 -0.0843428317563833
0 -0.286968648970379
0.11 -0.204897827723218
0.12 0.234508309839281
0.48 0.0359670787341086
2.18 -0.0870256484333988
3.61 -0.149517978019449
10.98 -0.00897894560440626
11.11 0.137116069189196
11.6 -0.150351158956611
12.43 0.270170433347882
13.02 -0.101538986753282
13.19 -0.0612197913573794
13.37 0.209893697228399
13.43 -0.0352114609577094
15.4 0.0181292362227661
17.98 0.271941130252064
20.18 0.0478045436691377
20.86 0.0210944896893826
};
\end{axis}

\end{tikzpicture}
        \caption{LC scenario}
        \label{fig:sub5}
    \end{subfigure}%
    \hfill
    \begin{subfigure}{0.48\textwidth}
        % This file was created with tikzplotlib v0.10.1.
\begin{tikzpicture}

\definecolor{darkslategray38}{RGB}{38,38,38}
\definecolor{lightgray204}{RGB}{204,204,204}

\begin{axis}[
height=\GapPlotsHeight,
width=\textwidth,
yticklabels={}, % Hide y tick labels
ytick style={draw=none}, % 
axis line style={darkslategray38},
legend style={fill opacity=0.8, draw opacity=1, text opacity=1, draw=none},
tick align=outside,
tick pos=left,
x grid style={lightgray204},
xlabel=\textcolor{darkslategray38}{\(\displaystyle {\Delta}[\mathrm{min}]\)},
xmajorgrids,
xmin=68.104, xmax=451.036,
xtick style={color=darkslategray38},
xtick={50,100,150,200,250,300,350,400,450,500},
xticklabels={
  \(\displaystyle {50}\),
  \(\displaystyle {100}\),
  \(\displaystyle {150}\),
  \(\displaystyle {200}\),
  \(\displaystyle {250}\),
  \(\displaystyle {300}\),
  \(\displaystyle {350}\),
  \(\displaystyle {400}\),
  \(\displaystyle {450}\),
  \(\displaystyle {500}\)
},
y dir=reverse,
y grid style={lightgray204},
ymin=-0.5, ymax=0.5,
ytick style={color=darkslategray38}
]
\path [draw=black, fill=white]
(axis cs:138.055,-0.4)
--(axis cs:138.055,0.4)
--(axis cs:223.99,0.4)
--(axis cs:223.99,-0.4)
--(axis cs:138.055,-0.4)
--cycle;
\addplot [black, forget plot]
table {%
138.055 0
85.51 0
};
\addplot [black, forget plot]
table {%
223.99 0
320.6 0
};
\addplot [black, forget plot]
table {%
85.51 -0.2
85.51 0.2
};
\addplot [black, forget plot]
table {%
320.6 -0.2
320.6 0.2
};
\addplot [black, forget plot]
table {%
155.62 -0.4
155.62 0.4
};
\addplot [draw=black, fill=white, forget plot, mark=*, only marks]
table{%
x  y
85.51 0.218617308272034
89.99 0.156910868409802
95.13 0.257010151014612
106.63 -0.0766784152350633
113.44 -0.287946817406041
113.89 0.205475209327093
113.95 0.138964294362786
134.14 0.16098629482183
141.97 -0.183060933600035
143.87 0.0597681957787199
144.35 0.287927212747933
148.01 0.0685400477201875
150.15 -0.139016787800254
151.27 0.170726438559788
151.3 -0.097348882865417
155.62 -0.146006057118792
164.84 -0.0163798941497648
169.95 -0.203218554024625
182.94 -0.129397478836434
190.81 -0.284065230211049
212.47 0.00748242410467215
213.4 0.234284853456625
221.54 -0.113972354668874
226.44 -0.0651379282853346
284.36 0.0734385823547131
297.29 0.148192120779628
301.16 0.232853024367807
320.6 0.175948160148437
372.71 0.254177774796753
372.83 0.243920368160611
433.63 -0.127338751175803
};
\end{axis}

\end{tikzpicture}
        \caption{HC scenario}
        \label{fig:sub6}
    \end{subfigure}
    \caption{Distribution of the distance between the upper and lower bound UB-LB [min] over all instances.}
    \label{fig:UBLBdist}
\end{figure}
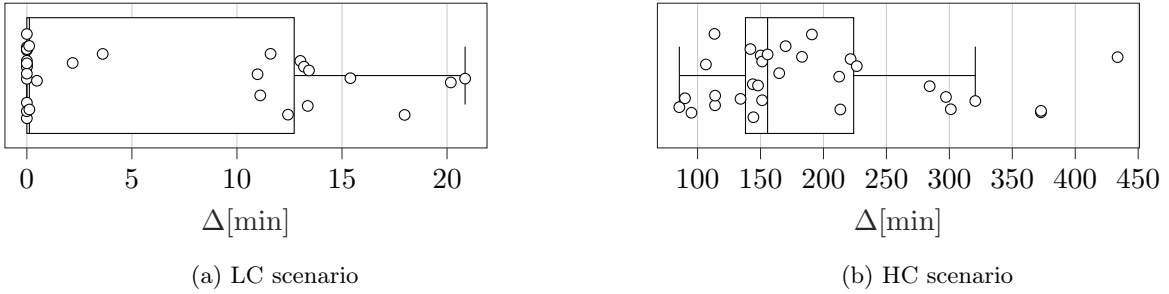

Focusing on the relative optimality gaps in Figure~\ref{fig:optimality_gap}, we observe a bipolar distribution in which optimality gaps accumulate either around 0\% or around 100\% for the \gls{lc} scenario. For the \gls{hc} scenario, all relative optimality gaps range between 90\% and 100\%. At first sight, these optimality gaps look anything but satisfying in terms of algorithmic performance, but unravel as follows at second sight: as our objective is to minimize the total congestion-related delay in the system, a trivial lower bound on our objective value is an objective value of zero. Clearly, our problem bears a huge degree of symmetry, which leads to weak lower bounds that cannot be significantly improved in our \gls{milp} formulation, see Figure~\ref{fig:LBabs}. In such cases, an absolute distance of even a few minutes of delay left in the system can create a large optimality gap. 

To support our claim that our matheuristic finds good solutions, we compare these optimality gaps to the delay reductions shown in Figure~\ref{fig:performance_results}. In fact, a comparison to the uncontrolled solution provides more insights compared to a lower bound comparison: relating Figure~\ref{fig:performance_results} to Figure~\ref{fig:optimality_gap}, it becomes obvious that the maximum improvement potential to reduce delay remains at 20\% and 50\% respectively for the \gls{lc} and \gls{hc} scenarios, although the optimality gap distribution shows values up to 100\%. 

\begin{remark}[lower bound]
    When minimizing the total congestion-related delay in a staggered routing setting, the relative optimality gap is an uninformative indicator of solution quality due to the nature of the objective function and the problem's weak lower bounds. 
\end{remark}
In this context, Figure~\ref{fig:UBLBdist} reports the difference between the solutions found and the lower bound, indicating that our matheuristic founds good solutions despite large optimality gaps. Specifically, we observe an average distance between the upper (UB) and lower (LB) bound of 6 and 193 minutes for the \gls{lc} and \gls{hc} scenarios. As the instances exhibit an average delay of 82 min and 577 min respectively in an uncontrolled setting, the absolute distance between the upper and lower bound appears rather small, indicating that our matheuristic effectively reduces congestion via staggering.

\end{appendices}
\end{document}